\documentclass[a4paper,11pt]{article}
\usepackage{amsmath,amssymb,amsthm,epsfig} 
\newcommand{\DD}{\mathbb{D}}

\newcommand{\HH}{\mathbb{H}}

\newcommand{\NN}{\mathbb{N}}

\newcommand{\PP}{\mathbb{P}}

\newcommand{\RR}{\mathbb{R}}
\renewcommand{\SS}{\mathbb{S}}

\def\cA{{\cal A}}

\def\cD{{\cal D}}   \def\cP{{\cal P}} 
  \def\cK{{\cal K}}

\renewcommand{\phi}{\varphi}
\renewcommand{\epsilon}{\varepsilon}

\newcommand{\Dom}{\mbox{Dom}}

\newtheorem{theo}{Theorem}[section]
\newtheorem{prop}[theo]{Proposition}
\newtheorem{coro}[theo]{Corollary}
\newtheorem{lemma}[theo]{Lemma}

\newtheorem{fact}[theo]{Fact}
\newtheorem*{adde}{Addendum}
\newtheorem{defi}[theo]{Definition}

\newtheorem{rema}[theo]{Remark}

\newcommand{\AdS}{\mathbb{A}d\mathbb{S}} 

\newcommand{\wt}{\widetilde}
\newcommand{\Conv}{\mbox{Conv}}
\newcommand{\Id}{\mbox{Id}}

\title{Constant  mean  curvature  foliations  of  globally  hyperbolic
  spacetimes locally modelled on $AdS_3$} 
\author{Thierry Barbot, Fran\c cois B\'eguin and Abdelghani Zeghib} 
\date{}

\begin{document}

\maketitle

\sloppy

\section{Introduction}

The purpose of this paper is  to prove the following result, which was
announced in~\cite{BarBegZeg}~:

\begin{theo}
\label{t.CMC-foliation}
Let  $(M,g)$  be   a  three-dimensional  maximal  globally  hyperbolic
spacetime, locally modelled
on  the anti-de Sitter  space $AdS_3$,  with closed  orientable Cauchy
surfaces.  Then, $M$ admits a CMC time function $\tau$.  
Moreover, the  function $\tau$ is unique and  real-analytic, and every
CMC spacelike compact surface in $M$ is a fiber of $\tau$. 
\end{theo}

%

Theorem~\ref{t.CMC-foliation} deals with three-dimensional spacetimes
whose  sectional curvature  is constant  and negative.   We  used the
equivalent formulation ``locally modelled  on the anti-de Sitter space
$AdS_3$'' to emphasize the fact that the geometry of $AdS_3$ and the
$(O(2,2),AdS_3)$-structure of  the spacetime will play  a crucial role
in    our   proof    of    the   Theorem~\ref{t.CMC-foliation}   (see
section~\ref{s.GX}).  

We  recall that  a  spacetime  $(M,g)$ is  said  to be  \emph{globally
  hyperbolic} if there exists a spacelike hypersurface $\Sigma$ in $M$ such
  that every  inextendable non-spacelike  curve intersects $\Sigma$  at one
  and  only  one   point.   Such  an  hypersurface  $\Sigma$   is  called  a
  \emph{Cauchy surface}.  A globally hyperbolic  spacetime $(M,g)$ locally
  modelled on $AdS_3$ is said to be \emph{maximal} if any embedding of
  $M$ in  a globally hyperbolic spacetime locally  modelled on $AdS_3$
  is surjective. Notice  that, if a spacetime $(M,g)$  admits a closed
  Cauchy hypersurface, then every Cauchy surface in $M$ is closed, and 
every closed spacelike hypersurface in $M$ is a Cauchy hypersurface. 
Moreover,   it   follows   from   Mess' work
  (\cite{Mes}) that a spacetime locally modeled on $AdS_3$ is maximal 
  globally hyperbolic with  compact Cauchy surfaces if and  only if it
  is  maximal with  respect to  the property  that there  is  a closed
  spacelike  surface  through  every  point.  

A \emph{time function} on a spacetime $(M,g)$ is a submersion 
$\tau:M\rightarrow\RR$ such  that $\tau$ is  strictly increasing along
every  future-directed  timelike   curve.  Every  globally  hyperbolic
spacetime  admits  (many)  time  functions.  Conversely,  a  spacetime
admitting a time  function which is surjective when  restricted to any
inextendable causal  curve is globally  hyperbolic; in this  case, the
level sets of $\tau$ are Cauchy hypersurfaces. 

A  \emph{CMC time  function} on  a spacetime  $(M,g)$ is  a  time function
$\tau:M\rightarrow\RR$  such that, for  every $\theta\in\RR$,  the set
$\tau^{-1}(\theta)$ is a spacelike hypersurface with constant mean
curvature $\theta$. In particular, a  spacetime which admits a CMC time
function  is foliated  by spacelike hypersurfaces with  constant mean
curvature. The foliation defined by a CMC time function is sometimes
called a \emph{York slicing}.  


\bigskip

Before    discussing   the implications   of
Theorem~\ref{t.CMC-foliation}, let  us say that there exist analogs of
this  theorem  for  spacetimes  with constant  non-negative  curvature
(see~\cite{And}and~\cite{Bar}  for  the flat  case  in any  dimension,
and~\cite{BarZeg}  for  the   positive  curvature  case  in  dimension
$3$).   In  fact,   three-dimensional   maximal  globally   hyperbolic
spacetimes with constant curvature  and compact Cauchy surfaces always
admits  a  CMC  time  function,  except for  three  special types  of
spacetimes~: up to finite coverings, these exceptionnal spacetimes are 
quotients of the Minkowski space $Min_3$ by a group of 
spacelike translations,  quotients of certain
domains of  the de Sitter space  $dS_3$ by rank~$2$  abelian groups of
parabolic isometries, and the de Sitter space
$dS_3$ itself.  
Even  in these special cases, there is a foliation by
compact closed CMC surfaces, which is unique except in the case of the
de Sitter space itself.  

\bigskip

The major motivation for proving Theorem~\ref{t.CMC-foliation} comes
from the  links of  this theorem with  the (vacuum)  Einstein equation.  

First  of  all, let  us  recall  that,  in dimension~$3$,  the  vacuum
Einstein  equation   (with  cosmological  constant)   reduces  to  the
requirement  that the  curvature  of the  spacetime  is constant.   In
particular, the solutions of the three-dimensional 
vacuum Einstein  with negative  cosmological constant are  exactly the
spacetimes with negative constant curvature.  

The notion of  global hyperbolicity is linked with  the most usual way
to find  solutions of the  Einstein equation: to solve  the associated
Cauchy  problem.   This  approach,  in dimension  $2+1$,  consists  in
considering a surface $\Sigma$  with a Riemannian metric $\bar{g}$ and
a symetric $2$-tensor $II$, and trying to find a Lorentzian metric $g$
on  $M=\Sigma\times ]-1,+1[$,  such  that $g$  satisfies the  Einstein
equation,  such   that  $\bar  g$   is  the  restriction  of   $g$  on
$\Sigma=\Sigma\times\{0\}$  and such that  $II$ represents  the second
fundamental  form of  $\Sigma=\Sigma\times\{0\}$  in $M=\Sigma  \times
]-1, +1[$.  For the problem to admit a solution, the initial data
$(\Sigma,\bar{g},II)$ must satisfy the \emph{constraint 
equations}  (for   geometers,  the Gauss-Codazzi  equations).  Conversely,
Choquet-Bruhat  theorem (\cite{Cho})  states that  every  initial data
satisfying the constraint  equation leads to a solution, which, by the
nature of the process, is globally hyperbolic.  Moreover, according to
Choquet-Bruhat and  Geroch (\cite{ChoGer}), there is a  unique maximal 
globally hyperbolic solution (up to isometry). 

The  main difficulty  when  dealing  with the  Cauchy  problem is  the
invariance of  Einstein equation under the  action of diffeomorphisms,
leading to an infinite dimensional space of local solutions. To bypass
this difficulty, one has to choose a gauge, \emph{i.e.} to reduce the
dimension  of   the  space   of  solution  by   imposing  additional
constraints. The method used by Choquet-Bruhat consists in considering
local coordinates $(x_1,x_2,x_3)$, such that the surface $\Sigma$
corresponds to $x_3=0$, and to  demand (with no loss of generality)
the  harmocity of  these  coordinates with  respect  to the  (unknown)
Lorentzian  metric $g$.   In such  coordinates, the  Einstein equation
becomes  a  quasi-linear   hyperbolic  equation  for  which  classical
techniques apply.  

Another similar method is to restrict to the case where each spacelike
surface  $\Sigma\times\{*\}$  is a  CMC  surface.  Then, the  equation
simplifies dramatically. The main drawback of this approach 
is that   one    has   to   assume   the   existence    of   a   CMC
surface.   Our theorem shows  that this  assumption, which  is \emph{a
  priori}  very  restrictive,   is  automatically  fulfilled  for  the
three-dimensional vacuum Einstein equation with negative cosmological
constant.   Hence,  the  remarkable  simplification  of  the  Einstein
equation  described  above, that  one  could  call ``CMC  reduction'',
applies in full generality.  

\bigskip

The CMC reduction is the  essential tool of the reduction described
by V.   Moncrief of Einstein equation to  a non-autonomous Hamiltonian
flow (that  we call \emph{Moncrief  flow}) on the cotangent  bundle of
the Teichm\"uller space of $\Sigma$ (\cite{Mon}). Moncrief flow can be
described as follows~: for every trajectory $\gamma:\RR\rightarrow
T^*\mbox{Teich}(\Sigma)$, there  exists a maximal  globally hyperbolic
space $M$  with CMC time function  $\tau$ such that  the projection of
$\gamma(t)$  on $\mbox{Teich}(\Sigma)$ is  the conformal  class $[\bar
g_t]$ of the Riemannian metric of the surface $\Sigma_t=\tau^{-1}(t)$,
and the 
cotangent vector $\gamma(t)$ is a holomorphic quadratic form extracted
from the  divergenceless and traceless part of  the second fundamental
form of $\Sigma_t$.  Our theorem shows that conversely every maximal
globally  hyperbolic spacetimes  corresponds  to a  trajectory of  the
Moncrief flow. Therefore,  maximal globally hyperbolic spacetimes with
constant  negative  curvature   and  Cauchy  surface  homeomorphic  to
$\Sigma$ are  in bijective  correspondance  with the  orbits of  the
Moncrief flow on $T^*\mbox{Teich}(\Sigma)$. 

\bigskip

Another important interest of Theorem~\ref{t.CMC-foliation}
is the uniqueness of the CMC time-function $\tau$.  In other words, 
Theorem~\ref{t.CMC-foliation}  provides a  canonical  time-function on
every  maximal globally  hyperbolic spacetime  with  constant negative
curvature and compact Cauchy surfaces.  

Note that,  we already know  another canonical time-function  on every
maximal   globally  hyperbolic   spacetimes  with   constant  negative
curvature    and    compact    Cauchy    surfaces:    the    so-called
\emph{cosmological   time   function}.     This   time   function   is
\emph{regular},  and thus,  shares nice  properties (it  is Lipschitz,
admits   first  and  second   derivatives  almost   everywhere,  etc.,
see~\cite{AndGalHow}). 
Nevertheless,   except   in  very   special   cases  (namely,   static
spacetimes), the  cosmological time function  is not differentiable everywhere,
whereas      the      CMC      time     function      provided      by
Theorem~\ref{t.CMC-foliation} is real-analytic.  

Benedetti and  Bonsante have  recently defined a  \emph{Wick rotation}
using  cosmological  time  functions  as  a key  ingredient.  In  this
context,  a Wick rotation  is a  procedure canonically  associating to
every  spacetime  locally  modelled  on $AdS_3$  a  spacetime  locally
modelled on  Minkowski space $Min_3$, or a  spacetime locally modelled
de Sitter  space $dS_3$, or a  hyperbolic manifold. One  may hope that
another Wick  rotation (the  same~?) could be  defined using  CMC time
functions.  

\bigskip

A by-product of the present article is to give new insights into the
colossal um-published work of G.\;Mess.  Indeed, a full proof of the 
classification of globally hyperbolic locally $AdS_3$ spacetimes,
with a new approach and tools, is a important step in the proof of our
principal  result. It  was  practically impossible  to  refer to  Mess
results without reproducing  ``everything''. Furthermore, we estimated
worthwhile and interesting (for the  community) to do the point on (at
least a part of) Mess work.

\subsection*{Sketch of the proof of Theorem~\ref{t.CMC-foliation}}

Consider a  maximal globally  hyperbolic $(M,g)$, locally  modelled on
$AdS_3$,    with   compact    Cauchy   surfaces.     The    proof   of
Theorem~\ref{t.CMC-foliation} 
essentially reduces to the existence  of a CMC time function $\tau$ on
$M$~: the  uniqueness of this  function follows easily from  a well-known
``maximum   principle'',   and  the
analiticity  of $\tau$  follows automatically  from  the Gauss-Codazzi
equation and from the uniqueness of the maximal solution to the Cauchy
problem for Einstein equation (see section~\ref{s.uniqueness}).  

In order to prove the existence of $\tau$, we will distinguish two
quite different cases according to whether Cauchy surfaces of $M$ have
genus~$1$  (\emph{i.e.} are two-tori),  or higher  genus (we  will see
that  a Cauchy  surface in  a locally  $AdS_3$ spacetime  cannot  be a
two-sphere).  

In the  case where $M$ admits  a Cauchy surface of  genus~$1$, we will
prove that $M$ is isometric to one of the model spacetimes known has
\emph{torus  universes} (see~\cite{Car}).  Since  such spacetimes  are
spatially homogeneous, it  is quite easy to exhibit  explicitely a CMC
time function (the level sets of the CMC times function are the orbits
of the isometry group of the  spacetime).  Note that in this case, the
CMC    time   function    co\"{\i}ncides    with   the    cosmological
time-function. This case is treated in section~\ref{s.g=1}.  

The  case of  spacetimes  with  higher genus  Cauchy  surface is  more
delicate.  We first observe that, in this case, the proof of
Theorem~\ref{t.CMC-foliation}  reduces  to  the  existence  of  a  CMC
compact spacelike surface in $M$. Indeed, using Moncrief's flow, and a
majoration of the Dirichlet  energy of CMC Cauchy surfaces, Andersson,
Moncrief  and Tromba  have proved  that the  existence of  a  CMC time
function  on $M$ follows  from the  existence of  a single  CMC Cauchy
surface in $M$ (see \cite{AndMonTro}).  

Now, a very classical and general method to prove the existence of CMC
surfaces consists in exhibiting a pair of surfaces called ``barriers''. 
In our setting, these barriers will be $C^2$ Cauchy surfaces 
$\Sigma^-,\Sigma^+$  in  $M$,  such  that the  mean
curvature of $\Sigma^+$ is everywhere negative, the mean curvature
of $\Sigma^-$ is everywhere positive, and $\Sigma^+$ is in the future of
$\Sigma^-$.   It   follows  e.g. from   a  result  of   C.  Gerhardt
(\cite{Ger}) that the existence of such barriers implies the existence
of a  Cauchy surface with  constant mean curvature (actually  a Cauchy
surface with zero mean curvature).  

So, we  are left to find  a pair of barriers  in $M$. The  way we 
construct  such  barriers  is  purely  geometrical.  One  of  the  key
ingredients of  our proof is  the locally projective structure  on the
anti de  Sitter space $AdS_3$,  which provides a notion  of convexity.
More precisely, using the time orientation and the locally 
projective  structure  of $AdS_3$,  we  will  define  some notions  of
convexity  and concavity for  spacelike surfaces  in $M$.  The key
point  is that convex  (resp. concave)  $C^2$ spacelike  surfaces have
negative (resp. positive) mean curvature. 

Mess' work implies that the spacetime $M$ can be seen as the quotient
of a domain $U$ of $AdS_3$ by a subgroup $\Gamma$ of 
$O(2,2)$.  We  give a very precise  description of the  domain $U$. In
this  description  appears  naturally  a  convex  set  $C_0$  (roughly
speaking,  $C_0$ is  the convex  hull of  the limit  set of  the group
$\Gamma$).   The boundary of  this convex  set $C_0$  is the  union of
two   disjoint  $\Gamma$-invariant   spacelike   surfaces  which   are
respectively  convex  and concave~;  the  projection $\Sigma_0^-$  and
$\Sigma_0^+$ of these surfaces in $M$ are natural candidates to be the
barriers.  

Unfortunately,  the surfaces  $\Sigma_0^-,\Sigma_0^+$  are not  smooth
(only Lipschitz).  Smoothness of  barriers is an essential requirement
in the  proof of existence of  CMC surfaces. So, the  remainder of our
proof   is   devoted   to    the   approximation   of   the   surfaces
$\Sigma_0^-,\Sigma_0^+$  by   smooth  convex  and   concave  spacelike
surfaces. Notice that this is not a so easy task as it could appear at
first glance~: standard convolution methods  can not be adapted to our
setting (see Remark~\ref{r.convolution}). 

\begin{rema}
The  notion  of \emph{convex  hypersurfaces}  can  be  defined in  any
locally  projective   space.   Hence,   the  problem  raised   by  the
non-smoothness of the surfaces $\Sigma_0^-,\Sigma_0^+$ can be seen as 
a  particular case of  a more  general question  (which, we  think, is
quite interesting)~: 
Can every (strictly) convex hypersurface in a locally
  projective space be approximated by a smooth one~?
\end{rema}

\section{Uniqueness and  analyticity of  CMC time functions} 
\label{s.uniqueness}

The purpose of this section is  to prove that, under the hypothesis of
Theorem~\ref{t.CMC-foliation},  the CMC  time function  $\tau$,  if it
exists, is unique  and real-analytic. First of all,  in order to avoid
any ambiguity on signs convention, we want to recall the definition of
the  mean  curvature  of  a  spacelike hypersurface  in  a  Lorentzian
manifold.  

\subsubsection*{Mean curvature of a spacelike hypersurface.}

Let  $\Sigma$ be a  smooth spacelike  hypersurface in  a time-oriented
Lorentzian 
manifold $M$, and  $p$ be a point of $\Sigma$.  Let  $n$ be the future
pointing unit normal vector field of $S$. We recall that the \emph{second
fundamental form} of  the surface $S$ is the  quadratic form ${II}_p$
on $T_p\Sigma$ defined by ${II}_p(X,Y)=-g(\nabla_X n,Y)$, where $g$ is
the Lorentzian  metric and $\nabla$  is the covariant  derivative. The
\emph{mean curvature}  of $S$  at $p$ is  the trace of  this quadratic
form.  

\begin{rema}
\label{r.curvature-hessian}
Let us identify the tangent space of $M$ at
$p$ with $\RR^n$, in such 
a way  that the tangent  space of $\Sigma$  at $p$ is  identified with
$\RR^{n-1}\times\{0\}$,  and   the  vector  $n$   is  identified  with
$(0,\dots,0,1)$.  Let $U$ be a neighbourhood  of $p$ in $M$. If $U$ is
small  enough, the  image  of  the surface  $\Sigma\cap  U$ under  the
inverse of the exponential map $\exp_p$ is the graph of a function
$f:\RR^{n-1}\rightarrow\RR$ such  that $f(0)=0$ and  $Df(0)=0$. 
The second fundamental form of $\Sigma$ at $p$ is
the opposite of  the hessian of $f$ at the  origin. In particular, the
mean curvature of $\Sigma$ at $p$ is the opposite of
the trace of the hessian of $f$ at the origin. 
\end{rema}

\subsubsection*{Uniqueness of the CMC time function $\tau$}

The    uniqueness   of    the  time    function   $\tau$    in
theorem~\ref{t.CMC-foliation} is a particular case of the following 
result: 

\begin{prop}
\label{p.uniqueness}
Let  $M$  be  a  globally  hyperbolic spacetime  with  compact  Cauchy
surfaces.  Assume  that $M$ admits  a CMC time function  $\tau$. Then,
every compact CMC spacelike surface in $M$ is a fiber of $\tau$.
\end{prop}
 
\begin{lemma}
\label{l.compare-curvature}
Let $\Sigma$ and  $\Sigma'$ be smooth spacelike  hypersurfaces in a
time-oriented Lorentzian manifold $M$.  Assume that $\Sigma$ and $\Sigma'$ are
tangent at  some point $p$, and  assume that $\Sigma'$  is contained in
the future of $\Sigma$.  
Then, the mean curvature of $\Sigma'$  at $p$ is smaller or equal than
those  of $\Sigma$.  Moreover,  the mean  curvatures  of $\Sigma$  and
$\Sigma'$ at  $p$ are  equal only if  $\Sigma$ and $\Sigma'$  have the
same $2$-jet at $p$.  
\end{lemma}

\begin{proof}
As  in  remark~\ref{r.curvature-hessian},  we  identify $T_p  M$  with
$\RR^n$, in such a  way that $T_p\Sigma=T_p\Sigma'$ is identified with
$\RR^{n-1}\times\{0\}$, and the future-pointing unit normal vector of
$\Sigma$ and $\Sigma'$ at  $p$ is identified with $(0,\dots,0,1)$. Let
$U$ be a neighbourhood of $p$ in $M$. If $U$ is 
small  enough, the image  of $\Sigma\cap  U$ (resp.   $\Sigma'\cap U$)
under the  inverse of  the exponential map  at $p$  is the graph  of a
function    $f:\RR^{n-1}\rightarrow\RR$   (resp.    of    a   function
$f':\RR^{n-1}\rightarrow\RR$), such that $f(0)=0$ and $Df(0)=0$ (resp.
$f'(0)=0$ and  $Df'(0)=0$). Since $\Sigma'$ is contained  in the future
of  $\Sigma$,  we  have  $f'\geq  f$. This  implies  that,  for  every
$v\in\RR^{n-1}$, we have $D^2f'(0).(v,v)\geq D^2f(0)(v,v)$.  According
to  Remark~\ref{r.curvature-hessian},  this   implies  that  the  mean
curvature  of $\Sigma'$  at  $p$ is  smaller  or equal  than those  of
$\Sigma'$.

The case  of equality is  a consequence of the  following observation:
given two functions   $f,f':\RR^{n-1}\rightarrow\RR$   such   that
$f(0)=f'(0)=0$  and $Df(0)=Df'(0)=0$,  and such  that $f'\geq  f$, then
the hessians of $f$ and $f'$ at $p$ are equal if and only if they have
the same trace. 
\end{proof}

\begin{proof}[Proof of the Proposition~\ref{p.uniqueness}]
For every $s\in\tau(\RR)$, denote by $\Sigma_s$ the Cauchy surface
$\tau^{-1}(s)$. Recall  that, for every  $s$, $\Sigma_s$ is  a compact
Cauchy surface  with constant  mean curvature equal  to $s$.  Now, let
$s_1:=\inf\{s\in\RR \mid 
\Sigma\cap\Sigma_s\neq\emptyset\}$ and $s_2:=\inf\{s\in\RR \mid 
\Sigma\cap\Sigma_s\neq\emptyset\}$. The  compactness of $\Sigma$ implies
that $s_1$  and $s_2$ do  exist (\emph{i.e.} are in  $\tau(\RR)$), and
that   $\Sigma$  does  intersect   the  surfaces   $\Sigma_{s_1}$  and
$\Sigma_{s_2}$. Moreover, 
by definition of $s_1$ and  $s_2$, the surface $\Sigma$ is contained in
the future the  surface $\Sigma_{s_1}$ and in the  past of the surface
$\Sigma_{s_2}$.      Let      $p_1$       be      a      point      in
$\Sigma\cap\Sigma_{s_1}$,     and    $p_2$     be    a     point    in
$\Sigma\cap\Sigma_{s_2}$. By Lemma~\ref{l.compare-curvature},  the
mean curvature of $\Sigma$ at $p_1$ is at most $s_1$, and the mean
curvature of $\Sigma$ at $p_2$ is  at least $s_2$. Since $\Sigma$ is a
CMC surface, and since $s_1\leq s_2$, this implies $s_1=s_2$. Moreover,
since $\Sigma$ is in the future of $\Sigma_{s_1}$ and in the past of
$\Sigma_{s_2}$,                      this                      implies
$\Sigma=\Sigma_{s_1}=\Sigma_{s_2}$. 
\end{proof}

\begin{rema}
\label{revet finis}
The uniqueness of CMC time function, when it exists, implies that it is 
preserved by isometries; in particular, by covering automorphisms
of isometric coverings. Hence, if a given spacetime admits a CMC time
function, the same is true for all its  finite quotients. This remark
enables us, for the proof of Theorem \ref{t.CMC-foliation}, 
to replace at every moment
the spacetime under consideration by any finite covering.
\end{rema}

\subsubsection*{Analiticity of the CMC time function $\tau$}

At  first  glance, uniqueness  of  CMC  foliations  suggests an  extra
regularity  of them.  However,  uniqueness seems  to come  from global
reasons, and  so only an automatic continuity  (i.e. $C^0$ regularity)
is  guaranteed by general  principles. One  knows, for  instance, many
situations  in  mathematics  (e.g.   dynamical systems  theory)  where
canonical objects  are defined by  an infinite limit process,  and are
therefore never smooth.   The situation is better here!   The point is
that,  due  to  the  formalism  of the  Cauchy  problem  for  Einstein
equations, one can have a double  vision. The first one is a spacetime
endowed with  a (local) CMC foliation.  The second one is  a CMC data,
that  is,   a  Riemannian  manifold  satisfying   a  ``CMC  constraint
equation'', which generates a spacetime having this manifold as a leaf
of a CMC foliation. The  regularity of the foliation derives thus from
that of the associated PDE system. 
More formally:

\begin{prop} 
\label{p.analytic}
Let $(M, g)$ be an {\em analytic} Lorentz manifold satisfying
vacuum Einstein equation with negative cosmological constant, that is
$Ricci_g= \Lambda g$ with $\Lambda<0$. Let $N \subset M$ be a compact
(spacelike) CMC hypersurface.  Then, there 
is a unique CMC foliation extending $N$, defined on a neighbourhood of 
it. This foliation is furthermore analytic.

In  particular,  any  (locally  defined) CMC  foliation  with  compact
leaves is analytic.  
\end{prop}

\begin{proof} 
Firstly,  a  folkloric  fact  on  Riemannian geometry  says  that  CMC
hypersurfaces in  analytic manifolds are analytic. The  reason is that
they  solve a  quasi-linear elliptic  PDE  of degree  2. This  extends
to the Lorentz case.  

Now, consider $N$, a CMC hypersurface in $M$, and let $h$ and $k$ be its
restricted   (Riemannian)   metric   and   second   fundamental   form
respectively.   Then,  $(N,h,k)$  is  a  CMC vacuum  data.   See,  for
instance \cite{AndMon}, for a  modern exposition on Einstein equations
in CMC gauges.   The authors write Einstein equation  in a gauge which
is harmonic on space, and CMC on time.  They show that the obtained
hyperbolic-elliptic PDE system is well-posed. In particular, solutions 
are analytic provided that initial data are. 
\end{proof}


\section{A short presentation of $(G,X)$-structures}
\label{s.GX}

Let $X$ be a manifold and $G$ be a group acting on $X$ with
  the following property:  if an element $g$ of  $G$ acts trivially on
  an open subset of $X$, then $g$ is the identity element of $G$.
A   \emph{$(G,X)$-structure}   on  a   manifold   $M$   is  an   atlas
  $(U_i,\phi_i)_{i\in I}$ where:\\
-- $(U_i)_{i\in I}$ is a covering of $M$ by open subsets, \\
-- for every $i$, the map $\phi_i$ is a homeomorphism from $U_i$ to an open
  subset of $X$,\\
-- for        every       $i,j$,       the        transition       map
  $\phi_i\circ\phi_j^{-1}:\phi_j(U_i\cap U_j)\rightarrow\phi_i(U_i\cap
  U_j)$ is the restriction of an element of $G$. 

To every manifold $M$ equipped with a $(G,X)$-structure are associated
two  natural  objects:   the  \emph{developping  map}  $\cD:\widetilde
M\rightarrow  X$, which is  a local  homeomorphism from  the universal
covering $\widetilde  M$ of $M$  to some open  subset of $X$,  and the
\emph{holonomy representation} $\rho:\pi_1(M)\rightarrow G$.
These natural objects satisfy the following equivariance 
property: for every  $x\in\widetilde M$ and every $\gamma\in\pi_1(M)$,
one has $\cD(\gamma.x)=\rho(\gamma).\cD(x)$. \\ 
A good reference for all these notions is~\cite{Gol2}. 

In   this  article,  we   are  interested   in  spacetimes   that  are
\emph{locally modelled} on the  anti-de Sitter space $AdS_3$, that is,
manifolds  equipped  with   a  $(G,X)$-structure  with  $X=AdS_3$  and
$G=\mbox{Isom}_0(AdS_3)=O_0(2,2)$.

\section{The three dimensional anti-de Sitter space}
\label{s.AdS}

In this section, we recall the construction of the different models of
the three-dimensional anti-de Sitter space, and we study the geometrical
properties of this space. 

\subsection{The linear model of the anti-de Sitter space}
\label{ss.linear-model}

We   denote  by  $(x_1,x_2,x_3,x_4)$   the  standard   coordinates  on
$\RR^4$.      We      will      also     use      the      coordinates
$(a,b,c,d)=(x_1-x_3,-x_2+x_4,x_2+x_4,x_1+x_3)$.    We   consider   the
quadratic form $Q=-x_1^2-x_2^2+x_3^2+x_4^2=-ad+bc$ and denote
by $B_Q$ the bilinear form associated to~$Q$.  

Let  $p$  be   a  point  on  the  quadric   of  equation  $(Q=-1)$  in
$\RR^4$. When  we identify  the tangent space  of $\RR^4$ at  $p$ with
$\RR^4$,  the  tangent  space  of  the  quadric  $(Q=-1)$  at  $p$  is
identified   with  the  $Q$-orthogonal   of  $p$.   Since  $Q$   is  a
non-degenerate  quadratic  form of  signature  $(-,-,+,+)$, and  since
$Q(p)=-1$, the  restriction of $Q$ to  the $Q$-orthogonal of  $p$ is a
non-degenerate quadratic form of signature $(-,+,+)$.  
This proves that the quadratic form $Q$ induces a Lorentzian metric of
signature  $(-,+,+)$ on  the  quadric $(Q=-1)$.  In  other words,  the
restriction        of        the       pseudo-Riemannian        metric
\hbox{$-dx_1^2-dx_2^2+dx_3^2+dx_4^2$}  to the  quadric  $(Q=-1)$ is  a
Lorentzian metric of signature $(-,+,+)$. 

\begin{defi}
The \emph{(linear model of the) three-dimensional anti-de Sitter space},
denoted by  $AdS_3$, is the  quadric $(Q=-1)$ in $\RR^4$ endowed with
the Lorentzian metric induced by $Q$.  
\end{defi}

One  can  easily verify  that  the  anti-de  Sitter space  $AdS_3$  is
diffeomorphic to  $\SS^1\times\RR^2$.  More precisely, one  can find a
diffeomorphism   $h:\SS^1\times\RR^2\rightarrow   AdS_3$   such   that
the  surface  $h(\{\theta\}\times\RR^2)$  is spacelike  for  every
$\theta$, and such that the circle $h(\SS^1\times\{x\})$ is timelike
for  every $x$.  In particular,  the anti-de  Sitter space  $AdS_3$ is
time-orientable; from now on,  we will assume that a time-orientation
has been chosen. 

The isometry  group of the anti-de  Sitter space $AdS_3$  is the group
$O(2,2)$ of  the linear transformations of $\RR^4$  which preserve the
quadratic form  $Q$. The group  $O(2,2)$ acts transitively  on $AdS_3$
and the stabilizer of  any point is isomorphic to  $O(2,1)$; hence, the
anti-de Sitter space $AdS_3$ can be seen as the homogenous space
$O(2,2)/O(2,1)$. We shall denote by $O_0(2,2)$ the connected component
of the identity  of $O(2,2)$; the elements of  $O_0(2,2)$ preserve the
three-dimensional orientation and the time-orientation of $AdS_3$.

\begin{prop}
\label{p.geodesics}
The  geodesics  of  $AdS_3$   are  the  connected  components  of  the
intersections of $AdS_3$ with  the two-dimensional vector subspaces of
$\RR^4$.
\end{prop}

\begin{proof}  
Let $P$ be a two-dimensional vector subspace of $\RR^4$. The geometry
of $P\cap AdS_3$ depends on the signature of the restriction of $Q$ to
the plane $P$:\\
-- If the restriction of $Q$ to the plane $P$ is a quadratic
form of  signature~$(-,-)$, then there  exists an element  $\sigma$ of
$O(2,2)$  which   maps  $P$  to  the   plane  $(x_3=0,x_4=0)$.   The
intersection  of $AdS_3$ with  the plane  $(x_3=0,x_4=0)$ is  a closed
timelike curve.  This curve has
to be a  geodesic of $AdS_3$, since  it is the fixed point  set of the
symetry  with  respect  to  the  plane $(x_3=0,x_4=0)$,  which  is  an
isometry of $AdS_3$.  Hence, the 
intersection of $AdS_3$  with the plane $P$ is  also a closed timelike
geodesic of $AdS_3$. \\
-- If the restriction of $Q$ to the plane $P$ is a
quadratic form of signature $(-,+)$, then there exists an element of
$O(2,2)$  which  maps  $P$  to  the plane  $(x_1=0,x_3=0)$.  The  same
arguments as above imply that
$P\cap  AdS_3$ is  the  union of  two  disjoint non-closed spacelike
geodesics of $AdS_3$.\\
-- If the restriction of $Q$ to the plane $P$ is a
degenerate quadratic form of signature $(0,-)$, then there exists an
element   of   $O(2,2)$   which    maps  $P$   to   the   plane
$(x_1=x_3,x_4=0)$. The same arguments as  in the first case imply that
$P\cap AdS_3$ is a non-closed lightlike geodesic 
of $AdS_3$. \\
-- Finally, if the restriction  of $Q$ to the plane $P$ is a
quadratic form of signature $(+,+)$,  $(0,-)$ or $(0,0)$, then one can
easily verify that the intersection $P\cap AdS_3$ is empty.  

The  discussion above  implies that  each connected  component  of the
intersection  of  $AdS_3$  with  a 2-dimensional  vector  subspace  of
$\RR^4$ is a  geodesic of $AdS_3$.  The converse  follows from the fact
that a geodesic  is uniquely determined by its  tangent vector at some
point. 
\end{proof}

\begin{rema}
\label{r.type-geodesics}
Let  $\gamma$  be a  geodesic  of  $AdS_3$.  According to  Proposition
\ref{p.geodesics},  there  exists  a $2$-dimensional  vector  subspace
$P_\gamma$ of $\RR^4$ such that $\gamma$ is a connected component of
$P_\gamma\cap   AdS_3$.  Moreover, reading  again  the   proof  of
Proposition~\ref{p.geodesics}, we notice that:

\noindent  --  if  $\gamma$  is  timelike, then  the  intersection  of
$P_\gamma$ with the quadric $(Q=0)$ is reduced to $(0,0,0,0)$;

\noindent -- if  $\gamma$ is lightlike, then $P_\gamma$  is tangent to
the quadric $(Q=0)$ along a line;

\noindent -- if $\gamma$ is spacelike, then $P_\gamma$ intersects 
transversally the quadric $(Q=0)$ along two lines. 
\end{rema}

\begin{rema}
\label{r.causal-structure-AdS}
The proof of Proposition \ref{p.geodesics} shows that all the timelike
geodesics of  $AdS_3$ are  closed, so  that a single  point is  not an
``achronal''  set  in $AdS_3$.  Moreover,  one  can  prove that  the  past
and the future in $AdS_3$ of  any point $p\in AdS_3$ are both equal to
the whole of $AdS_3$.  So, the causal structure of $AdS_3$ is not very
interesting. This  is the  reason why, instead  of working  in $AdS_3$
itself, we  shall work in some  ``large'' subsets of  $AdS_3$ which do
not contain any closed geodesics (see subsection \ref{ss.affine-domains}).  
\end{rema}

Using the same kind of arguments as in the proof of Proposition
\ref{p.geodesics}, one can prove the following:

\begin{prop}
\label{p.totally-geodesic}
The  two-dimensional  totally  geodesic   subspaces  of  $AdS_3$  are  the
connected  components  of  the   intersections  of  $AdS_3$  with  the
three-dimensional vector subspaces of $\RR^4$. 
\end{prop}

\begin{rema}
In particular, given  any point $p\in AdS_3$ and any vector plane $P$
in $T_p AdS_3$, there exists a totally geodesic subspace of $AdS_3$
whose tangent space at $p$ is the plane $P$.
\end{rema}

Let $p$ be a point in $AdS_3$. We call \emph{dual surface of the point
  $p$} the intersection $p^*$ of  the   hyperplane
$p^\perp=\{q\in\RR^4\mid   B_Q(p,q)=0\}$  with   $AdS_3$;   hence,  by
Proposition~\ref{p.totally-geodesic},  each   connected  component  of
$p^*$   is   a    two-dimensional   totally   geodesic   subspace   of
$AdS_3$. One  can easily  verify that $p^*$  is made of  two connected
components, and that the restriction 
  of $Q$  to $p^*$  is a  quadratic form of  signature $(+,+)$  (it is
  enough to consider the case where $p$ is the point $(1,0,0,0)$ since
  $O_0(2,2)$ acts transitively on  $AdS_3$). Hence, the surface $p^*$ is
  the union  of two disjoint spacelike totally  geodesic subspaces of
  $AdS_3$.  

\begin{rema}
\label{r.from-p-to-p*}
Every point  of the surface $p^*$ can  be joined from $p$  by a
timelike geodesic segment.  
\end{rema}

\begin{proof}
Let $q$ be a point in  $p^*$. We denote by $P$ the $2$-dimensional
vector  subspace   spanned  by  $p$   and  $q$  in  $\RR^4$.   We  have
$Q(p)=Q(q)=-1$ and $B_Q(p,q)=0$; this  implies that the restriction of
the  quadratic form  $Q$  to the  plane  $P$ is  a  quadratic form  of
signature    $(-,-)$.    Hence,    according   to    the    proof   of
Proposition~\ref{p.geodesics}, the intersection  of the plane $P$ with
$AdS_3$ is  a timelike  geodesic. This proves  in particular  that the
points $p$ and $q$ are joined by a timelike geodesic segment.  
\end{proof}

\subsection{The Klein model of the anti-de Sitter space}
\label{ss.projective-model}

We  shall now  define  the  ``Klein model  of  the anti-de  Sitter
space". An interesting feature of this model is that it allows us to attach 
a boundary to the anti-de Sitter space. This boundary
will   play  a   fundamental  role   in  the   proof   of  Theorem
\ref{t.CMC-foliation}.  

We see the sphere $\SS^3$  as the quotient of $\RR^4\setminus\{0\}$ by
positive homotheties.   We denote by  $\pi$ the natural  projection of
$\RR^4\setminus\{0\}$  on $\SS^3$.   We denote  by $[x_1:x_2:x_3:x_4]$
the ``positively  homogenous'' coordinates  on $\SS^3$ induced  by the
coordinates     $(x_1,x_2,x_3,x_4)$     on     $\RR^4$:    one     has
$[x_1:x_2:x_3:x_4]=[y_1:y_2:y_3:y_4]$  if  and  only if  there  exists
$\lambda>0$         such        that        $(x_1,x_2,x_3,x_4)=\lambda
(y_1,y_2,y_3,y_4)$. Similarly, we denote by $[a:b:c:d]$ the positively
homogenous   coordinates  on  $\SS^3$   induced  by   the  coordinates
$(a,b,c,d)$  on   $\RR^4$.   We  endow  $\SS^3$   with  its  canonical
Riemannian metric.  

\begin{rema}
\label{r.sign-well-defined}
Given  a point  $p\in\SS^3$,  the  quantity $Q(p)$  is  defined up  to
multiplication  by a  positive number;  this  means that  the sign  of
$Q(p)$ is well-defined. Similarly, given two points $p,q\in\SS^3$, the
sign of $B_Q(p,q)$ is well-defined.  
\end{rema}

\begin{defi}
The projection $\pi$ maps diffeomorphically $AdS_3$ on its image 
$\pi(AdS_3)\subset\SS^3$.  The  \emph{Klein model of the anti-de
  Sitter space}, that  we denote by $\AdS_3$, is  the image of $AdS_3$
under  $\pi$, equipped  with the  image  of the  Lorentzian metric  of
$AdS_3$.  We  denote by $\partial\AdS_3$  the boundary of  $\AdS_3$ in
$\SS^3$.  
\end{defi}

Observe that $\AdS_3$  is made of the points  of $\SS^3$ which satisfy
the  inequation $(Q<0)$.   Hence, $\partial\AdS_3$  is the  quadric of
equation  $(Q=0)$  in $\SS^3$.  This  quadric  admits two  transversal
rulings by  families of  great circles of  $\SS^3$. The  first ruling,
that  we call  \emph{left  ruling},  is the  family  of great  circles
$\{L_{(\lambda:\mu)}\}_{(\lambda:\mu)\in\RR\PP^1}$                where
$L_{(\lambda:\mu)}=\{[a:b:c:d]\in\partial\AdS_3\mid 
(a:c)=(b:d)=(\lambda:\mu)\mbox{  in }\RR\PP^1\}$.  The  second ruling,
that  we call  \emph{right ruling},  is  the family  of great  circles
$\{R_{(\lambda:\mu)}\}_{(\lambda:\mu)\in\RR\PP^1}$                where
$R_{(\lambda:\mu)}=\{[a:b:c:d]\in\partial\AdS_3\mid
(a:b)=(c:d)=(\lambda:\mu)\mbox{  in }\RR\PP^1\}$.  Through  each point
of $\partial\AdS_3$ passes one leaf 
of the left  ruling and one leaf of the right  ruling. Any leaf of
the  left ruling  intersects any  leaf of  the right  ruling  at two
antipodal points. 

The  elements   of  $O_0(2,2)$  preserve   the  left  and   the  right
ruling  of  $\partial\AdS_3$.  Hence,  for  each  element $\sigma$  of
$O_0(2,2)$, we can consider the action of $\sigma$ on the left and the
right   rulings.  This   defines   a  morphism   from  $O_0(2,2)$   to
$PSL(2,\RR)\times PSL(2,\RR)$. It is easy to see that this morphism is
onto, and that the kernel of  this morphism is a subgroup of order~$2$
of $O_0(2,2)$. As a consequence, we obtain an isomorphism from $O_0(2,2)$
to $SL(2,\RR)\times  SL(2,\RR)/(-\Id,-\Id)$ such that  the elements of
$SL(2,\RR)\times \{\pm Id\}/(-\Id,-\Id)$ preserve individually each circle
of the right ruling, and the elements of $\{\pm Id\}\times
SL(2,\RR)/(-\Id,-\Id)$ preserve  individually each leaf  of the left
ruling.

\begin{prop}
\label{p.geodesics-Klein}
The  geodesics  of  $\AdS_3$  are  the  connected  components  of  the
intersections of $\AdS_3$ with the great circles of $\SS^3$.  
\end{prop}

\begin{proof} 
By construction of $\AdS_3$, the geodesics of $\AdS_3$ are the images 
under     $\pi$     of    the     geodesics     of~$AdS_3$.      By
Proposition~\ref{p.geodesics},  the  geodesics   of  $AdS_3$  are  the
connected components of the intersections of $AdS_3$ with the
two-dimensional vector subspaces of $\RR^4$. The image under
$\pi$  of a  two-dimensional vector  subspace  of $\RR^4$  is a  great
circle   of    $\SS^3$.   Putting   everything    together,   we   get
Proposition~\ref{p.geodesics-Klein}.  
\end{proof} 

\begin{rema}
\label{r.type-geodesics-Klein}
Let     $\gamma$     be     a     geodesic     of     $\AdS_3$.     By
Proposition~\ref{p.geodesics-Klein},   $\gamma$  is   a  connected
component of  $\AdS_3\cap\widehat\gamma$, where $\widehat\gamma$  is a
geodesic of  $\SS^3$.  Moreover, Remark~\ref{r.type-geodesics} and
the proof of Proposition~\ref{p.geodesics-Klein} imply that:

\noindent -- if $\gamma$ is a timelike geodesic, then the great circle
$\widehat\gamma$ is contained in $\AdS_3$ and $\gamma=\widehat\gamma$,

\noindent  --  if  $\gamma$  is  lightlike,  then  the  great  circle
$\widehat\gamma$  is  tangent  to  $\partial\AdS_3$ at  two  antipodal
points  points  $p,-p$, and  $\gamma$  is  one  of the  two  connected
components of $\widehat\gamma\setminus\{p,-p\}$,

\noindent  --  if  $\gamma$   is  spacelike,  then  the  great  circle
$\widehat\gamma$ intersects $\partial\AdS_3$ transversally at four points 
$\{p_1,-p_1,p_2,-p_2\}$,  and $\gamma$  is one  of the  four connected
components of $\widehat\gamma\setminus\{p_1,-p_1,p_2,-p_2\}$. 
\end{rema}

\begin{rema}
\label{r.lightlike-geodesic}
Let  $q$  be a  point  of  $\partial\AdS_3$, and  $p$  be  a point  in
$\AdS_3$.  The great  $2$-sphere $S_q$ of $\SS^3$ which  is tangent to
the   quadric   $\partial\AdS_3$   at  $q$   is   $S_q=\{r\in\SS^3\mid
B_Q(q,r)=0\}$.   Consequently,   Remark~\ref{r.type-geodesics-Klein}
implies  that  there exists  a  lightlike  geodesic $\gamma$  passing
through $p$ and such that the ends of $\gamma$ in $\partial\AdS_3$ are
the points $q$ and $-q$ if and only if $B_Q(q,p)=0$. 
\end{rema}

Using Proposition  \ref{p.totally-geodesic} and the  same arguments as
in the proof of Proposition \ref{p.geodesics-Klein}, we obtain:

\begin{prop}
\label{p.totally-geodesic-Klein}
The  two-dimensional totally  geodesic subspaces  of $\AdS_3$  are the
connected components  of the intersections of $\AdS_3$  with the great
$2$-spheres of the sphere $\SS^3$. 
\end{prop}

Given  a point  $p$  in  $\AdS_3$, we  define  the \emph{dual  surface
  $p^\star$   of  $p$}   just  as   we  did   in  the   linear  model:
  $p^\star=\{q\in\AdS_3\mid B_Q(p,q)=0\}$. Note that the 
definitions  in  the linear  model  and  in  the Klein  model  are
coherent: if $\hat p$ is a point in $AdS_3$ such that $\pi(\hat p)=p$,
then the  dual surface  of $p$ is  the image  under $\pi$ of  the dual
surface       of       $\hat        p$.       We       denote       by
  $\overline{p^*}=\{q\in\AdS_3\cup\partial\AdS_3\mid B_Q(p,q)=0\}$ the
  closure on $p^*$ in $\AdS_3\cup\partial \AdS_3$.

\begin{rema}
In the sequel, we will indifferently denote the anti-de Sitter space by 
$AdS_3$ or $\AdS_3$. Mainly, we will have a preference to the first notation
when concerned with metric properties, and to the second one while discussing
convexity (see section~\ref{ss.convex-subsets}) or properties 
of the boundary at infinity $\partial\AdS_3$.
\end{rema}

\subsection{Affine domains in the anti-de Sitter space}
\label{ss.affine-domains}

By  an  \emph{open  hemisphere}   of  $\SS^3$,  we  mean  a  connected
component of $\SS^3$ minus a great $2$-sphere. Given an 
open    hemisphere    $U$,    we    say    that    a    diffeomorphism
$\phi:U\rightarrow\RR^3$  is  a \emph{projective  chart}  if $\phi$  maps
the  great circles  of $\SS^3$  (intersected with  $U$) to  the affine
lines of $\RR^3$. It is well-known that, for every open hemisphere $U$
of  $\SS^3$, there  exists an  projective chart $\phi:U\rightarrow\RR^3$.
This defines a locally projective structure on $\SS^3$, which induces a
locally projective structure on $\AdS_3$. The purpose of this subsection
is to define some particular projective charts of $\AdS_3$.

\bigskip

For every $p\in\AdS_3$, we consider the open hemisphere 
$U_p:=\{q\in\SS^3\mid B_Q(p,q)<0\}$, and the sets 
$$\begin{array}{lllll}
\cA_p & := & \{q\in\AdS_3\mid B_Q(p,q)<0\} & = & \AdS_3\cap U_p\\
\partial\cA_p & := & \{q\in\partial\AdS_3\mid B_Q(p,q)<0\} & = & 
\partial\AdS_3\cap U_p
\end{array}$$
Note that $\partial\cA_p$ is \emph{not} the boundary of $\cA_p$ in
$\SS^3$:  it is  the boundary  of $\cA_p$  in $U_p$.   Also  note that
$\cA_p$ is the connected component of $\AdS_3\setminus p^*$ containing
$p$, and  that $\cA_p\cup\partial\cA_p$ is the  connected component of
$(\AdS_3\cup\partial\AdS_3)\setminus\overline{p^*}$ containing $p$.  

Let  $p_0$ be  the point  of coordinates  $[1:0:0:0]$ in  $\SS^3$.  We
 observe  that 
$$U_{p_0}=\{[x_1:x_2:x_3:x_4]\in\SS^3 \mid  x_1>0\}$$
and we consider the diffeomorphism 
$$\begin{array}[t]{rrcl} 
\Phi_{p_0} : & U_{p_0} & 
\longrightarrow & \RR^3\\ 
& [x_1:x_2:x_3:x_4]   &  \longmapsto  &
(x,y,z)=\left(\frac{x_3}{x_1},\frac{x_4}{x_1},\frac{x_2}{x_1}\right)
\end{array}$$\\
Now,  given any point  $p\in\AdS_3$, we can find  an element
$\sigma_p$  of  $O_0(2,2)$,  such  that  $\sigma_p(p)=p_0$.  Then,  we
consider the diffeomorphism $\Phi_p:U_p\rightarrow\RR^3$ defined by 
$\Phi_p=\Phi_{p_0}\circ\sigma_p$.  

For every  $p\in\AdS_3$, the  diffeomorphism $\Phi_p$ maps  the domain
$\cA_p$  on   the  region  of   $\RR^3$  defined  by   the  inequation
$(x^2+y^2-z^2<1)$,  and   maps  $\partial\cA_p$  on   the  one-sheeted
hyperboloid of equation $(x^2+y^2-z^2=-1)$.  Moreover, $\Phi_{p_0}$ is
a   projective  chart   (as  the   usual   stereographic  projection),
\emph{i.e.} it maps  the great circles of $\SS^3$  to the affine lines
of $\RR^3$. 
Combining this with Proposition~\ref{p.geodesics-Klein}, we obtain
that,  for every  $p\in\AdS_3$, the  diffeomorphism $\Phi_p$  maps the
geodesics  of $\AdS_3$  to the  intersections of  the affine  lines of
$\RR^3$ with the set  $(x^2+y^2-z^2<1)$.  Similarly, $\Phi_p$ maps the
totally  geodesic subspaces of  $\AdS_3$ to  the intersections  of the
affine planes of $\RR^3$ with the set $(x^2+y^2-z^2<1)$. 

\begin{rema}
\label{r.type-geodesics-AdS+}
Let $\gamma$ be a geodesic of $\AdS_3$. Let $\gamma_p$ be the image
under   $\Phi_p$  of  $\gamma\cap\cA_p$.    According  to   the  above
remark, $\gamma_p$ is contained in 
an affine line $\widehat\gamma_p$ of $\RR^3$. Moreover, using
Remark~\ref{r.type-geodesics-Klein},  we  see  that:\\ 
-- if $\gamma$ is timelike, then the line $\widehat\gamma_p$
does not intersect the hyperboloid $(-x^2+y^2+z^2=1)$ and
$\gamma_p=\widehat\gamma_p$,\\ 
- if $\gamma$ is lightlike, then the affine line $\widehat\gamma_p$ is
tangent  to the hyperboloid  $(-x^2+y^2+z^2=1)$ at  one point  $q$ and
$\gamma_p$   is    one   of   the   two    connected   components   of
$\widehat\gamma_p\setminus q$,\\ 
- if  $\gamma$  is  spacelike,   then  the  line  $\widehat\gamma_p$  intersects
transversally  the   hyperboloid  $(-x^2+y^2+z^2=1)$  at   two  points
$q_1,q_2$  and   $\gamma$  is  the  bounded   connected  component  of
$\widehat\gamma\setminus\{q_1,q_2\}$.  
\end{rema}

The image under $\Phi_p$ of any geodesic of $\AdS_3$ is contained in an
affine line  of $\RR^3$. This implies  in particular that  there is no
closed  geodesic of $\AdS_3$  contained in  $\cA_p$. Moreover,  one can
prove that there  is no closed timelike curve in  $\cA_p$, so that the
causal structure of $\cA_p$ is more interesting than those of $\AdS_3$
(see Remark~\ref{r.causal-structure-AdS}).

\subsection{Convex subsets of $\AdS_3$}
\label{ss.convex-subsets}

Using the local projective structure of $\AdS_3$, we will define a
notion of convex subsets of $\AdS_3$.  

First,  we  define a  \emph{convex  subset of  $\SS^3$}  to  be a  set
$C\subset\SS^3$ such that: $C$ is contained in some open hemisphere $U$ 
of     $\SS^3$,    and    there     exists    some projective   chart
$\phi:U\rightarrow\RR^3$  such  that the  set  $\phi(C)$  is a  convex
subset of $\RR^3$.  

Note  that,  if  $C$ is  a  convex  subset  of $\SS^3$,  then,  for
\emph{every}  open  hemisphere  $V$  of $\SS^3$  containing  $C$,  and
\emph{every} projective chart  $\psi:V\rightarrow\RR^3$, the set $\phi(C)$
is a convex  subset of $\RR^3$. Moreover,  a set $C$ contained
in some  open hemisphere of $\SS^3$  is a convex subset  of $\SS^3$ if
and  only if the  positive cone  $\pi^{-1}(C)$ is  a convex  subset of
$\RR^4$   (recall   that   $\pi$   is  the   natural   projection   of
$\RR^4\setminus\{0\}$ on $\SS^3$).

Now, given a subset $E$ of $\SS^3$ such that $C$ is contained in some
open  hemisphere   of  $\SS^3$,  we  define   the  \emph{convex  hull}
$\Conv(C)$ of  the set $C$  to be the  intersection of all  the convex
subsets of $\SS^3$ containing $C$. Note
that,   if   $U$   is   an   open  hemisphere   containing   $C$   and
$\Phi:U\rightarrow\RR^3$ is a projective chart, the set $\Conv(C)$ is the 
image  under $\Phi^{-1}$  of the  convex hull  in $\RR^3$  of  the set
$\Phi(C)$.  Moreover, $\Conv(C)$ is also the image under
$\pi$ of the convex hull in $\RR^4$ of the positive cone $\pi^{-1}(C)$. 

Now, recall that $\AdS_3$ is contained in the sphere $\SS^3$, and let
$C$ be a subset of $\AdS_3$. We say that $C$ is a \emph{convex subset
of $\AdS_3$} if it  is convex as a subset of $\SS^3$.   We say that $C$
is a \emph{relatively convex subset $C$ of $\AdS_3$} if $C$ is the
intersection    of    $\AdS_3$     with    a    convex    subset    of
$\SS^3$.  Equivalently,  $C$  is   a  convex  subset  of  $\AdS_3$  if
$C=\Conv(C)$, and  $C$ is  a relatively convex  subset of  $\AdS_3$ if
$C=\Conv(C)\cap\AdS_3$.

\subsection{The $SL(2,\RR)$-model of the anti-de Sitter space}
\label{ss.SL2}

The linear model of the $3$-dimensional anti-de Sitter space is the
quadric   $\{(a,b,c,d)\in\RR^4\mid  -ad+bc=-1\}$   endowed   with  the
Lorentzian metric  induced by the  quadratic form $Q(a,b,c,d)=-ad+bc$.
Therefore, the anti-de 
Sitter   space  can  be   identified  with   the  group   of  matrices
$SL(2,\RR)=\left\{\left(\begin{array}{cc}a & b\\ c& d 
\end{array}\right)\in M(2,\RR)  \mid ad-bc=1\right\}$ endowed with
the Lorentzian metric induced by the quadratic form $-\det$ defined on
$M(2,\RR)$    by   $-\det\left(\begin{array}{cc}a    &   b\\    c&   d
  \end{array}\right)=ad-bc$.  

The quadratic form  $-\det$ on $M(2,\RR)$ is invariant  under left and
right  multiplication  by   elements  of  $SL(2,\RR)$  (actually,  the
Lorentzian metric induced by $-\det$  is a multiple of the Killing form
of the Lie group $SL(2,\RR)$). This implies that the isometry group of
$(SL(2,\RR),-\det)$   is    $SL(2,\RR)\times   SL(2,\RR)$   acting   on
$SL(2,\RR)$ by left and right multiplication, \emph{i.e.}  acting by 
$(g_1,g_2).g=g_1gg_2^{-1}$.  

\subsection{Causal structure of the anti-de Sitter space}
\label{ss.conformal-structure}

Denote $dt^2$ the standard Riemannian metric on the circle $\SS^1$, by
$ds^2$ the standard Riemannian metric on the $2$-dimensional 
sphere $\SS^2$, by $\DD^2$  the open upper-hemisphere of $\SS^2$,
and by $\overline{\DD^2}$ the closure of $\DD^2$.  We will prove that
$\AdS_3$      has      the      same     causal      structure      as
$(\SS^1\times\DD^2,-dt^2+ds^2)$.  More precisely:

\begin{prop}
\label{p.conformal-structure}
There exists a diffeomorphism $\Psi:\AdS_3\rightarrow\SS^1\times\DD^2$
such  that  the   pull  back  by  $\Psi$  of   the  Lorentzian  metric
$-dt^2+ds^2$ defines  the same causal structure as  the original metric
of  $\AdS_3$, that  is,  the two  metrics  are in  the same  conformal
class. Moreover, the diffeomorphism $\Psi$ can be extended to a 
diffeomorphism $\overline{\Psi}:\AdS_3\cup\partial\AdS_3 \rightarrow
\SS^1\times\overline{\DD^2}$.  
\end{prop}

To   prove  this,  we   will  embed
$\AdS_3$ in the so-called \emph{three-dimensional Einstein universe}. 
Denote by $(x_1,x_2,x_3,x_4,x_5)$ the standard coordinates on
$\RR^5$, consider the quadratic form $\widetilde Q$ on $\RR^5$ defined by
$\widetilde     Q(x_1,x_2,x_3,x_4,x_5)=-x_1^2-x_2^2+x_3^2+x_4^2+x_5^2$,
denote by $\SS^4$ the quotient of $\RR^5\setminus\{0\}$ by positive
homotheties, and by $\widetilde\pi$ the natural projection of
$\RR^5\setminus\{0\}$ on $\SS^4$.  Then, the \emph{three-dimensional
 Einstein space}, 
  denoted  by $Ein_3$, is the image  under $\widetilde\pi$ of
the quadric $(\widetilde Q=0)$.  There is a natural conformal class of
Lorentzian metrics on $Ein_3$, defined as follows: \\
--- Given  an  open  subset  $U$  of  $Ein_3$,  and  a  local  section
$\sigma:U\rightarrow    \RR^5\setminus\{0\}$    of   the    projection
$\widetilde\pi$, we 
define a  Lorentzian metric  $g_\sigma$ on $U$  as follows.  For every
point $p\in U$  and every vector $v\in T_p Ein_3$,  we choose a vector
$\hat v\in T_{\sigma(p)}\RR^5$ such that $d\widetilde\pi(\sigma(p)).\hat v=v$.
The quantity $\widetilde Q(\hat v)$  does not depend  on the choice  of the
vector $\hat v$: indeed, the vector $\hat v$ is tangent to the quadric
$(\widetilde Q=0)$, the  vector $\hat v$  is defined up  to the addition  of an
element     of    $\widetilde\pi^{-1}(p)$,    and     the    half-line
$\widetilde\pi^{-1}(p)$ is contained in the $\widetilde Q$-orthogonal of
the tangent space of the quadric $(\widetilde Q=0)$ at $\sigma_p$.  We
set $g_\sigma(v):=\widetilde Q(\hat v)$. \\
--- The conformal class of the  metric $g_\sigma$ does not  depend on
the section $\sigma$. Indeed, if $\sigma$ and $\sigma'$ are
two sections of the projection $\widetilde\pi$ defined on $U$, then we
have $g_{\sigma'}=\lambda^2.g_\sigma$,  where  $\lambda:U\rightarrow\RR$ is
the function such that $\sigma'=\lambda.\sigma$. 

\begin{proof}[Proof of Proposition~\ref{p.conformal-structure}]
Let   $A=\{[x_1:x_2:x_3:x_4:x_5]\in   Ein_3\mid   x_5>0\}$,  and   let
$\partial A$ be the boundary  of $A$.  We will consider two particular
sections  of the  projection $\widetilde\pi$.  First, we  consider the
section 
$\sigma$, defined on $A$, whose image is contained in the
affine hyperplane $x_5=1$.  The anti-de Sitter space
$\AdS_3$  is isometric  to the  set $A$  equipped with  the Lorentzian
metric $g_{\sigma}$:  the most natural isometry  is the diffeomorphism
$\Phi$                            defined                           by
$\Phi([x_1:x_2:x_3:x_4])=[x_1:x_2:x_3:x_4:1]$.   Now, we  consider the
section $\sigma'$, defined on the
whole of $Ein_3$, whose image is contained in the Euclidean sphere
$x_1^2+x_2^2+x_3^2+x_4^2+x_5^2=2$.  The set 
$A$ equipped with the Lorentzian metric $g_{\sigma'}$ is isometric to
the set $\{(x_1,x_2,x_3,x_4,x_5)\in\RR^5\mid x_1^2+x_2^2=1\mbox{ ,
}x_3^2+x_4^2+x_5^2=1\mbox{    ,    }x_5>0\}\simeq   \SS^1\times\DD^2$
equipped          with          the         Lorentzian          metric
$-(dx_1^2+dx_2^2)+(dx_3^2+dx_4^2+dx_5^2)\simeq  -dt^2+ds^2$:  the most
natural  isometry  is   the  diffeomorphism  $\Phi'=\sigma'_{|A}$.  We
consider                       the                      diffeomorphism
$\Psi:=\Phi'\circ\Phi:\AdS_3\rightarrow\SS^1\times\DD^2$.   Since  the
metric $g_\sigma$  and $g_{\sigma'}$  are conformally equivalent,
the  pull back  by $\Psi$  of the  metric $-dt^2+ds^2$  is conformally
equivalent to the original metric of $\AdS_3$. 

The  diffeomorphism  $\Phi$  can   be  extended  to  a  diffeomorphism
$\overline{\Phi}:\AdS_3\cup\partial\AdS_3\rightarrow A\cup\partial A$:
for   every   $[x_1:x_2:x_3:x_4]$   in   $\partial\AdS_3$,   we   have
$\overline{\Phi}([x_1:x_2:x_3:x_4])=[x_1:x_2:x_3:x_4:0]$.         The
diffeomorphism   $\Phi'$   can  be   extended   to  a   diffeomorphism
$\overline{\Phi'}:A\cup\partial                            A\rightarrow
\SS^1\times\overline{\DD^2}$: we have
$\overline{\Phi'}=\sigma'_{A\cup\partial      A}$.      Hence,     the
diffeomorphism $\Psi$ can be extended to a
diffeomorphism $\overline{\Psi}=\overline{\Phi}\circ\overline{\Phi'} :
\AdS_3\cup\partial\AdS_3 \rightarrow \SS^1\times\overline{\DD^2}$. 
\end{proof}

\paragraph{Causal structure on $\AdS_3\cup\partial\AdS_3$.}  
Let $\overline{g}$     be     the     Lorentzian     metric     on
$\AdS_3\cup\partial\AdS_3$,  obtained by  pulling back  the Lorentzian
metric $-dt^2+ds^2$ defined on $\SS^1\times\overline{\DD^2}$ by the
diffeomorphism $\overline{\Psi}$. The Lorentzian metric $\overline{g}$
defines  the same  causal  structure on  $\AdS_3$  as the  original
metric of  $\AdS_3$. From now on,  we endow $\AdS_3\cup\partial\AdS_3$
with  the  causal  structure defined  by  the  metric
$\overline{g}$.   This causal structure allows  us to
speak of timelike,   lightlike    and   spacelike    objects   in
$\AdS_3\cup\partial\AdS_3$.  In   particular,  we  can   consider  the
causal  structure induced  on  the quadric  $\partial\AdS_3$. Given  a
point $q\in\partial\AdS_3$, it is easy to verify that the lightcone of
$q$ for this conformally Lorentzian
structure is the union of the leaf of the left ruling and of the circle
of the right ruling passing through $q$. 

\begin{rema}
\label{r.conformal-structure}
Let  $p_0$ be  the point  of coordinates  $[1:0:0:0]$ in  $\SS^3$. 
Recall   that  $\cA_{p_0}\cup\partial\cA_{p_0}$   is  the   subset  of
$\AdS_3\cup\partial\AdS_3$ defined by the inequation $(x_1>0)$. Hence, 
the diffeomorphism $\overline{\Psi}$ defined above             maps
$\cA_{p_0}\cup\partial\cA_{p_0}$    on    $\{(x_1,x_2,x_3,x_4,x_5)\mid
x_1^2+x_2^2=1\mbox{ , }x_1>0\mbox{ , }x_3^2+x_4^2+x_5^2=1\mbox{ , }x_5
\geq 0\}\simeq (-\pi/2,\pi/2)\times\overline{\DD^2}$. 
\end{rema}

\begin{coro}
\label{c.conformal-structure}
For every $p\in\AdS_3$, the domain $\cA_p\cup\partial\cA_p$
has the same causal structure as the Lorentzian space 
$\left((-\pi/2,\pi/2)\times\overline{\DD^2},-dt^2+ds^2\right)$.  
\end{coro}

\begin{proof}
Since $O(2,2)$ acts transitively on $\AdS_3$, it is enough to consider
the case where $p$ is the point of coordinates $[1:0:0:0]$.  This case
follows from Proposition~\ref{p.conformal-structure} and
Remark~\ref{r.conformal-structure}. 
\end{proof}

The two following propositions will play some fundamental roles in the
proof of Theorem~\ref{t.CMC-foliation}:

\begin{prop}
\label{p.lightcone}
Let  $p$   be  a   point  in   $\AdS_3$,  and  $q$   be  a   point  in
$\partial\cA_p$.  A  point $r\in\cA_p\cup\partial\cA_p$ can  be joined
from $q$ by a timelike (resp.  causal) curve if and only if $B_Q(q,r)$
is positive (resp.  non-negative).  
\end{prop}

\begin{proof}
Since  $O(2,2)$ acts  transitively  on $\AdS_3$,  we  can assume  that
$p=[1:0:0:0]$.  There  exists a timelike  curve joining $q$ to  $r$ in
$\cA_p\cup\partial\cA_p$ if and only  if there exists a timelike curve
joining         $\Psi_p(q)$          to         $\Psi_p(r)$         in
$\left((-\pi/2,\pi/2)\times\overline{\DD^2},-dt^2+ds^2\right)$.  
We see $\left ((-\pi/2,\pi/2)\times\overline{\DD^2},-dt^2+ds^2\right)$ as
the set $\{(x_1,x_2,x_3,x_4,x_5)\in\RR^5\mid x_1^2+x_2^2=1\mbox{ ,
}x_1>0\mbox{  , }x_3^2+x_4^2+x_5^2=1\mbox{  , }x_5>0\}$  equipped with
the  metric $-(dx_1^2+dx_2^2)+(dx_3^2+dx_4^2+dx_5^2)$. Coming  back to
the  definition  of  the  diffeomorphism  $\Psi_p$  (see  the  proof of
Proposition~\ref{p.conformal-structure}),  we observe  that $B_Q(q,r)$
and $B_{\widetilde Q}\left(\Psi_p(q),\Psi_p(r)\right)$ have the same sign. 
Moreover, it is clear that the points $\Psi_p(q)$ and $\Psi_p(r)$
can be joined by a timelike (resp.  causal) curve in
$\left((-\pi/2,\pi/2)\times\overline{\DD^2},-dt^2+ds^2\right)$  if and
only  if   $\widetilde  Q(\Psi_p(q)-\Psi_p(r))$  is   negative  (resp.
non-positive).  Finally, notice that the quantity $\widetilde
Q\left(\Psi_p(q)-\Psi_p(r)\right)$          and         $B_{\widetilde
  Q}\left(\Psi_p(q),\Psi_p(r)\right)$   have  opposite   signs  (since
$\widetilde                          Q\left(\Psi_p(q)\right)=\widetilde
Q\left(\Psi_p(r)\right)=0$).  Putting everything  together,  we obtain
the proposition. 
\end{proof}

\begin{rema}
\label{r.past-future-plane}
Let $p$ be a point in $\AdS_3$.  Let $P$ be a totally
geodesic spacelike subspace of  $\cA_p$ (by  such we  mean the  intersection of
$\cA_p$ with a totally geodesic spacelike subspace of $\AdS_3$). Then, 
$P$ divides $\cA_p$ into two closed regions: the past of $P$ in
$\cA_p$ and the future of $P$ in $\cA_p$.  
\end{rema}

\begin{proof}
We identify  $\cA_p$ and $P$  with their images under  the embedding
$\Phi_p$.  Then,
$P$  is   the  intersection  of  $\cA_p$  (\emph{i.e.}    of  the  set
$(-x^2+y^2+z^2<1)$) with an affine  plane $\widehat P$ of $\RR^3$.  We
consider the  two regions of  $\cA_p$ defined as the  intersections of
$\cA_p$    with   the   closures    two   connected    components   of
$\RR^3\setminus\widehat  P$. Since $P$  is spacelike  and connected,
the  past  (resp.  the  future)  of  $P$  in $\cA_p$  is  necessarly
contained    in     one    of    these     two    regions.     Finally,
Remark~\ref{r.type-geodesics-AdS+} implies that, for every point
$q\in\cA_p$, there exists  a timelike geodesic joining $q$  to a point
of $P$.  Hence,  the union of the  past and the future of  $P$ must be
equal to $\cA_p$. The proposition follows.  
\end{proof}

\section{Globally hyperbolic spacetimes}
\label{mgh}

All  along this  section,  we consider  a  maximal globally  hyperbolic
spacetime $M$, 
locally modelled on $AdS_3$, with closed orientable Cauchy surfaces. 
All the Cauchy surfaces have the same genus, that we denote by $g$.
We denote by $\widetilde M$ the universal covering
of  $M$. We  choose a  Cauchy surface  $\Sigma_0$ in  $M$, and the lift
$\widetilde\Sigma_0$  of $\Sigma_0$  in $\widetilde  M$. Since  $M$ is
locally modelled on $AdS_3$, we can consider the developping
map   $\cD:\widetilde{M}   \rightarrow    AdS_3$   and   the  holonomy
representation   $\rho:\pi_1(M)=\pi_1(\Sigma_0)\rightarrow   O_0(2,2)$
(see section~\ref{s.GX}).  

Let                 $S_0=\cD(\widetilde\Sigma_0)$,                 and
$\Gamma=\rho(\pi_1(M))$. Identifying $O_0(2,2)$ with $SL(2,\RR)\times 
SL(2,\RR)/(-\Id,-\Id)$ (see subsection~\ref{ss.projective-model}), we
can see $\rho$ as a 
representation              of              $\pi_1(M)$              in
$SL(2,\RR)\times SL(2,\RR)$.  Then, we will denote by $\rho_L$ and
$\rho_R$ the  representations of  $\pi_1(M)$ in $SL(2,\RR)$  such that
$\rho=\rho_L\times\rho_R$. 

In subsection~\ref{s.S_0},  we  will  study  the surface  $S_0$  and  its
boundary    $\partial   S_0$    in    $\AdS_3\cup\partial\AdS_3$.   
In particular,  we will show  that $S$ cannot  be a sphere,  i.e., its
genus  $g$  is positive.   The  results  of  this subsection  are  not
original: most of them are contained in Mess preprint
(\cite{Mes}). Yet, we will provide a  proof of each result to keep our
paper as self-contained as possible (by the way, using the
conformal equivalence of $\AdS_3\cup\partial\AdS_3$ with
$(\overline{\DD}^2\times\SS^1,-dt^2+ds^2)$, we were able to simplify 
  some of the proofs of Mess). 

In
subsection~\ref{s.Cauchy-development}, we study the Cauchy development
$D(S_0)$ of the 
surface $S_0$. In particular, we prove that $M$ is isometric to the
quotient            $\Gamma\backslash            D(S_0)$.

\subsection{The spacelike surface $S_0$}
\label{s.S_0}

The  purpose of  this subsection  is to  collect as  many  information as
possible on the surface $S_0$. In particular, we will prove that $S_0$
is  an open  disc properly  embedded  in $\AdS_3$,  that the  closure
$\overline{S_0}$ of $S_0$ in $\AdS_3\cup\partial\AdS_3$ is a closed 
topological   disc,   and  that   $\overline{S_0}$   is  an   achronal
set.

The  Lorentzian   metric  of  $M$  induces  a   Riemannian  metric  on
the Cauchy surface $\Sigma_0$, which can be lifted to get a Riemannian
metric on  $\wt\Sigma_0$. Since $\Sigma_0$ is  compact, the Riemannian
metrics on $\Sigma_0$ and $\wt\Sigma_0$ are complete. The developping
map  $\cD$  induces  a  locally  isometric immersion  of  the  surface
$\wt\Sigma_0$ in $AdS_3$. It turns out that this immersion is
automatically a proper embedding: 

\begin{prop}
\label{p.embedded}
The   surface   $S_0$  is   an   open   disc   properly  embedded   in
$AdS_3$. Moreover, every timelike geodesic of $AdS_3$ intersects the
surface $S_0$ at exactly one point. 
\end{prop}

\begin{proof}
We consider the projection $\zeta:AdS_3\rightarrow\RR^2$, defined by 
$\zeta(x_1,x_2,x_3,x_4)=(x_3,x_4)$.  Observe  that  the  fibers  of  the
projection $\zeta$ are the orbits  of a timelike killing vector field of
$\AdS_3$. We endow $\RR^2$ with the Riemannian metric $g_{\zeta}$ defined
as follows. Given a point $q\in\RR^2$ and a vector $v\in T_q\RR^2$, we
choose a  point $\hat  q\in \zeta^{-1}(q)$, and  we consider  the unique
vector $\hat v\in T_{\hat q} AdS_3$ such that $d\zeta_{\hat q}.\hat v=v$
and such that  $\hat v$ is orthogonal to  the fibers $\zeta^{-1}(q)$. We
define $g_{\zeta}(v)$  to be to the norm  of the vector $\hat  v$ for the
Lorentzian metric of  $AdS_3$. This definition does not  depend on the
choice of the point $\hat q$, since the fibers of $\zeta$ are the orbits
of a killing  vector field. It is easy to  verify that $\RR^2$ endowed
with   the   metric   $g_\zeta$   is   isometric   to   the   hyperbolic
plane.

\smallskip

\noindent  \emph{Claim  1.  Given  any  point  $q\in  AdS_3$  and  any
  spacelike  vector $v$  in  $T_q  AdS_3$, the norm of the  vector
  $d\zeta_q(v)$ for  the metric $g_{\zeta}$ is  bigger than the  norm of $v$
  in $AdS_3$.}  

\smallskip

\noindent Indeed, write $v=u+w$ where  $u$ is tangent to
the fiber of the projection $\zeta$ (in particular, $u$ is timelike) and $w$ is
orthogonal to this fiber. On the one hand, by definition of $g_{\zeta}$, the
norm of the vector $d\zeta_q(v)$ for the metric $g_{\zeta}$ is equal to the
norm of $w$ in $AdS_3$.  On the other hand, the norm of $v$ in $AdS_3$
is less  than the norm of  $w$, since $u$ is  timelike. This completes
the proof of claim~1.  

\smallskip

\noindent  \emph{Claim  2.   For  every  locally  isometric  immersion
  $f:\wt\Sigma_0\rightarrow   AdS_3$,   the    map   $\zeta\circ   f   :
  \wt\Sigma_0\rightarrow \RR^2$ is an homeomorphism. In particular, the
  surface $f(\wt\Sigma_0)$  intersects each fiber of  $\zeta$ at exactly
  one point.}  

\smallskip

\noindent  By the first claim, the map $\zeta\circ f$ is locally distance
increasing (when the surface $\wt\Sigma_0$ is endowed with its Riemannian
metric, and $\RR^2$  is endowed with the metric  $g_{\zeta}$).  Since the
Riemannian metric of $\Sigma_0$ is complete, this implies 
that  $\zeta\circ f:\wt\Sigma_0\rightarrow \RR^2$  has the  path lifting
property, and thus is a covering map.  Since $H$ is simply connected, this
implies   that    $\zeta\circ   f:\wt\Sigma_0\rightarrow \RR^2$    is   an
homeomorphism.  This completes the proof of claim~2. 

\smallskip

\noindent  Applying claim  2 with  $f$ being  the developping  map
$\cD$, we obtain that $\cD:\wt\Sigma_0\rightarrow AdS_3$  is a proper
embedding, and that $\wt\Sigma_0$ is homeomorphic to $\RR^2$ (and thus
homeomorphic    to    an    open    disc).    Hence,    the    surface
$S_0:=\cD(\wt\Sigma_0)$   is  an  open   disc  properly   embedded  in
$AdS_3$. Now, let $\gamma$ be a timelike geodesic of $AdS_3$. Observe that
the circle $\zeta^{-1}(0,0)$ is a timelike geodesic of $AdS_3$. Since
$O(2,2)$  acts  transitively  on  the  set  of  timelike  geodesic  of
$AdS_3$,    there     exists    $\sigma\in    O(2,2)$     such    that
$\sigma(\gamma)=\zeta^{-1}(0,0)$;  in particular, $\sigma(\gamma)$  is a
fiber   of   the   projection    $\zeta$.    Applying   claim   2   with
$f=\sigma^{-1}\circ \cD$, we obtain that the surface 
$\sigma^{-1}(S_0)=\sigma^{-1}\circ \cD(\wt\Sigma_0)$ intersects each 
fiber  of  $\zeta$ at  exactly  one  point.   Hence, the  surface  $S_0$
intersects the geodesic $\gamma$ at exactly point. 
\end{proof}

\begin{rema}
\label{r.other-surface}
Proposition~\ref{p.embedded} is still valid if 
$\Sigma_0$ is replaced by another Cauchy surface of $M$. 
\end{rema} 

\begin{rema}
\label{r.sphere-impossible}
The     proof    of     Proposition~\ref{p.embedded}     shows    that
$\widetilde\Sigma_0$ is homeomorphic to a disc.  Hence, there does not
exist any globally hyperbolic  spacetime, locally modelled on $AdS_3$,
with closed orientable Cauchy surfaces of genus~$0$. 
\end{rema}


Now,   we    will   use   the    conformal   equivalence   between
$\AdS_3\cup\partial\AdS_3$ and
$(\SS^1\times\overline{\DD^2},-dt^2+ds^2)$. Let us start by some remarks:

\begin{rema}
\label{r.spacelike-surface} 
\noindent (i) Let $S$ be a spacelike (resp. non-timelike) surface
in $(\SS^1\times\overline{\DD^2},-dt^2+ds^2)$. Then, every point of
$S$  has   a  neighbourhood in  $S$  which  is  the   graph  of  a
contracting\footnote{We recall that,  given two metric
  spaces  $(E,d)$   and  $(E',d')$,  a   mapping  $f:(E,d)\rightarrow
  (E',d')$ is said  to be \emph{contracting} if $d'(f(x),f(y))<d(x,y)$
  for every $x\neq y$.}  (resp.   $1$-Lipschitz)  mapping  $f:(U,ds^2)\rightarrow
(\SS^1,dt^2)$, where $U$ is an open subset of $\overline{\DD^2}$.  

\smallskip

\noindent  (ii)  Every  properly  embedded  spacelike  (resp.  non
timelike) surface in $(\SS^1\times\DD^2,-dt^2+ds^2)$ is the graph of a
contracting (resp.  $1$-Lipschitz) mapping $f:(\DD^2,ds^2)\rightarrow
(\SS^1,dt^2)$. 

\smallskip

\noindent  (iii) Of course,  (i) and  (ii) remain  true if  we replace
$\SS^1$ by $(-\pi/2,\pi/2)$.
\end{rema}

\begin{proof}
Item (i) is an immediate consequence of the product structure of
$(\SS^1\times\overline{\DD^2},-dt^2+ds^2)$.  To  prove  (ii),  we
consider  a  properly  embedded  spacelike (resp.   non-timelike)
surface $S$  in $(\SS^1\times\DD^2,-dt^2+ds^2)$. Let  $p_2$ be the
projection of  $\SS^1\times\DD^2$ on $\DD^2$.  Using item (i)  and the
fact  that  $S$  is  properly  embedded,  it  is  easy  to  show  that
$p_2:S\rightarrow\DD^2$ is a covering map. Hence, $p_2:S\rightarrow\DD^2$
is a homeomorphism, and the surface $S$ is the graph of a mapping
$f:\DD^2\rightarrow\SS^1$. By item (i), the mapping $f$ is contracting
(resp. $1$-Lipschitz).  
\end{proof} 

\begin{rema}
\label{r.causal-curve} 
In the same  vein, we observe that timelike  (resp. causal) curves are
represented in $(\SS^1\times\overline{\DD^2},-dt^2+ds^2)$  by
graphs    of   contracting   (resp.    $1$-Lipschitz)    mappings
$g:(J,dt^2)\rightarrow (\overline{\DD^2},ds^2)$, where $J$ is a 
subinterval of $\SS^1$. 
\end{rema}

Putting                Proposition~\ref{p.embedded}                and
Remark~\ref{r.spacelike-surface} together, we obtain the following: 

\begin{prop}
\label{graphe}
Any conformal equivalence between $\AdS_3$ and 
$(\SS^1\times\DD^2,-dt^2+ds^2)$ maps the surface $S_{0}$ to
the graph  of a contracting mapping $f: \DD^{2} \rightarrow \SS^{1}$.  
\end{prop}

Now, let us denote by $\overline{S_0}$ the closure of the surface
$S_0$ in $\AdS_3\cup\partial\AdS_3$.

\begin{coro}
\label{p.closed-topological-disc}
Any conformal equivalence between $\AdS_3\cup\partial\AdS_3$ and
$(\SS^1\times\overline{\DD^2},-dt^2+ds^2)$     maps     the    closure
$\overline{S_0}$ of the surface $S_0$  to the graph of a $1$-Lipschitz
mapping               $\overline{f}:(\overline{\DD^2},ds^2)\rightarrow
(\SS^1,dt^2)$, which is contracting in restriction to the open
disc $\DD^2$.  In particular, $\overline{S}_0$ is a closed topological
disc.  
\end{coro}

\begin{proof}
The  result follows  from Proposition~\ref{graphe}  and from  the fact
that  any contracting mapping  from $(\DD^2,ds^2)$  to $(\SS^1,dt^2)$
can    be    extended     as    a    $1$-Lipschtiz    mapping    from
$(\overline{\DD^2},ds^2)$ to $(\SS^1,dt^2)$.  
\end{proof}

Proposition~\ref{p.embedded}                                        and
corollary~\ref{p.closed-topological-disc} imply 
that   the  boundary   $\partial  S_0$   of  the   surface   $S_0$  in
$\AdS_3\cup\partial\AdS_3$  is  a   topological  simple  closed  curve
contained  in $\partial\AdS_3$.   Of course,  the curve  $\partial S_0$
must be invariant by the holonomy group $\Gamma=\rho(\pi_1(M))$.  

\begin{rema}
\label{r.curve}
According to the proof of Proposition~\ref{p.embedded}, the surface
$S_0$     intersects     each      fiber     of     the     projection
$\zeta:AdS_3\rightarrow\RR^2$                defined                by
$\zeta((x_1,x_2,x_3,x_4))=(x_3,x_4)$.   This  implies  that the  curve
$\partial S_0$
intersects       each        fiber       of       the       projection
$\zeta:\partial\AdS_3\rightarrow\SS^1$            defined           by
$\zeta([x_1:x_2:x_3:x_4])=[x_3:x_4]$.  

Futhermore,    if    we    identify 
$\AdS_3\cup\partial\AdS_3$    with
$(\SS^1\times\overline{\DD^2},-dt^2+ds^2)$,   then
the   curve
$\partial  S_0$  is identified  with  the  graph  of a  mapping  from
$\partial\DD^2$  to $\SS^1$.  This  implies, in  particular, that  the
curve $\partial S_0$ is not null-homotopic in $\partial\AdS_3$. 
\end{rema}

Thanks to Remark~\ref{r.spacelike-surface}, we can define a notion
of spacelike topological surface in $\AdS_3\cup\partial\AdS_3$:

\begin{defi}
\label{d.spacelike}
Let  $S$  be a  topological  surface  (with  or without  boundary)  in
$\AdS_3\cup\partial\AdS_3$.
Using the conformal equivalence between $\AdS_3\cup\partial\AdS_3$ and
$(\SS^1\times\overline{\DD^2},-dt^2+ds^2)$,  we  can   see  $S$  as  a
surface  in  $\SS^1\times\overline{\DD^2}$.   We  will  say  that  the
topological  surface  $S$  is  \emph{spacelike}  (resp.  \emph{non
  timelike}) if every point of $S$ has a neighbourhood in $S$ which is
the   graph   of  a   contracting   (resp.   $1$-Lipschitz)   mapping
$f:(U,ds^2)\rightarrow (\SS^1,dt^2)$,  where $U$ is an  open subset of
$\overline{\DD^2}$.  
\end{defi}

With   this  definition,  $\overline{S_0}$   is  a   non-timelike
topological    surface    in   $\AdS_3\cup\partial\AdS_3$.  

\begin{prop}
\label{pasdelien}
Every lightlike geodesic intersects the surface $S_0$ at most
once.  Moreover, if  a lightlike geodesic has one  of its endpoints on
the curve $\partial S_0$, then this geodesic does not intersect $S_0$.
\end{prop}

\begin{proof}
Let $p$ be a point on the surface $S_0$, and $\gamma$
a lightlike geodesic containing $p$. Denote by 
$d$ the distance function on the hemisphere $\overline{\DD}^{2}$, and 
let $p_{0}$ be the center of the hemisphere, i.e. the unique point
for which $d(p_{0}, q) = \pi/{2}$ for any point $q$ in $\partial \DD^{2}$.
Select a conformal equivalence 
$\AdS_3\cup\partial\AdS_3 \approx (\SS^1\times\overline{\DD^2},-dt^2+ds^2)$
for which $p$ is identified with $(0, p_{0})$ and $\cA_p$ with $]-\pi/2, \pi/2[\times \DD^2$.
Then, $\overline{S_0}$ is represented as the graph of a $1$-Lipschitz mapping
$f$ for which $f(p_0)=0$. On the other hand, as every lightlike geodesic
containing $p$, $\gamma$ is contained in $\cA_p$ and is represented by 
a curve $(d(p_{0}, r), r)$, where $r$ describes
a geodesic in $\DD^2$ containing $p_{0}$. Since the restriction of $f$
to $\DD^2$ is contracting, it follows immediatly that $\gamma$ does not
contain another point of $S_{0}$ than $p$. The first statement in the proposition follows.

Assume  now   that  one  of  the   two  end  points   of  $\gamma$  is
$(f(q),q)\in\partial S_0$. Then, $d(q,p_{0}) = \pi/2 = f(q)$, and 
since $f$ is $1$-Lipschitz, for any point $r$ on the geodesic of $\DD^2$ under
consideration, we must have $d(p_{0}, r) = f(r)$. This is impossible, since the restriction
of $f$ to $\DD^{2}$ is contracting. 
\end{proof}

\begin{prop}
\label{p.curve-in-affine-domain}
For every $p\in S_0$, the surface $\overline{S_0}$ is contained
in the affine domain $\cA_p\cup\partial\cA_p$. 
\end{prop}

\begin{proof}
We keep the notation used in the proof of the previous lemma.
 It follows immediatly that the maximum value of $f$ is at
most $\pi/2$, and its minimum value is at least $-\pi/{2}$. In other words,
$\overline{S}_{0}$ is contained in the closure of $\cA_p$. Moreover, 
in the proof above we have actually shown that $f$ does not attain the
values $\pi/{2}$, $-\pi/2$. The proposition follows.
\end{proof}

\begin{prop}
\label{p.ajout}
For           every          $p\in\AdS_3$           such          that
$\overline{S_0}\subset\cA_p\cup\partial\cA_p$, the surface 
$\overline{S_0}$ is an   achronal     subset    of
$\cA_p\cup\partial\cA_p$  (\emph{i.e.}  a  timelike curve  contained  in
$\cA_p\cup\partial\cA_p$  cannot  intersect  $\overline{S_0}$  at  two
distinct points). Moreover, if two points in $\overline{S_{0}}$ are causally
 related,  then  they belong  to  a  lightlike  geodesic of  $\partial
 \AdS_3$ contained in $\partial \overline{S_0}$. 
\end{prop}

\begin{proof}
We keep the notations used in the proof of Proposition \ref{pasdelien}
(except that $(0, p_0)$ is not assumed now to belong to $S_0$, i.e.,
the mapping $f$ admitting $\overline{S_0}$ as graph does not necessarly
vanish at $p_0$). 
A future oriented causal curve in $\cA_p$ is represented by a curve 
$(g(t), r(t))$ where $g$ satisfies: $g(t)-g(s) \geq d(r(t), r(s))$. 
Assume the existence of $t < t'$ such that $g(t) = f(r(t))$ and
$g(t') = f(r(t'))$. Then: 
$$\mid f(r(t')) - f(r(t))\mid \leq  d(r(t), r(t')) \leq g(t') - g(t) =
f(r(t')) - f(r(t))$$ 
Therefore,   all  these   inequations  are   equalities.   According  to
Proposition \ref{pasdelien}, 
it follows that $(g(t), r(t))$ and $(g(t'), r(t'))$ belong both to $\cA_p$.
Moreover, it follows that for every $s$ in $[t,t']$, $f(r(s)) = g(r(s)) = f(r(t)) + d(r(s),r(t))$.
The proposition follows.
\end{proof}

\begin{rema}
\label{r.past-future-surface}
Let  $p$  be a  point  such  that the  surface  $S_0$  is contained  in
$\cA_p$.   Proposition~\ref{p.embedded} implies  that every  point of
$\cA_p$ is either in the past\footnote{Here, by ``past'', we mean the ``past in
$\cA_p$'': a  point $q$ is in the  past of the surface  $S_0$ if there
exists a future-directed causal  curve contained in $\cA_p$ going from
$S_0$ to $q$. Similarly for the future.} or in the future of 
the  surface  $S_0$.   Moreover, it should be clear
to the reader that, according to corollary~\ref{p.closed-topological-disc} 
and Proposition~\ref{p.ajout}, 
a  point of $\cA_p$ cannot  be simultaneously in the  past and in
the future of the surface, except if it is on the surface $S_0$. 
\end{rema}

\subsection{Cauchy  development of  the surface $S_0$}
\label{s.Cauchy-development}

In this subsection, we study the Cauchy development $D(S_0)$ of the surface
$S_0$ in  $AdS_3$. The main goal  of the subsection is  to prove that
$M$ is isometric to a quotient $\Gamma\backslash D(S_0)$. 

Let  us first recall  the definition  of the  Cauchy development of a
spacelike  surface.  Given a  spacelike surface  $S$ in  $AdS_3$, the
\emph{past  Cauchy development}  $D^-(S)$ of  $S$  is the  set of  all
points $p\in AdS_3$  such that every  future-inextendable causal curve
through  $p$  intersects  $S$.  The \emph{future  Cauchy  development}
$D^+(S)$ of $S$ is defined similarly. The \emph{Cauchy development} of
$S$ is  the set $D(S):=D^-(S)\cup  D^+(S)$.  It is well-known  and not
difficult to prove 
that $D(S)$ is a connected open domain. The following lemma provides a
more tractable definition of $D(S)$: 

\begin{lemma}
\label{l.def-alternative}
Let $S\subset AdS_3$ be a spacelike surface.  The past Cauchy
development  of $S$  is the  set  of all  points $p$  such that  every
inextendable future-directed lightlike geodesic ray through $p$ intersects $S$.  
\end{lemma}

\begin{proof}
Let $p\in AdS_3$ be a  point such that every past-directed lightlike
geodesic ray  through $p$  intersects the surface  $S$. Then, every
past-directed   lightlike  geodesic   ray  through   $p$  intersects
(transversally)   the  surface   $S$   at  exactly   one  point   (see
Proposition~\ref{pasdelien}). Hence, the set $C$ of
all the points of $S$ that can be joined from $p$ by a past-directed
lightlike geodesic ray is homeomorphic  to a circle. Therefore, $C$ is
the  boundary of  a closed  disk $D\subset  S$ (recall  that $S$  is a
properly embedded disc,  see Proposition~\ref{p.embedded}). Let $L$ be
the union of all the segments  of lightlike geodesics joining $p$ to a
point of $C$.  The union of $D$ and $L$  is a non-pathological sphere.
By Jordan-Schoenflies theorem, this topological sphere is the boundary
of a  ball $B\subset AdS_3$.  A non-spacelike curve  cannot escape $B$
through   $L$;   as   a  consequence,   every   past-inextendable
non-spacelike curve through $p$ must  escape from $B$ through $D$ ; in
particular, every past-inextendable  non-spacelike curve through $p$
must intersect $S$. Hence, the point $p$ is in $D^+(S)$.
\end{proof}

\begin{rema}
\label{r.range-contained-Cauchy-development}
Since the  surface $\Sigma_0$  is a Cauchy  surface in $M$,  the range
$\cD(\widetilde M)$ of  the developping map $\cD$ must  be contained in
the Cauchy development of the surface $S_0=\cD(\widetilde\Sigma_0)$.  
\end{rema}

We now define another domain, the \emph{black domain}
$E(\partial S_0)$, which, as we will prove later, coincides with
the Cauchy development $D(S)$.

\paragraph{Definition of the set $E(\partial S_0)$.} 
The   set
$$E(\partial    S_0)=\{r\in\SS^3\mid   B_Q(r,q)<0\mbox{    for   every
}q\in\partial S_0\}$$ 
is called the \emph{black domain} of  the curve $\partial S_0$
(explanations on this terminology are provided below).

\begin{rema}
\label{r.black-domain}
Here  are  a  few  observations   about  the  definition  of  the  set
$E(\partial S_0)$:

\smallskip

\noindent       \emph{(i)}       We       will       prove       below
(Proposition~\ref{c.black-domain-in-AdS}))   that  the   black  domain
$E(\partial S_0)$  (which is defined above  as a subset  of the sphere
$\SS^3$) is  actually contained in  the anti-de Sitter  space $\AdS_3$.
Moreover,  we will  prove that,  for a  suitable choice  of  the point
$p_0$,  the set  $E(\partial S_0)$  is contained  in the  affine domain
$\cA_{p_0}$ (Proposition~\ref{p.black-domain-in-affine-domain}).  

\smallskip

\noindent  \emph{(ii)}  Consider  a  point $p_0\in  AdS_3$  such  that
$E(\partial   S_0)$  is   contained  in   $\cA_{p_0}$.    According  to
Proposition~\ref{p.lightcone}, the  set $E(\partial S_0)$  is made of
the points $r\in\cA_{p_0}$  such that there does not  exist any causal
curve joining $r$ to the  curve $\partial S_0$ within $\cA_{p_0}$.  In
other  words, $E(\partial  S_0)$ is  the set  of ``all  the  points of
$\cA_{p_0}$ that cannot be seen  from any point of the curve $\partial
S_0$''.  This is the reason why we call $E(\partial S_0)$
the \emph{black domain} of the curve~$\partial S_0$.  

\smallskip

\noindent \emph{(iii)} The black domain $E(\partial S_0)$ is clearly a
convex subset  of $\SS^3$ (by  construction, it is an  intersection of
convex  subsets  of  $\SS^3$).  In particular,  $E(\partial  S_0)$  is
connected.  

\smallskip

\noindent \emph{(iv)} Here is a nice way to visualize $E(\partial
S_0)$. Consider a point  $p_0\in\AdS_3$ such that $E(\partial S_0)$ is
contained     in     the     affine     domain     $\cA_{p_0}$     (see
Proposition~\ref{p.black-domain-in-affine-domain}).       Using     the
diffeomorphism    $\Phi_{p_0}$,   we    can    identify   $\cA_{p_0}$,
$\partial\cA_{p_0}$,  $\partial  S_0$,  $E(\partial  S_0)$  with  some
subsets of  $\RR^3$ (in particular,  $\partial\cA_{p_0}$ is identified
with   the   hyperboloid   of  equation   $(x^2+y^2-z^2=1)$).    Given
$q\in\partial S_0$, the  set $T_q=\{r\in\cA_p\mid B_Q(q,r)=0\}$ is the
affine  plane   of  $\RR^3$  which  is  tangent   to  the  hyperboloid
$\partial\cA_{p_0}$  at  $q$.  If we define the  set  $E_q=\{r\in\cA_p\mid
B_Q(q,r)<0\}$   as the connected   component   of
$\RR^3\setminus    T_q$ containing at least one point of 
$\partial S_0$,  $\partial S_0$
is contained in the closure of $E_q$, and  the set $E(\partial S_0)$ is
the  intersection  over  all  $q\in\partial  S_0$,  of  the  $E_q$'s.

\smallskip

\noindent \emph{(v)} Let $r$ be  a point on the boundary (in $\AdS_3$)
of $E(\partial S_0)$.  The definition of the set $E(\partial S_0)$ and
the  compactness  of  the  curve  $\partial  S_0$  imply  that  we  have
$B_Q(r,q)=0$ for some point $q$ on the curve $\partial S_0$. Hence, by
Remark~\ref{r.lightlike-geodesic},  there exists a  lightlike geodesic
$\gamma$ passing 
through $r$, such that  one of the two ends of $\gamma$  is a point of
the curve $\partial S_0$.  
\end{rema}

\begin{prop}
\label{p.surface-in-black-domain}
The surface $S_0$ is contained in $E(\partial S_0)$. 
\end{prop}

\begin{proof}
Let $p$ be a point in $S_{0}$. By Proposition \ref{p.curve-in-affine-domain},
 the surface $\overline{S_0}$ is contained
in the affine domain $\cA_p\cup\partial\cA_p$.  
By Proposition~\ref{p.lightcone}, if for some $q$ in
$\partial S_0$ we have $B_Q(p,q)\geq 0$, there is a causal curve
in $\cA_p$ joining $p$ to $q$. But such a curve cannot exist according to
Proposition \ref{p.ajout}.
The proposition follows.
\end{proof}

\begin{prop}
\label{p.B dans E}
The black domain $E(\partial S_0)$ contains the Cauchy
development $D(S_0)$. 
\end{prop}

\begin{proof}
Assume the contrary.  Since $D(S_0)$ and $E(\partial S_0)$ have a non-empty
intersection  (the surface  $S_0$  is contained  in  both $D(S_0)$  and
$E(\partial  S_0)$), and  since $D(S_0)$  is connected,  $D(S_0)$ must
contain  some point  $r$  of  the boundary  of  $E(\partial S_0)$.  By
item~(v)  of  Remark~\ref{r.black-domain},  there exists  a  lightlike
geodesic $\gamma$  passing through $r$, such  that one of  the ends of
$\gamma$ is a point $q$ on the curve $\partial
S_0$. Since $r$
is  in $D(S_0)$, the  lightlike geodesic  $\gamma$ must  intersect the
surface $S_0$.  But, this is impossible according to 
Proposition \ref{pasdelien}.
\end{proof}

\begin{coro}
\label{l.no-timelike-geodesic-in-Cauchy-dev}
The black domain $E(\partial S_0)$ and 
the Cauchy development $D(S_0)$ do not contain any timelike
geodesic. 
\end{coro}

\begin{proof}
Let $\gamma$ be a timelike geodesic. Recall that
$\gamma$ is a closed  geodesic. Consider all future oriented lightlike
geodesic rays
starting from a point of $\gamma$: the union of their future extremities
covers the whole $\partial \AdS_3$, in particular, it contains $\partial S_0$.
It follows that $\gamma$ cannot be contained in the black domain 
$E(\partial S_0)$. Therefore, the corollary follows from
Proposition \ref{p.B dans E}.
\end{proof}

\begin{prop}
\label{p.one-to-one}
The developping map $\cD:\widetilde M\rightarrow AdS_3$ is one-to-one.  
\end{prop}

\begin{proof}
Consider the lifting $\tau: \widetilde{M} \rightarrow \mathbb R$
of any time function on $M$. Select 
any timelike geodesic $\Delta_0$ of $AdS_3$. According to the
corollary~\ref{l.no-timelike-geodesic-in-Cauchy-dev}, the intersection between
$\Delta_0$ and $E(\partial S_0)$ is a subarc $I \approx \mathbb R$ 
(it is connected since $E(\partial S_0)$ is convex). 
Every level set of $\tau$ is the  lift of a Cauchy surface of $M$. So,
by Proposition~\ref{p.embedded} and Remark~\ref{r.other-surface}, for
every $t$ in $\RR$, the image of $\tau^{-1}(t)$ under $\cD$ is a spacelike
surface  that  intersects  $\Delta_0$   at  one  and  only  one  point
$d(t)$. Clearly, $d$ is a strictly increasing function, hence, 
it is injective. Therefore, for any $p$  and $q$ in  $\widetilde M$,
if $\cD(p)=\cD(q)$, then $\tau(p)=\tau(q)$: $p$ and $q$ belongs to the same
spacelike  level   set  of  $\tau$.   According  to  (the   proof  of)
proposition~\ref{p.embedded}, the restriction  of $\cD$ to every level
of $\tau$ is injective. Hence, $p=q$.
\end{proof}

\begin{prop}
\label{p.proper-discontinuity}
The holonomy group $\Gamma=\rho(\pi_1(M))$ acts freely, and properly
discontinuously  on the  Cauchy  development $D(S_0)$  of the  surface
$S_0$. 
\end{prop}

\begin{proof}
First  note   that  the  group  $\Gamma$  acts   freely  and  properly
discontinuously   on   the   surface   $S_0=\cD(\wt\Sigma_0)$   (since
$\cD:\wt\Sigma_0\rightarrow AdS_3$ is a proper embedding). 

Suppose that the group $\Gamma$ does not act freely on the future Cauchy 
development $D^+(S_0)$.  Then, there exists an element $\gamma$ of
$\Gamma$ which fixes a point $p$  of $D^+(S_0)$. Then, as in the proof
of Lemma~\ref{l.def-alternative}, we consider the set $C$ of all the 
points  of $S_0$  that can  be joined  from $p$  by  a past-directed
lightlike geodesic ray.  The set $C$ is homeomorphic  to a circle, and
thus, it is the boundary of a closed disc $D\subset S_0$. The disc $D$
must  be  invariant  under   $\gamma$  (since  the  surface  $S_0$  is
$\Gamma$-invariant, and since $\gamma$ fixes the point $p$). Hence, by
Brouwer's  theorem,  $\gamma$ fixes  a  point  in  $D$.  In  particular,
$\gamma$  fixes a  point  in  $S_0$. This  contradicts  the fact  that
$\Gamma$ acts  freely on  $S_0$.  Hence, $\Gamma$  must act  freely on
$D^+(S_0)$.  The same  arguments  show that  $\Gamma$  acts freely  on
$D^-(S_0)$. 

Now,  let $K$ be  a compact  subset contained  in $D^+(S_0)$.   All the
points of intersection of the past-directed lightlike geodesic rays emanating 
from the points  of $K$ with the surface $S_0$  belong to some compact
subset  $K'$  of the  surface  $S_0$.  Since  $\Gamma$ maps  lightlike
geodesic    rays    to lightlike    geodesic    rays,    the    set
$\{\gamma\in\Gamma\mid \gamma  K\cap K\neq\emptyset\}$ is  contained in
the        set       $\{\gamma\in\Gamma\mid        \gamma       K'\cap
K'\neq\emptyset\}$. Hence, the proper discontinuity of the action of
$\Gamma$ on  $D^+(S_0)$ follows from  the proper discontinuity  of the
action on  $S_0$. The same  arguments show  that $\Gamma$
acts properly discontinuously on $D^-(S_0)$.
\end{proof}

\begin{prop}
\label{p.quotient}
The  spacetime  $M$ is  isometric  to  the quotient  $\Gamma\backslash
D(S_0)$ (the isometry being induced by the developping map $\cD$). 
\end{prop}

\begin{proof}
By      Proposition~\ref{p.proper-discontinuity},     the     quotient
$\Gamma\backslash D(S_0)$ is a manifold (which is automatically a globally
hyperbolic,  since  it  is  the  quotient of  the  Cauchy  development
$D(S_0)$).   By  Remark~\ref{r.range-contained-Cauchy-development}  and
Proposition~\ref{p.one-to-one}, the  developping map $\cD$  induces an
isometric embedding of $M$ in $\Gamma\backslash D(S_0)$.  Since $M$ is
assumed  to  be  maximal  as  a  globally  hyperbolic  manifold,  this
embedding must be onto.  
\end{proof}

According to  Proposition~\ref{p.quotient}, constructing a  surface in
$M$ with some specified geometrical properties amounts to constructing
a $\Gamma$-invariant  surface in $D(S_0)$. In particular,  we will use
the following remark several times: 
 
\begin{rema}
If  $S$ is  a  $\Gamma$-invariant spacelike  surface  contained in  the
Cauchy development  $D(S_0)$, then $\Gamma\backslash S$  is a Cauchy  surface in
$M=\Gamma\backslash D(S_0)$. Indeed,  $\Gamma\backslash S$ is a  spacelike compact surface
in $M=\Gamma\backslash D(S_0)$, and every compact  spacelike surface in $M$ is a
Cauchy surface. 
\end{rema}

\section{Proof of Theorem~\ref{t.CMC-foliation} in the case $g\geq 2$}
\label{s.higher-genus}

We  have   to  prove  that  $M$   admits  a  CMC   time  function.  
In this section, we  give the proof in the case $g  \geq 2$; the proof
in  the other case  $g =1$  (see Remark  \ref{r.sphere-impossible}) is
completely different and will be achieved in 
section~\ref{s.g=1}.

In
subsection~\ref{s.reduction}, we will explain why in the case $g\geq 2$,
this problem reduces to the
proof  of   the  existence   of  a  pair   of  barriers  in   $M$.

In subsection~\ref{s.timelike}, we prove that when $\Sigma$ has higher genus,
then the compactified surface $\overline{S_0}$ is strictly achronal.
In subsection~\ref{ss.black-domain}, we study the 
intersection $C_0$ of $\AdS_3$ with the convex hull of the
curve $\partial S_0$. In particular,  we prove that $C_0$ is contained
in the Cauchy development $D(S_0)$, so that we may 
consider  the  projection  $\Gamma\backslash  C_0$ of  $C_0$  in
$\Gamma\backslash            D(S_0)\simeq            M$.    We also complete
the study in the previous section above by proving, for example, that the
Cauchy development and the black domain coincide\footnote{This last statement
remains true in the case $g=1$,  but the proof is quite different than
those of the case $g \geq 2$.}. 

In
subsection~\ref{ss.quasi-barriers}, we define 
the  notion of  convexity  and concavity  for  spacelike surfaces  in
$\AdS_3$, and  we prove that the  boundary of $C(S_0)$  in $\AdS_3$ is
the union of two  disjoint spacelike topological surfaces $S_0^-$ and
$S_0^+$,   repectively   convex    and   concave.    The   projections
$\Sigma_0^-=\Gamma\backslash S_0^-$ and $\Sigma_0^+=\Gamma\backslash 
S_0^+$  of  these surfaces  in  $\Gamma\backslash  D(S_0)\simeq M$  is
``almost a  pair of barriers''.   There are still two  small problems:
in general,  the surfaces  $\Sigma_0^-$ and $\Sigma_0^+$  have totally
geodesic  regions  (whereas,  for  barriers,  we  need  surfaces  with
positive and negative mean curvature),  and in general, these are only
topological 
surfaces (whereas, for barriers, we need surfaces of class $C^2$). The
purpose of subsections~\ref{s.uniformisation} and~\ref{s.smoothing} is to
approximate the surfaces 
$\Sigma_0^-$ and $\Sigma_0^+$ by a true pair of barriers.

\subsection{Reduction  of Theorem~\ref{t.CMC-foliation} to the existence
  of a pair of barriers} 
\label{s.reduction}

V.  Moncrief has  proved that the  solutions of  the vacuum  Einstein
equation in dimension $2+1$ with a compact Cauchy surface can be
described  as the  orbits of  a non-autonomous  hamiltonian flow  on a
finite-dimensional   space  (namely  the   cotangent  bundle   of  the
Teichm\"uller space of the Cauchy 
surface).  Using this hamiltonian flow, L.  Andersson, Moncrief and A.
Tromba   have   obtained   the  following   theorem   (\cite[corollary
7]{AndMonTro}): 

\begin{theo}[Andersson, Moncrief, Tromba]
\label{t.Moncrief}
Let $N$ be a 3-dimensional maximal globally hyperbolic spacetime, with
constant curvature,  and with closed  Cauchy surfaces of  genus $g\geq
2$. If  $N$ admits  a CMC Cauchy  surface, then  it admits a  CMC time
function. 
\end{theo}

Thanks     to      Theorem~\ref{t.Moncrief},     the     proof     of
Theorem~\ref{t.CMC-foliation} is reduced to the proof of the existence
of a CMC Cauchy surface.  The existence of CMC surfaces, in particular
the existence of  surfaces with zero mean curvature,  has been studied
in many contexts. The problem  usually splits into two disjoint steps~:
a geometrical  step which consists in constructing  some surfaces with
(non-constant)   negative   and   positive   mean   curvature   called
\emph{barriers}, and an analytical  step which consists in solving the
appropriate PDE  to prove  the existence of  a surface with  zero mean
curvature  assuming the existence  of barriers.   In our  context, the
needed    statement    for    the    second    step    is    due    to
C.~Gerhardt (see \cite[Theorem 6.1]{Ger}\footnote{The result proved by
  Gerhardt is  actually more general  than the statement that  we give
  below.}): 
 
\begin{defi}
A \emph{pair of  barriers} in a  three-dimensional globally
hyperbolic Lorentzian manifold $N$ is a pair of disjoint Cauchy
surfaces $\Sigma^-$  and $\Sigma^+$ in $N$,  such that $\Sigma^+$
is in the future of $\Sigma^-$,  the supremum of the mean curvature of
$\Sigma^-$  is negative,  and the  infimum  of the  mean curvature  of
$\Sigma^+$ is positive.  
\end{defi}

\begin{theo}[Gerhardt] 
\label{t.Gerhardt}
Let  $N$   be  a  three-dimensional   globally  hyperbolic  Lorentzian
manifold, with  compact Cauchy surfaces.   Assume that there  exists a
pair of barriers in $N$.  Then, $N$ admits a Cauchy surface with
zero mean curvature in $N$ (\emph{i.e.}, a maximal Cauchy surface).
\end{theo}

Using the results of
Andersson-Moncrief-Tromba and Gerhardt stated  above, the proof of our
main  theorem reduces  to the  proof  of the  existence of  a pair  of
barriers in $M$.  

\subsection{Strict achronality}
\label{s.timelike}

\begin{prop}
\label{p.curve-spacelike}
The topological surface $\overline{S_0}$ is spacelike.  
\end{prop}

\begin{rema}
This  Proposition is  false  without the  assumption  that the  Cauchy
surface $\Sigma_0$ has genus $g\geq 2$, see Remark~\ref{r.not-spacelike}. 
\end{rema}

\begin{proof}
We already  know that $\overline{S_0}$  is non-timelike,  and that
$S_0$ is  spacelike. Hence, $\overline{S_0}$ is spacelike  if and only
if the curve $\partial S_0$ does not contain any non-trivial lightlike
arc.  Therefore, $\overline{S_0}$ is  spacelike if  and only  if $\partial
S_0$ does not  contain any non-trivial arc of some leaf of the left
or the right ruling of $\partial\AdS_3$.  

Let us denote  by $\RR\PP^1_L$ (resp.  $\RR\PP^1_R$) the  space of the
leaves  of the  left (resp.  right) ruling  of  $\partial\AdS_3$.  We
recall that the action of  the holonomy $\rho$ on $\RR\PP^1_L$ reduces
to the action of $\rho_R$ (since, $\rho_L$ preserves individually each
cicrle  of  the left  ruling).  Similarly,  the  action of  $\rho$  on
$\RR\PP^1_R$ reduces to the action of $\rho_L$. 

\begin{lemma}
\label{l.action-minimale}
The actions  of the representations $\rho_L$  and $\rho_R$ respectively
on $\RR\PP^1_R$ and $\RR\PP^1_L$ are minimal. 
\end{lemma}

\begin{proof}
Let $p$ be  a point of the surface $S_0$,  and $n$ the future-pointing
unitary normal  vector of  $S_0$ at  $p$. If $v$  is a  unitary vector
tangent to  $S_0$ at  $p$, then $n+v$  is a future  pointing lightlike
vector.   The  lightlike geodesic  directed  by  $n+v$  is tangent  to
$\partial       AdS_3$       at       two       antipodal       points
(Remark~\ref{r.type-geodesics-Klein}).  These two antipodal points
lie   on  the   same   leaf   of  the   right   ruling;  denote   by
$R_{[\lambda:\mu]}$            this           leaf            (with
$[\lambda:\mu]\in\RR\PP^1_L$). 
The map $(p,v)\rightarrow (p,R_{[\lambda:\mu]})$ identifies the unitary
tangent bundle of the surface $\Sigma_0$ with the flat $\RR\PP^1$
bundle over  $\Sigma_0$ given by $\pi_1(\Sigma_0)\backslash(S_0\times\RR\PP^1)$
where $\gamma\in\pi_1(M)=\pi_1(\Sigma_0)$ acts by
$\gamma.(p,[\lambda:\mu])=(\rho(\gamma)(p),\rho_L(\gamma)([\lambda:\mu]))$.
Hence, the  Euler class  of the representation  $\rho_L$ is  the Euler
class of the unitary tangent bundle of $\Sigma_0$. By a theorem of
Goldman (see \cite{Gol})\footnote{Here, we use the fact that the genus
of $\Sigma_0$ is at least $2$.}, this implies $\rho_L(\pi_1(M))$ is a
cocompact Fuchsian subgroup of $SL(2,\RR)\times Id\simeq SL(2,\RR)$. 
In particular, the action of $\rho_L$ on $\RR\PP^1_R$ is minimal.  
\end{proof}

\textit{End of the proof of Proposition~\ref{p.curve-spacelike}.} 
Denote by $U$ the open subset  of $\partial S_0$, defined as the union
of the interiors of all the non-trivial arcs of leaves of left ruling
contained in  $\partial S_0$.  Note that the  holonomy $\rho$ preserves
the open  set $U$. Now, let  $U_R\subset\RR\PP^1_R$ be the  set of all
leaves of the 
right  ruling that  intersect $U$.  Then $U_R$  is an  open  subset of
$\RR\PP^1_R$ which is preserved by $\rho_L$. Hence, $U_R$ is either empty
or equal to $\RR\PP^1_R$. But the equality $U_R=\RR\PP^1_R$
would imply that $\partial S_0$ is  a leaf of the left ruling, which
is  impossible by  Proposition~\ref{p.curve-in-affine-domain}.  Hence,
$U_R$ is empty, \emph{i.e.} the  curve $\partial S_0$ does not contain
any non-trivial arc  of leaf of the left  ruling. Similarly, for the
right ruling. This completes the proof. 
\end{proof}

\begin{rema}
\label{r.ensemble-limite}
On the one hand,  Proposition~\ref{p.embedded} implies that the action
of $\Gamma$ on the surface $S_0$ is free and properly discontinous. On
the other hand,  Lemma~\ref{l.action-minimale} implies that the action
of $\Gamma$ on $\partial S_0$  is minimal. As a consequence, the curve
$\partial  S_0$ is  the limit  set of  the action  of $\Gamma$  on the
surface $S_0$.
\end{rema}

We thus obtain a more powerfull version of Proposition~\ref{p.ajout}:

\begin{coro}
\label{p.achronal}
For           every          $p\in\AdS_3$           such          that
$\overline{S_0}\subset\cA_p\cup\partial\cA_p$, the surface 
$\overline{S_0}$ is a   strictly    achronal     subset    of
$\cA_p\cup\partial\cA_p$  (\emph{i.e.}  a  causal  curve  contained  in
$\cA_p\cup\partial\cA_p$ can not intersect $\overline{S_0}$ at two
distinct points). 
\end{coro}

\subsection{The convex hull of  the curve $\partial  S_0$} 
\label{ss.black-domain}

In  this subsection,  we will  consider the  convex  hull $\Conv(\partial
S_0)$ of the curve $\partial S_0$. The main goal is to
prove  that  the  set  $\Conv(\partial S_0)\setminus\partial  S_0$  is
contained in the Cauchy development of the surface $S_0$. We will also 
prove that the black domain and the Cauchy development coincide.

\paragraph{Definition of the set $C_0$.}
Denote  by $\Conv(\partial S_0)$ the convex hull
in     $\SS^3$     of     the     curve    $\partial     S_0$     (see
subsection~\ref{ss.convex-subsets}), and consider the set 
$$C_0=\Conv(\partial S_0)\cap\AdS_3$$

\begin{prop}
\label{p.convex-hull-in-black-domain}
The  set  $\Conv(\partial S_0)\setminus\partial  S_0$  is contained  in
$E(\partial S_0)$.  
\end{prop}

\begin{proof}
Let  $q$ be  a point  $\Conv(\partial S_0)\setminus\partial  S_0$, and
let $\widehat q$ be any  point in $\pi^{-1}(\{q\})$ (recall that $\pi$
is the radial projection of $\RR^4\setminus\{0\}$ on $\SS^3$). Let $r$ be a point
in   $\partial  S_0$,   and  let   $\widehat  r$   be  any   point  in
$\pi^{-1}(\{r\})$.  We  have to  prove  that  $B_Q(q,r)$ is  negative,
\emph{i.e.} that $B_Q(\hat q,\hat  r)$ is negative. Since $\widehat q$
is  in  $\pi^{-1}(\Conv(\partial  S_0))$,  one can  find  points
$\widehat  q_1,\dots,\widehat q_n\in\pi^{-1}(\partial S_0)$,  and 
positive numbers $\alpha_1,\dots,\alpha_n$, such that
$\alpha_1+\dots+\alpha_n=1$, and such that $\widehat q=\alpha_1\widehat
q_1+\dots+\alpha_n\widehat q_n$. We  denote by $q_1,\dots,q_n$ the projections
of   the   points  $\widehat   q_1,\dots,\widehat   q_n$.   
For   each $i\in\{1,\dots,n\}$, there are two possibities:\\ 
--- either  $q_i=r$,  and  then  we  have  $B_Q(\widehat  q_i,\widehat
r)=B_Q(\widehat  r,\widehat  r)=0$ (since  $\widehat r$  is on  the
quadric $(Q=0)$),\\ 
--- or    $q_i\neq   r$,    and    then   corollary~\ref{p.achronal}
and Proposition~\ref{p.lightcone}  imply that  $B_Q(\widehat q_i,\widehat  r)$ is
negative.\\ 
Moreover, at  least one $q_i$'s  is different from $r$  (otherwise, we
would have $q_1=\dots=q_n=q$, which is absurd since $q$ is not on
$\partial S_0$). Hence, the quantity $B_Q(\widehat q,\widehat
r)=\alpha_1  B_Q(\widehat q_1,\widehat  r)+\dots+\alpha_n B_Q(\widehat
q_n,\widehat r)$ is negative. The proposition follows.
\end{proof}

\begin{lemma}
\label{l.black-domain-and-boundary-AdS}
For every point $q\in\partial\AdS_3$, there exists a point
$r\in\partial S_0$, such that $B_Q(q,r)$ is non-negative. Moreover, if
the point $q$  is not on the curve $\partial S_0$,  then the point $r$
can be choosen such that $B_Q(q,r)$ is positive. 
\end{lemma}

\begin{proof}
Let $q$ be a  point in $\partial\AdS_3$. Denote by $[x_1:x_2:x_3:x_4]$
the coordinates  of $q$  in $\SS^3$.  Remark~\ref{r.curve}  imply that
there exists $x_1',x_2'$ such that the point $r$ of coordinates
$[x_1':x_2',x_3,x_4]$ is on the curve $\partial S_0$.  
The   sign   of   $B_Q(q,r)$   is   the   sign   of   the   expression
$-x_1x_1'-x_2x_2'+x_3^2+x_4^2$  (we  recall  that  only  the  sign  of
$B_Q(q,r)$  is  well-defined,  see  Remark~\ref{r.sign-well-defined}).
Since the points $q$ and $r$ are both on $\partial\AdS_3$, we have
$Q([x_1:x_2:x_3:x_4])=Q([x_1':x_2':x_3:x_4])=0$. Hence, we have
$-x_1x_1'-x_2x_2'+x_3^2+x_4^2=\frac{1}{2}((x_1-x_1')^2+(x_2-x_2')^2)$. 
As a consequence, $B_Q(q,r)$ is non-negative. Moreover, if $q$ is not
on  the  curve $\partial  S_0$,  then  $(x_1,x_2)$  is different  from
$(x_1',x_2')$, and thus, $B_Q(q,r)$ is positive.  
\end{proof}

\begin{coro}
\label{c.black-domain-in-AdS}
The black domain $E(\partial S_0)$ is contained in $\AdS_3$.
\end{coro}

\begin{proof}
Lemma~\ref{l.black-domain-and-boundary-AdS}    says   that   the
intersection  of  $\partial\AdS_3$ with  $E(\partial  S_0)$ is  empty.
Since  $E(\partial S_0)$  is connected,  this implies
that $E(\partial  S_0)$ is either  contained in $\AdS_3$,  or disjoint
from $\AdS_3$. But, the intersection  of $E(\partial S_0)$ with $\AdS_3$ is
non-empty       (by Proposition~\ref{p.convex-hull-in-black-domain}, for example). 
Hence, $E(\partial S_0)$ is contained in $\AdS_3$.  
\end{proof}

\begin{coro}
\label{c.convex-hull}
The  set $\Conv(\partial    S_0)\setminus\partial   S_0$    is    contained   in
$\AdS_3$, i.e., $C_0=\Conv(\partial S_0)\setminus\partial S_0$. 
\end{coro}

\begin{proof}
The corollary follows immediately from Proposition~\ref{p.convex-hull-in-black-domain} and
corollary~\ref{c.black-domain-in-AdS}.
\end{proof}

We  will denote  by $\overline{E(\partial  S_0)}$ the  closure  of the
black domain $E(\partial S_0)$ in $\AdS_3\cup\partial\AdS_3$.

\begin{coro}
\label{c.boundary-black-domain}
The intersection of $\overline{E(\partial S_0)}$ with $\partial\AdS_3$
is the curve $\partial S_0$.
\end{coro}

\begin{proof}
Proposition~\ref{p.convex-hull-in-black-domain}   implies  that  every
point of the curve  $\partial S_0$ is in $\overline{E(\partial S_0)}$.
Conversely,  let $q$  be a  point  in $\partial\AdS_3\setminus\partial
S_0$. According  to Lemma~\ref{l.black-domain-and-boundary-AdS}, there
exists  a  point  $r\in\partial   S_0$  such  that  $B_Q(q,r)>0$.   By
continuity of  the bilinear form  $B_Q$, there exists  a neighbourhood
$U$ of  $q$ in $\SS^3$, such  that $B_Q(q',r)>0$ for  every $q'\in U$.
In  particular, there  exists  a  neighbourhood $U$  of  $q$ which  is
disjoint   from   $E(\partial   S_0)$.    Hence,   $q$   is   not   in
$\overline{E(\partial S_0)}$.  
\end{proof}

\begin{prop}
\label{p.black-domain-in-affine-domain}
There exists a point $p_0\in\AdS_3$ such that $E(\partial S_0)$ is contained
in the affine domain $\cA_{p_0}$.  
\end{prop}

\begin{adde}
If the curve $\partial S_0$ is not flat\footnote{We say that the curve
  $\partial S_0$  is \emph{flat}  if it is  the boundary of  a totally
  geodesic subspace of $\AdS_3$, or equivalently, if it is contained in
  a great $2$-sphere in $\SS^3$.}, then one can choose the point $p_0$
  such    that   $\overline{E(\partial    S_0)}$   is    contained   in
  $\cA_{p_0}\cup\partial\cA_{p_0}$.  
\end{adde}

\begin{lemma}
\label{l.black-domain-in-affine-domain}
For     every    point     $p\in     C(\partial    S_0)=\Conv(\partial
S_0)\setminus\partial  S_0$,  the black  domain  $E(\partial S_0)$  is
disjoint  from  the  totally  geodesic  surface $p^*$  (and  thus,  is
disjoint from the closed surface~$\overline{p^*}$).  
\end{lemma}

\begin{proof}
Let $p$ be  a point in $C(\partial  S_0)$, and $\hat p$ be  a point in
$\RR^4\setminus\{0\}$  such  that $\pi(\hat  p)=p$.  Since  $p$ is  in
$\Conv(\partial   S_0)$,   one   can   find  some   points   $\widehat
p_1,\dots,\widehat    p_n\in   \pi^{-1}(\partial   S_0)$    and   some
positive  numbers  $\alpha_1,\dots,\alpha_n$  such that  $\widehat
p=\alpha_1\widehat p_1+\dots+\alpha_n\widehat p_n$.  
Let $q$  be a point in  $E(\partial S_0)$ and  $\hat q$ be a  point in
$\RR^4\setminus\{0\}$ such that $\pi(\widehat q)=q$.  Since $q$ is in
$E(\partial  S_0)$, the  quantity  $B_Q(\widehat{p_i},\widehat q)$  is
negative for  every $i$. Hence, the  quantity $B_Q(\widehat p,\widehat
q)=\alpha_1  B_Q(\widehat p_1,\widehat  q)+\dots+\alpha_n B_Q(\widehat
p_n,\widehat q)$ is negative. In
particular, the point $q$  is not on the surface $p^*=\{r\in\AdS_3\mid
B_Q(\widehat p,\widehat r)=0\}$.  This proves 
that set $E(\partial S_0)$ is disjoint from the totally geodesic surface $p^*$. 
Since $E(\partial S_0)$  is contained in $\AdS_3$, it  is also disjoint
from the closed surface $\overline{p^*}$.
\end{proof}

\begin{proof}[Proof of Proposition~\ref{p.black-domain-in-affine-domain}]
Let $p_0$ be a point in $C_0$.  By
Lemma~\ref{l.black-domain-in-affine-domain},   $E(\partial  S_0)$  is
disjoint from the totally geodesic surface $p_0^*$.  Since $E(\partial
S_0)$ is connected, this implies that $E(\partial S_0)$ is contained in
one  of the two  connected components  of $\AdS_3\setminus  p_0^*$. By
Proposition~\ref{p.convex-hull-in-black-domain}, the point $p_0$ is in
$E(\partial  S_0)$.  Hence,  $E(\partial  S_0)$  is  contained  in  the
connected  component   of  $\AdS_3\setminus  p_0^*$  containing~$p_0$,
that is, in $\cA_{p_0}$. 
\end{proof}

\begin{proof}[Proof of the addendum]
If $\partial S_0$ is not flat, then the set $C_0$ has
non-empty  interior.  Let  $p_0$  be   a  point  in  the  interior  of
$C(\partial  S_0)$. On  the one  hand,  the set  $E(\partial S_0)$  is
disjoint from the closed surface $\overline{p^*}$ for every
$p\in  C_0$.  On  the other  hand,  the union  of all  the
surfaces $\overline{p^*}$ when $p$ ranges over $C_0$ is a 
neighbourhood   (in  $\AdS_3\cup\partial\AdS_3$)   of  the   surface
$\overline{p_0^*}$.  Hence, 
$E(\partial  S_0)$ is disjoint  from a  neighbourhood of  the surface
$\overline{p_0^*}$.  Hence, $\overline{E(\partial S_0)}$ is disjoint
from     the     surface     $\overline{p_0^*}$.      Moreover,     by
Proposition~\ref{p.convex-hull-in-black-domain}, the point $p_0$ is in 
$\overline{E(\partial S_0)}$.   Therefore, $\overline{E(\partial S_0)}$ is
contained        in       the       connected        component       of
$(\AdS_3\cup\partial\AdS_3)\setminus    \overline{p_0^*}$   containing
$p_0$, \emph{i.e.} is contained in $\cA_{p_0}\cup\partial\cA_{p_0}$.
\end{proof}

From now on, we fix a point $p_0\in\AdS_3$, such that
$\overline{E(\partial        S_0)}$        is       contained        in
$\cA_{p_0}\cup\partial\cA_{p_0}$.  

\begin{prop}
\label{p.B=D}
The black domain $E(\partial S_0)$ coincides with the Cauchy
development $D(S_0)$. 
\end{prop}

\begin{proof}
Proposition \ref{p.B dans E} provides an inclusion.
To  prove  the   other  inclusion,  we  work  in   the  affine  domain
$\cA_{p_0}$.  Let $p$ be a point in $E(\partial
S_0)$.  By  Remark~\ref{r.past-future-surface}, every point  of $\cA_{p_0}$ is
either in the past, or in the future of the surface $S_0$.  We assume,
for example, that  $p$ is in the future of $S_0$.   We will prove that
$p$ is in  $D^+(S_0)$.  For that purpose, we  consider a past-directed
lightlike geodesic ray  $\gamma$ emanating from $p$, and  we denote by
$q$ the past end of $\gamma$.  

\smallskip

\noindent  \emph{Claim.  The  geodesic  ray  $\gamma$  intersects  the
  boundary of $E(\partial S_0)$ at some point $r$ in the past
  of $S_0$.} 

\smallskip

\noindent To prove this claim,  we argue by contradiction. First,
we suppose that the
geodesic ray $\gamma$ is contained in $E(\partial S_0)$. Then,
by                              Proposition~\ref{c.black-domain-in-AdS}
and corollary~\ref{c.boundary-black-domain}, the past end of $\gamma$ must be a
point  $q$   of  the  curve   $\partial  S_0$.   But  then,   we  have
$B_Q(p,q)=0$, and this contradicts the fact that $p$ is in $E(\partial
S_0)$. Now, we  suppose that the geodesic ray  $\gamma$ intersects the
boundary $E(\partial S_0)$ at some point $r$  in the future of
the surface  $S_0$. By item (v)  of Remark~\ref{r.black-domain}, there
exists a lightlike geodesic
ray $\gamma'$ emanating from $r$, such that the end of $\gamma'$ is a
point $q$ of the curve $\partial S_0$. The geodesic ray $\gamma'$ must
be past-directed  from $r$ to $q$, since  $r$ is in the  future of the
surface $S_0$. So, we  have a past-directed lightlike geodesic segment
going from $p$ to $r$, and a past-directed geodesic ray going from $r$
to $q$;  concatenating these two  curves, we obtain a  piecewise $C^1$
causal curve  going from $p$  to $q\in\partial S_0$.  This contradicts
the  fact  that  $p$  is  in  $E(\partial  S_0)$  (see  item  (ii)  of
Remark~\ref{r.black-domain}) and completes the proof of the claim. 

Since the point  $p$ is in the future of the  surface $S_0$, and since
the point $r$ given by the claim  is in the past of the surface $S_0$,
the geodesic ray  $\gamma$ must intersect the  surface $S_0$. So,
we have proved that every past-directed geodesic ray emanating from $p$
intersects the  surface $S_0$. Hence,  the point $p$ is  in $D^+(S_0)$
(Lemma~\ref{l.def-alternative}). This proves that $E(\partial S_0)$ is
contained in $D(S_0)$. 
\end{proof}

\begin{rema}
Proposition~\ref{p.B=D}   implies  in   particular  that   the  Cauchy
development  $D(S_0)$  depends  only  on  the  curve  $\partial  S_0$,
\emph{i.e.} if  $S$ is another  complete spacelike surface  in $\AdS_3$
such that $\partial S=\partial S_0$, then $D(S)=D(S_0)$.  
\end{rema}

\begin{rema}
\label{r.meme-bord}
Let   $\Sigma$    be   any   Cauchy   surface   in    $M$,   and   let
$S:=\cD(\wt\Sigma)$. On 
the one hand, we have $D(S)=D(S_0)=\cD(\widetilde M)$. On the other hand,
Propositions~\ref{c.boundary-black-domain}  and~\ref{p.B=D} imply that
the  curve  $\partial S_0$  is  the  intersection  of the  closure  in
$\AdS_3\cup\partial\AdS_3$          of          $D(S_0)$          with
$\partial\AdS_3$.   Similarly,   the  curve   $\partial   S$  is   the
intersection of the  closure in $\AdS_3\cup\partial\AdS_3$ of $D(S)$
with $\partial\AdS_3$. As a  consequence, we have $\partial S=\partial
S_0$. 
\end{rema}

\begin{rema}
\label{r.ensemble-limite-2}
For every point $p\in D(S_0)=\cD(\widetilde M)$, one can find a
Cauchy surface $\Sigma$ in $M$ such that $p\in\cD(\widetilde\Sigma)$.
By Remark~\ref{r.ensemble-limite} and~\ref{r.meme-bord}, the limit set
of the  action of $\Gamma$ on the surface $S$ is  the curve $\partial
S=\partial S_0$.  As a consequence, the limit set of the action of
$\Gamma$ on $D(S_0)$ is also the curve $\partial S_0$.  
\end{rema}

\begin{figure}[ht]
\centerline{\input{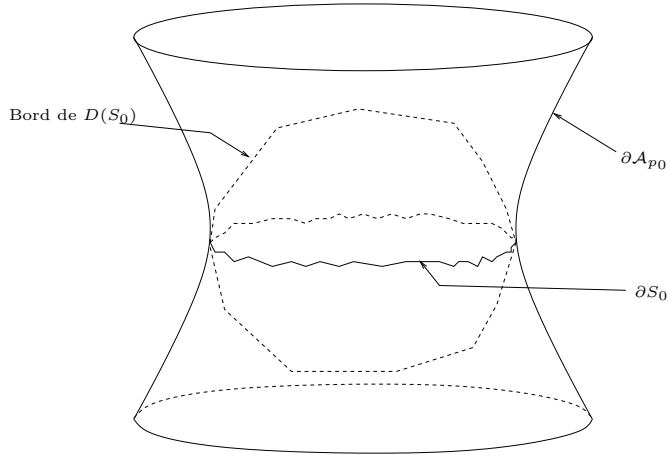}}
\caption{The affine domain $\cA_{p_0}$, the curve $\partial S_0$ 
and the Cauchy development $D(S_0)$.}
\end{figure}

\subsection*{Interlude: proof of Theorem~\ref{t.CMC-foliation} in the 
case where $\partial S_0$ is flat} 

Our strategy  for proving the existence  of a pair of  barriers in $M$
does not work in the particular case where $\partial S_0$ is
flat, mostly  because        the       addendum        of
Proposition~\ref{p.black-domain-in-affine-domain} is false when 
$\partial S_0$  is flat. This is not  a big problem, since  there is a
direct and  very short  proof of Theorem~\ref{t.CMC-foliation}  in this
particular case:

\begin{proof}[Proof  of  Theorem~\ref{t.CMC-foliation}  in  the  case
  where $\partial S_0$ is flat]
Assume  that $\partial  S_0$ is  flat. Then it is  the boundary  of a
totally geodesic  subspace $P_0$  of $\AdS_3$.  This  totally geodesic
subspace is  necessarly spacelike, since  the curve $\partial  S_0$ is
spacelike. By  construction, $P_0$  is contained in  $C_0$;
hence,   it   is  contained   in   the   Cauchy  development   $D(S_0)$
(Proposition~\ref{p.B=D}      and~\ref{p.convex-hull-in-black-domain}).
Moreover, the holonomy group 
$\Gamma=\rho(\pi_1(M))$ preserves $P_0$ (since it preserves the curve
$\partial S_0$).  As a consequence, $\Gamma\backslash P_0$ is a totally
geodesic  compact spacelike  surface in  $\Gamma\backslash D(S_0)\simeq  M$.  In
particular, $\Gamma\backslash P_0$ is  a Cauchy surface with zero mean
curvature in 
$M$.       Applying      Theorem~\ref{t.Moncrief},      we      obtain
Theorem~\ref{t.CMC-foliation}.  
\end{proof}

\paragraph{Assumption.}  From  now  on,   we  assume  that  the  curve
$\partial S_0$ is not flat.

\subsection{A pair of convex/concave topological Cauchy surfaces}
\label{ss.quasi-barriers}

In this subsection, we will first define some notions of convexity and 
concavity for spacelike surfaces  in $M$. The main interesting feature
of this notion for our purpose is the fact that the mean curvature of a
smooth convex (resp.  
concave)     spacelike      surface     is     always     non-positive
(resp.   non-negative).  Then, we  will  exhibit  a  pair of  disjoint
topological  Cauchy surfaces  $(\Sigma_0^-,\Sigma_0^+)$  in $M$,  such
that $\Sigma_0^-$ is convex, $\Sigma_0^+$ is concave, and $\Sigma_0^+$
is in the future of $\Sigma_0^-$.

\subsubsection{Convex and concave surfaces in $\AdS_3$}

Let $S$ be a topological surface in $\cA_{p_0}$, and $q$ be a point of
$S$. A \emph{support plane} of $S$ at $q$ is a ($2$-dimensional) totally
geodesic  subspace\footnote{By a  \emph{totally  geodesic subspace  of
    $\cA_{p_0}$},  we  mean  the  intersection of  a  totally  geodesic
    subspace of $\AdS_3$ with $\cA_{p_0}$. Note that, with this
 definition,   the  degenerated   totally  geodesic   subspaces  of
    $\cA_{p_0}$ are not connected (although their closure in 
    $\cA_{p_0}\cup\partial\cA_{p_0}$ is connected), but this does not play any role in the subsequent.} $P$ of
$\cA_{p_0}$, such that $q\in P$, and  such that $S$ is contained in the
closure of one of the connected components of $\cA_{p_0}\setminus P$.

\begin{rema}
\label{r.spacelike-support-plane}
Let $S$ be  a topological surface in $\cA_{p_0}$.  If $S$ is spacelike
(in  the sense  of definition~\ref{d.spacelike}),  then  every support
plane  of $S$  is spacelike.   Conversely, if  $S$ admits  a spacelike
support plane at every point, then $S$ is spacelike.  
\end{rema}

\begin{rema}
\label{r.spacelike-support-plane-2} 
Let $S$ be a topological surface in $\cA_{p_0}$ and $P$ be a spacelike
support  plane  of   $S$.   Then,  $S$  is  contained   in  the  causal
past\footnote{By causal  past, we mean causal past  in $\cA_{p_0}$} of
$P$,   or   $S$   is   contained    in   the   future   of   $P$   (see
Remark~\ref{r.past-future-plane}).  
\end{rema}

Let $S$ be a topological spacelike surface in $\cA_{p_0}$. We say that
$S$ is  \emph{convex}, if it  admits a support  plane at
each of  its points, and if  it is contained  in the future of  all its
support  planes. We  say that  $S$ is  \emph{concave}, if  it  admits a
support plane a each of its points,  and if it is contained in the past
of all its support planes.  

Now,  let $\Sigma$  be a  topological  spacelike surface  in $M$,
$\widetilde\Sigma$ be  a lift of  $\Sigma$ in $\widetilde M$,  and 
$S=\cD(\widetilde\Sigma)$.  Note that $S$ is a topological spacelike 
surface   contained   in   $\cD(\widetilde   M)\subset\cA_{p_0}$   (see
section~\ref{s.Cauchy-development}).    We   say   that  $\Sigma$   is
\emph{convex} (resp. \emph{concave}) if $S$ is convex (resp. concave).

\begin{prop}
\label{p.link-convexity-curvature}
Let $\Sigma$ be a $C^2$ spacelike surface in $M$. If $\Sigma$ is convex,
then $\Sigma$ has non-positive mean curvature. If $\Sigma$ is concave,
then $\Sigma$ has non-negative mean curvature.  
\end{prop}

\begin{proof}
Let $\widetilde\Sigma$ be a lift of $\Sigma$ in $\widetilde M$, and
let $S=\cD(\widetilde M)$. Assume that $\Sigma$ is convex. Then $S$ is
convex. Hence, for every $q\in  S$, the surface $S$ admits a spacelike
support plane  $P_q$ at $q$, and  is contained in the  future of $P_q$.
By  Lemma~\ref{l.compare-curvature}, the  mean  curvature of  the
surface $S$ at $q$ is smaller  or equal than the mean curvature of the
support plane  $P_q$.  But,  since $P_q$ is  totally geodesic,  it has
zero  mean curvature.  Hence,  the surface  $S$ has  non-positive mean
curvature.  Hence, the surface $\Sigma$ also has non-positive mean
curvature (since the developping map $\cD$ is locally isometric).  
\end{proof}

The notions of convexity and concavity defined above can only help us
in finding  spacelike surfaces with  non-positive (resp. non-positive)
mean curvature.  Yet, to apply Gerhardt's Theorem~\ref{t.Gerhardt}, we
need to find spacelike surfaces with positive (resp. negative) mean
curvature.  This is the reason why we will define below a notion of
\emph{uniformly curved surface} in $M$.

Let $S$  be a topological  surface in $\RR^3$,  and $q$ be a  point on
$S$. We  fix a Euclidean metric  on $\RR^3$. We say  that the surface
$S$ is \emph{more curved than a sphere of radius $R$ at $q$}, if there 
exists a closed Euclidean ball $B$  of radius $R$, such that $q$ is on
the boundary of $B$, and such that $B$ contains a neighbourhood of $q$
in $S$.  

\begin{rema}
\label{r.osculating-quadric}
Assume that the surface $S$ is  $C^2$. Then, $S$ is more curved than a
sphere of radius  $R$ at $q$ if and only if  the osculating quadric of
$S$ at $q$ is an ellipsoid of diameter smaller than $2R$. 
\end{rema}

Consider a topological surface $\Sigma$ in $M$, and a lift $\widetilde \Sigma$
 of $\Sigma$. Let $S=\cD(\widetilde\Sigma)$.  We see $\Sigma$
as a surface in $\RR^3$. Let
$\Delta\subset\widetilde\Sigma$ be a fundamental domain of the covering
$\widetilde\Sigma\rightarrow\Sigma$,  and let $D=\cD(\Delta)$.   
We say that the surface  $\Sigma$ is \emph{uniformly curved}, if there
exists $R\in  (0,+\infty)$ such  that the surface  $S$ is  more curved
than a sphere of radius $R$ at each point of $D$. It is easy to verify 
that this definition depends neither on the choice of the fundamental
domain $\Delta$,  nor on  the choice of the  Euclidean metric  on $\RR^3$
(although  one has  to  change  the constant  $R$,  when changing  the
fundamental domain $\Delta$ or the Euclidean metric on $\RR^3$). 

\begin{prop}
\label{p.link-convexity-curvature-2}
Let $\Sigma$ be a $C^2$ spacelike surface in $M$.  If $\Sigma$
is  convex  and uniformly  curved,  then  $\Sigma$  has negative  mean
curvature. If $\Sigma$ is  concave and uniformly curved, then $\Sigma$
has positive mean curvature.  
\end{prop}

\begin{proof} 
Let $\widetilde\Sigma$ be a lift of $\Sigma$ in $\widetilde M$, and let
$S:=\cD(\widetilde M)$.  Assume that $\Sigma$ is  convex and uniformly
curved. Then,  $S$ is convex. So,  for every $q\in  S$, the surface
$S$ admits a support plane $P_q$ at $q$, and is contained in the future
of $P_q$.  Moreover, since $\Sigma$  is uniformly curved,  the surface
$S$ and the  plane $P_q$ do not have the  same osculating quadric (see
Remark~\ref{r.osculating-quadric}).  
By  Lemma~\ref{l.compare-curvature},   this  implies  that   the  mean
curvature of $S$ at $q$ is strictly smaller than the mean curvature of
the plane $P_q$. Since $P_q$  is totally geodesic, $P_q$ has zero mean
curvature.  Hence,  $S$ has  negative mean curvature.  Therefore, $\Sigma$
also has negative mean curvature. 
\end{proof}

\subsubsection{Boundary  of $\Gamma$-invariant  convex  sets contained in
  $D(S_0)$}

\begin{prop}
\label{p.spacelike-support-plane}
Let $S$  be a topological surface  in $\cA_{p_0}$. Assume  that $S$ is
contained   in   $D(S_0)$,   and   that   the  boundary   of   $S$   in
$\cA_{p_0}\cup\partial\cA_{p_0}$ is equal to the curve $\partial S_0$.
Then every  support plane of $S$ is  spacelike\footnote{Note that, 
in general, the surface $S$ does not admit any support plane.}.  
\end{prop}

\begin{proof}
Using  the diffeomorphism $\Phi_{p_0}$,  we identify  $\cA_{p_0}$ with
the region of $\RR^3$ defined by the inequation $(x^2+y^2-z^2<1)$, and
$\partial\cA_{p_0}$      with     the      one-sheeted     hyperboloid
$(x^2+y^2-z^2=1)$. Let $q$ be a point  on the surface $S$ and $P$ be a
support plane of $S$ at $q$.  The totally geodesic subspace $P$ is the
intersection  of $\cA_{p_0}$  with  an affine  plane  $\widehat P$  of
$\RR^3$.

On the one hand,  since $P$ is a support plane of  $S$, the closure of
$S$  must be  contained in  the  closure of  one of  the two  connected
components  of $\RR^3\setminus\widehat  P$. In  particular,  the curve
$\partial  S_0$ must  be contained  in the  closure of  one of  the two
connected components of $\RR^3\setminus\widehat P$. On the other hand,
$\partial  S_0$ is  a simple  closed  curve  on the  hyperboloid
$\partial\cA_{p_0}$,     which     is     not    null-homotopic     in
$\partial\cA_{p_0}$ (see Remark~\ref{r.curve}).  Consequently: 

\smallskip

\noindent   \emph{Fact    1.  The   support    plane   $P=\widehat
P\cap\cA_{p_0}$ does not contain any affine line of $\RR^3$.} Indeed,
if $\widehat P\cap\cA_{p_0}$ contains  an affine line of $\RR^3$, then
it is  easy to see  $\widehat P\cap\partial\cA_{p_0}$ is  a hyperbola,
and    that    the    two    connected   components    of    $\partial
\cA_{p_0}\setminus\widehat P$  are contractible in $\partial\cA_{p_0}$
(we  recall  that  $\cA_{p_0}$  is  the  region  $(x^2+y^2-z^2<1)$  in
$\RR^3$). Hence,  every curve contained  in the closure of  a connected
component of $\partial \cA_{p_0}\setminus\widehat P$ is null-homotopic
in $\partial\cA_{p_0}$.  

\smallskip

\noindent \emph{Fact 2.   If the plane $\widehat  P$ is tangent
  to  the hyperboloid  $\partial\cA_p$  at some  point  $r$, then  $r$
  belongs to
  the curve $\partial S_0$.}  Indeed, if $\widehat P$ is tangent
  to the hyperboloid $\partial\cA_{p_0}$ at some point $r$, then every
curve contained in  the closure of one of  the two connected components
of $\partial\cA_{p_0}\setminus\widehat P$  which is not null-homotopic
in $\partial\cA_{p_0}$ contains $r$.  

\medskip

Now, we  argue by contradiction:  we assume that the  totally geodesic
plane  $P$ is  not spacelike.   Then,  $P$ is  either timelike  (the
Lorentzian  metric  restricted  to  $P$  has  signature  $(+,-)$),  or
degenerated (the Lorentzian metric  restricted to $P$ is degenerated).
We will show that the two possibilities lead to a contradiction. 

\noindent  --  If  $P$  is  timelike,  then  $P$  contains  timelike
geodesics.  By Remark~\ref{r.type-geodesics-AdS+}, a timelike geodesic
of $\cA_{p_0}$ 
is an affine line of $\RR^3$ which is contained in $\cA_{p_0}$.  Hence,
$P=\widehat P\cap\cA_{p_0}$  contains an affine line  of $\RR^3$. This
is absurd according to Fact 1 above. 

\noindent --  If $P$  is degenerated then  $P$ contains  lightlike and
spacelike geodesic, but does not contain any timelike geodesic. By
Remark~\ref{r.type-geodesics-AdS+},  this  implies  that $\widehat  P$
is tangent to the hyperboloid $\partial\cA_{p_0}$ at 
some point $r$.  According to Fact 2, the point $r$ must belong to 
the curve $\partial  S_0$.  But then, Remark~\ref{r.black-domain} item
(iv) implies that $P$ is disjoint from
$E(\partial S_0)$.  In particular, the point $q$ is not in $E(\partial
S_0)$.   This is  absurd  since,  by hypothesis,  the  surface $S$  is
contained in $E(\partial S_0)=D(S_0)$.   
\end{proof}

\begin{prop}
\label{p.boundary-convex-set}
Let $C$ be a non-empty $\Gamma$-invariant closed\footnote{By such, we mean that
  $C$   is    closed   in    $\AdS_3$,   but   not    necessarily   in
  $\AdS_3\cup\partial\AdS_3$. Actually, a non-empty $\Gamma$-invariant 
  subset of  $AdS_3$ cannot be closed in  $AdS_3\cup\partial AdS_3$}
  convex subset of $\AdS_3$, contained in $D(S_0)$. Then: 

\smallskip

\noindent \emph{(i)}  The boundary of $C$  in $\AdS_3$ is made of
two disjoint $\Gamma$-invariant topological surfaces $S^-$ and $S^+$,
such that $S^-$ is convex, $S^+$ is concave, $C$ is in the future of
$S^-$ and in the past of $S^+$. 

\smallskip

\noindent    \emph{(ii)}    $\Sigma^-:=\Gamma\backslash    S^-$    and
$\Sigma^+:=\Gamma\backslash S^+$ are two disjoint Cauchy surfaces in 
$\Gamma\backslash  D(S_0)\simeq M$.   Moreover, $\Sigma^-$  is convex,
$\Sigma^+$   is  concave,  and   $\Sigma^+$  is   in  the   future  of
$\Sigma^-$. Of course, the boundary  of the set $\Gamma\backslash C$ in
$M$ is the union of the surfaces $\Sigma^-$ and $\Sigma^+$.
\end{prop}

\begin{proof}
Since $C$ is  contained in $D(S_0)$, it is also  contained in the affine
domain $\cA_{p_0}$.
We  denote by  $\partial C$  the boundary  of $C$  in  $\cA_{p_0}$,  by $\overline{C}$ the closure of $C$ in $\cA_{p_0}\cup\partial
\cA_{p_0}$, and by $\overline{\partial C}$ the boundary of
$\overline{C}$  in  $\cA_{p_0}\cup\partial\cA_{p_0}$.   Of course,  we
have $\partial C=\overline{\partial C}\cap\cA_{p_0}=\overline{\partial
C}\setminus\partial\cA_{p_0}$.  

The set $\overline{C}$ is a compact convex subset of
$\cA_{p_0}\cup\partial\cA_{p_0}$. So,  the diffeomorphism $\Phi_{p_0}$
maps $\overline{C}$ to a compact convex subset of $\RR^3$. 
 Hence,  $\overline{\partial   C}$  is  a  $\Gamma$-invariant
topological sphere. 
We have to understand the intersection of $\overline{\partial C}$
with  $\partial\cA_{p_0}$.  On  the one  hand, by  hypothesis,  $C$ is
contained   in   $D(S_0)$;  hence,   $\overline{C}$   is  contained   in
$\overline{D(S_0)}$.   The  intersection  of $\overline{D(S_0)}$  with
$\partial\cA_{p_0}$  is  equal  to   the  curve  $\partial  S_0$  (see
Propositions~\ref{c.boundary-black-domain} and~\ref{p.B=D}).  Hence, the 
intersection  of $\overline{\partial  C}$ with  $\partial\cA_{p_0}$ is
contained in the curve $\partial S_0$.  On the other hand, $C$ is a non-empty
$\Gamma$-invariant subset of $D(S_0)$. Hence, the closure of $C$ must
contain the curve $\partial S_0$ (since this curve is the limit set of
the action of $\Gamma$ on $D(S_0)$). As a consequence,
we have $\overline{\partial C}\cap\partial\AdS_3=\partial S_0$.  
  
We have proved that $\partial
C=\overline{\partial      C}\setminus\partial\cA_{p_0}$      is      a
$\Gamma$-invariant  topological   sphere  minus  the  $\Gamma$-invariant
Jordan curve $\partial S_0$.  Hence,  $\partial C$ is the union of two
disjoint $\Gamma$-invariant topological discs $S^-$ and $S^+$, such that
$\partial S^-=\partial S^+=\partial S_0$.   Since the surfaces $S^-$ and
$S^+$  are contained  in the  boundary of  a convex  set, they  admit a
support    plane    at   each    of    their    points.   Hence,    by
Proposition~\ref{p.spacelike-support-plane}                         and
Remark~\ref{r.spacelike-support-plane},  the surfaces  $S^-$  and $S^+$
are  spacelike.   Since  $S^-$  is  a spacelike  disc  with  $\partial
S^-=\partial  S_0$,  it   separates  $\cA_{p_0}$  into  two  connected
components: the  past and  the future  of $S^-$. The  set $C$  must be
contained in one of these  two connected components, so $C$ is contained
either  in  the  past  or  in  the  future  of  $S^-$.  Similarly,  for
$S^+$. Moreover, $C$ can not be in the future (resp. the past) of both
$S^-$ and $S^+$.  So, up to exchanging $S^-$ and $S^+$,  the set $C$ is
in the future of $S^-$ and  in the past of $S^+$. In particular, $S^+$
is in the  future of $S^-$. Since  $C$ is in the future  of $S^-$, the
surface $S^-$  must be in  the future of  each of its  support planes.
Hence, the surface $S^-$ is  convex. Similar arguments show that $S^+$
is concave.  This completes the proof of $(i)$. 

Now, since  $S^-$ and $S^+$ are  $\Gamma$-invariant spacelike surfaces
in $D(S_0)$, their projections $\Sigma^-:=\Gamma\backslash S^-$ and
$\Sigma^+:=\Gamma\backslash    S^+$    are    Cauchy    surfaces    in
$\Gamma\backslash D(S_0)\simeq M$ (recall that every compact spacelike
surface in $M$ is a Cauchy 
surface). Of course, $\Sigma^+$ is in the future of $\Sigma^-$, since
$S^+$ is in the future of $S^-$.  Finally, the convexity of $\Sigma^-$
and  the  concavity of  $\Sigma^+$  follow,  by  definition, from  the
convexity of $S^-$ and the concavity of $S^-$. 
\end{proof}

\subsubsection{Definition of the topological Cauchy surfaces $\Sigma_0^-$
  and $\Sigma_0^+$}

The  set  $C(\partial  S_0)=\Conv(\partial S_0)\setminus\partial  S_0$
satisfies  the hypothesis  of Proposition~\ref{p.boundary-convex-set}.
Hence, the  boundary in  $\AdS_3$ of $C(\partial  S_0)$ is made  of two
disjoint $\Gamma$-invariant spacelike topological surfaces $S_0^-$ and
$S_0^+$, such that $S_0^-$ is  convex, $S_0^+$ is concave, and $S_0^+$
is   in    the   future   of   $S_0^-$.     
Moreover,  the   surfaces
$\Sigma_0^-:=\Gamma\backslash S_0^-$ and $\Sigma_0^+:=\Gamma\backslash
S_0^+$   are    two   disjoint   topological    Cauchy   surfaces   in
$\Gamma\backslash D(S_0)\simeq M$, such that 
$\Sigma_0^-$ is  convex, $\Sigma_0^+$ is concave,  and $\Sigma_0^+$ is
in the future of $\Sigma_0^-$.  

\begin{defi}
A pair $(S^-,S^+)$ of 
disjoint $\Gamma$-invariant spacelike topological surfaces in $AdS_3$
such that $S_0^-$ is  convex, $S^+$ is concave, and $S^+$
is   in    the   future   of   $S^-$ is called a \emph{convex trap}.

Similarly, a pair $(\Sigma^-,\Sigma^+)$ of 
disjoint spacelike topological surfaces in $M$
such that $\Sigma^-$ is  convex, $\Sigma^+$ is concave, and $\Sigma^+$
is   in    the   future   of   $\Sigma^-$ is called a \emph{convex trap}.

In both circonstances, a convex trap is uniformly curved if the boundary surfaces 
$S^{-}$, $S^{+}$ (or $\Sigma^{-}$, $\Sigma^{+}$) ) are uniformly curved.
The convex trap is smooth if the boundary surfaces are smooth.
\end{defi}

\subsection{A  pair of  uniformly curved  convex/concave topological
  Cauchy surfaces} 
\label{s.uniformisation}

Our  goal  is  to find  a  pair of  barriers  in $M$.  By
Proposition~\ref{p.link-convexity-curvature-2},  this   goal  will  be
achieved if we find a smooth uniformely curved convex trap. 
For the moment, the convex trap 
$(\Sigma_{0}^-, \Sigma_{0}^+)$ is not smooth, and not uniformely curved.
The purpose of  this subsection is to
prove the following proposition: 

\begin{prop}
\label{p.uniformisation}
Arbitrarily  close to  $\Sigma_{0}^-$ (resp. $\Sigma_0^+$), there
exists     a     topological     Cauchy    surfaces     $\Sigma_{1}^-$
(resp. $\Sigma_{1}^+$), which is  convex (resp. concave) and uniformly
curved. 
\end{prop}

The  idea of  the  proof of  Proposition~\ref{p.uniformisation} is  to
replace  the  convex set  $C_0=C(\partial  S_0)$  by its  ``Lorentzian
$\varepsilon$-neighbourhood''.    This  idea  comes   from  Riemannian
geometry.      Indeed,      it     is     well-known      that     the
$\varepsilon$-neighbourhood of a convex subset of the hyperbolic space 
$\HH^n$ is  uniformly convex. We  will prove that a  similar phenomenon
occurs in $\AdS_3$ (although technical problems appear). 

The  \emph{length}   of  a  $C^1$ causal  curve
$\gamma:[0,1]\rightarrow      \AdS_{3}$     is     $l(\gamma)=\int_0^1
(-g(\dot\gamma(t),\dot\gamma(t)))^{1/2}dt$,    where   $g$    is   the
Lorentzian  metric of  $\AdS_{3}$.  Given an  achronal  subset $E$  of
$\cA_{p_0}$ and  a point $p$  in $\cA_{p_0}$, the  \emph{distance from
  $p$ to $E$} is the supremum of the lengths of
all  the $C^1$ causal  curves joining  $p$ to  $E$ in  $\cA_{p_0}$ (if
there  is no such  curve, then  the distance  from $p$  to $E$  is not
defined)\footnote{The same  definition work in the case  where the set
  $E$ is not achronal. But then, the distance from $p$ to $E$ might be
  positive even  if $p\in E$~!}  The  distance from $p$  to $E$, when
finite, is lower semi-continuous  in $p$.  Moreover, the distance from
$p$ to $E$ is continuous in $p$, when $p$ is in the Cauchy development
of $E$ (see, for instance,~\cite[page 215]{Haw}).  

Given  an achronal  subset $E$  of $\cA_{p_0}$  and  $\epsilon>0$, the
\emph{$\epsilon$-future} of  $E$ is the  set made of the  points $p\in
\cA_{p_0}$, such that $p$ is in the future of $E$ and such that the
distance from $p$ to $E$ is at most
$\epsilon$.   We define similarly  the \emph{$\epsilon$-past}  of $E$.
We   denote    by   $I_\epsilon^-(E)$   and    $I_\epsilon^+(E)$   the
$\epsilon$-past and the $\epsilon$-future of the set $E$.  

\begin{lemma}
\label{l.inclus-dev-Cauchy}
There exists $\epsilon>0$ such that the $\epsilon$-past and the
$\epsilon$-future of the surface $S_0^+$ are contained in~$D(S_{0})$.  
\end{lemma}

\begin{proof}
Since the  set $D(S_{0})$ is  a neighbourhood of the  surface $S_0^+$,
and since the surface $\Sigma_0^+=\Gamma\backslash S_0^+$ is 
compact, one  can find  a $\Gamma$-invariant neighbourhood  $U_0^-$ of
the surface $S_0^+$, such that  $U_0^+$ is contained in $D(S_{0})$, and
such that $\Gamma\backslash U_0^+$ is compact.  

\smallskip

\noindent \emph{Claim. There exists $\epsilon>0$ such that the 
distance  from any  point $p\notin  U_0^+$ to  the surface  $S_0^+$ is
bigger than~$\epsilon$.} 

\smallskip

\noindent By  contradiction, suppose that, for  every $n\in\NN$, there
exists a point $x_n\in 
\cA_{p_0}\setminus  U_0^+$ such that  the distance  from $x_n$  to the
surface $S_0^+$ is less than $1/n$. Then, for each $n$, we 
consider a causal curve $\gamma_n$  joining the point $x_n$  to the
surface $S_0^+$. This curve $\gamma_n$ must intersect the boundary of
$U_0^+$; let $z_n$ be a point in $\gamma_n\cap \partial U_0^+$.  Since
$z_n$ is on a causal curve joining $x_n$ to the surface $S_0^+$, the
distance from $z_n$ to $S_0^+$ must be smaller than $1/n$. Now, recall
that $\Gamma\backslash U_0^+$ is compact.  Hence, up to replacing each
$z_n$ by 
its image under some element of $\Gamma$,  we may  assume that  all the  $z_n$'s are  in  a compact
subset of  the boundary of $U_0^+$.   Then, we consider  a limit point
$z$ of the sequence $(z_n)_{n\in\NN}$.  By lower semi-continuity of the
distance, the  distance from  $z$ to the  surface $S_0^+$ is  equal to
zero  (note that the  distance from  $z$  to the  surface $S_0^+$  is
well-defined,  since every  point of  $\cA_{p_0}$ can  be  joined from
the     surface     $S_0^+$    by     a     timelike    curve,     see
Remark~\ref{r.past-future-surface}). Hence, the point $z$ is on the 
surface $S_0^+$. This is absurd, since $z$ must be on
the  boundary of~$U_0^+$,  and  since $U_0^+$  is  a neighbourhood  of
$S_0^+$. This completes the proof of the claim. The lemma follows
immediately.  
\end{proof}

\paragraph{Definition of the set $C_1$.} 
From   now  on,   we  fix   a  number   $\epsilon>0$  such   that  the
$\epsilon$-pasts  and $\epsilon$-futures of  the surfaces  $S_0^-$ and
$S_0^+$ are contained in $D(S_0)$. We consider the set
$$C_1:=C_0\cup I_\epsilon^-(S_0^-)\cup I_\epsilon^+(S_0^+)$$
Obviously,  $C_1$ is  a $\Gamma$-neighbourhood  of $C_0$  contained in
$D(S_0)$. Actually, $C_1$ should be thought as a ``Lorentzian
$\epsilon$-neighbourhood" of $C_0$.

\bigskip
Our aim  is to  prove that the  boundary of the  set $\Gamma\backslash
C_1$ is made
of  two  topological  Cauchy  surfaces which  are  convex/concave  and
uniformly curved. For that purpose,  we first need to prove that $C_1$
is a convex set. Let us introduce some notations.  We denote by $\cP(S_0^-)$ 
(resp. by $\cP(S_0^+)$)  the set of the support  planes of the surface
$S_0^-$ (resp.  the surface $S_0^+$). 

\begin{lemma}
\label{l.convex-comme-intersection}
The set $C_1$ is made of the points $p\in\cA_{p_0}$ such that:\\
- for every plane $P$ in $\cP(S_0^+)$, the point $p$ is in the past or
in the $\epsilon$-future of $P$, \\
- for every plane $P$ in $\cP(S_0^-)$,  the point $p$ is in the future
or in the $\epsilon$-past of $P$. \\
In other words: 
\begin{equation}
\label{eq.convex-comme-intersection}
C_1=\left(\bigcap_{P\in\cP(S_0^-)}I^-_\epsilon(P)\cup
I^+(P)\right)\cap\left(\bigcap_{P\in\cP(S_0^+)}I^-(P)\cup
I^+_\epsilon(P)\right)
\end{equation}
\end{lemma}

\begin{proof}
We     denote     by     $C_1'$     the     right-hand     term     of
equality~(\ref{eq.convex-comme-intersection}). Let  $p$ be a  point of
$\cA_{p_0}$ which is not in $C_1'$. Assume for instance that there
exists a  plane $P\in\cP(S_0^+)$,  such that $p$  is in the  future of
$P$,  and the  distance from  $p$ to  $P$ is  bigger  than $\epsilon$.
Since the surface $S_0^+$
is in the past of $P$, this 
implies that  $p$ is in  the future of  $S_0^+$ and that  the distance
from $p$ to $\Sigma_0^+$ is  bigger than~$\epsilon$. Hence, $p$ is not
in $C_1$.

Conversely, let $p$  be a point of $\cA_{p_0}$ which  is not in $C_1$.
Assume for instance that $p$ is in the future of the surface $S_0^+$
and  the distance  from  $p$  to $S_0^+$  is  bigger than  $\epsilon$.
Then there  exists a  timelike curve $\gamma$  joining $p$ to  a point
$q\in\Sigma_0^+$,  such that  the length  of $\gamma$  is  bigger than
$\epsilon$. Let $P$ be a support  of $C_0$ such that $q\in P\cap C_0$.
By definition, $P$ is an element  of $\cP(S_0^+)$, the point $p$ is in
the future of $P$, and the distance from $p$ to $P$ is bigger than the
length of $\gamma$. Hence, $p$ is not in $C_1'$. 
\end{proof}

Using the diffeomorphism $\Phi_{p_0}$ (see
subsection~\ref{ss.affine-domains}),   we identify   the  domain
$\cA_{p_0}$ with the region of $\RR^3$ where $x^2+y^2-z^2<1$.  
Let $P_0$ be the totally geodesic subspace of
$\cA_{p_0}$ defined as the intersection of $\cA_{p_0}$ with the affine
plane $(z=0)$ in $\RR^3$. Obviously, $P_0$ is spacelike.

\begin{lemma}
\label{l.epsilon-voisinage}
The set $I^-(P_0)\cup I^+_\epsilon(P_0)$  is the region of $\cA_{p_0}$
defined by the inequation
$$z\leq \tan\epsilon.\sqrt{1-x^2-y^2}$$
\end{lemma}

\begin{proof}
All  the calculations  have to  be  made in  the linear  model of  the
anti-de Sitter space, using the coordinates $x_1,x_2,x_3,x_4$ (because
in this model the Lorentzian metric is simply 
the restriction of a global quadratic form).   
The equation  of $P_0$  in  this system  of coordinates  is
$(x_1=0)$.      The    equation    $(z=\tan\epsilon.\sqrt{1-x^2-y^2})$
corresponds to  the equation  $(x_1=\sin\epsilon)$.  On the  one hand,
since $P_0$ is a smooth spacelike surface, the distance from a point $q\in
D(P_0)$ to  the plane $P_0$  is realized as  the length of  a geodesic
segment  joining  $q$ to  $P_0$  and  orthogonal  to $P_0$  (see,  for
instance,~\cite{Haw}).        On        the       other        hand,
Proposition~\ref{p.geodesics} implies every point $q$ on the surface
$(x_1=\sin\epsilon)$ belongs to a  unique geodesic which is orthogonal
to  $P_0$.  So, we  are left  to prove  that, for  every point  $p$ on
$P_0$,  the length of the  unique segment  of geodesic  orthogonal to
$P_0$ and joining $p$ to  the surface $(x_1=\sin\epsilon)$ is equal to
$\epsilon$.  This follows  from Proposition~\ref{p.geodesics} and from
an elementary calculation. 
\end{proof}

\begin{rema}
\label{r.epsilon-voisinage}
Lemma~\ref{l.epsilon-voisinage}   shows    that   $I^-(P_0)\cup
I^+_\epsilon(P_0)$ is a relatively convex subset of $\AdS_3$ (that is,
the intersection of a convex subset of $\SS^3$ with $\AdS_3$). 
Moreover, it shows that there exists $R$ such that the boundary of the
set $I^-(P_0)\cup I^+_\epsilon(P_0)$ is more curved than  a sphere of
radius $R$ at every point: if we consider the Euclidean metric on
$\RR^3$ for  which $(x,y,z)$ is an orthonormal  system of coordinates,
then we can take $R=(\tan\epsilon)^{-1}$.  
Although   this   does   not   clearly   appear  in   the   proof   of
Lemma~\ref{l.epsilon-voisinage}, this  phenomenon is related  with the
negativity of the curvature of $\AdS_3$.  
\end{rema}

\begin{coro}
\label{c.convex}
The set $C_1$ is convex.
\end{coro}

\begin{proof}
Consider  a  totally  geodesic  subspace $P\in  \cP^+(S_0^+)$.   There
exists $\sigma_P\in O_0(2,2)$, such that $\gamma_P(P_0)=P$. Then, 
$\sigma_P$ maps  the set  $I^-(P_0)\cup I^+_\epsilon(P_0)$ to  the set
$I^-(P)\cup I^+_\epsilon(P)$. By Remark~\ref{r.epsilon-voisinage}, the
set $I^-(P_0)\cup I^+_\epsilon(P_0)$  is relatively convex. Hence, the
set $I^-(P)\cup I^+_\epsilon(P)$ is also relatively convex.  The same
arguments show that, for every 
$P\in\cP^-(C_0)$, the set $I^-_\epsilon(P)\cup I^+(P)$ is relatively
convex.  Together with  Lemma  \ref{l.convex-comme-intersection}, this
shows that the set $C_1$ is a relatively convex subset of $\AdS_3$. 
Moreover, $C_1$ is  contained in $D(S_0)$, which is  a convex subset of
$\AdS_3$   (see    item~(iii)   of   Remark~\ref{r.black-domain}   and
Proposition~\ref{p.B=D}). Therefore, 
$C_1$ is a convex subset of $\AdS_3$.  
\end{proof}

\paragraph{Definition of the surfaces $S_1^-$, $S_1^+$, $\Sigma_1^-$
  and  $\Sigma_1^+$.}  The set  $C_1$  is  a $\Gamma$-invariant  closed
convex subset of $\AdS_3$, containing $C_0$, and contained in $D(S_0)$.  
By  Proposition~\ref{p.boundary-convex-set}, the
boundary of $C_1$  in $\AdS_3$ is the union  of two $\Gamma$-invariant
spacelike topological surfaces $S_1^-$  and $S_1^+$, such that 
$(S_1^-, S_1^+)$ is a convex trap. Also by Proposition~\ref{p.boundary-convex-set},
$(\Sigma_1^-:=\Gamma\backslash     S_1^-, \Sigma_1^+:=\Gamma\backslash 
S_1^+)$ is a convex trap.

\begin{rema}
\label{r.def-alternative-S_1}
The surface $S_1^-$ (resp. $S_1^+$) is  the set made of all the points
of  $\cA_{p_0}$  which  are  in   the  past  of  the  surface  $S_0^-$
(resp.   $S_0^+$),   at  distance   exactly   $\epsilon$  of   $S_0^-$
(resp. $S_0^+$):  this follows from  the definition of the  set $C_1$,
and  from the  continuity of  the  distance from  a point  $p$ to  the
surface    $S_0^-$    (resp.    $S_0^+$)    when   $p$    ranges    in
$D(S_0)=D(S_0^-)=D(S_0^+)$.    Thus,    the    surface    $\Sigma_1^-$
(resp. $\Sigma_1^+$) is the set made of all the points of $M$ which 
are in the past of  the surface $\Sigma_0^-$ (resp.  $\Sigma_0^+$), at
distance exactly $\epsilon$ of $\Sigma_0^-$ (resp. $\Sigma_0^+$).
\end{rema}

\begin{prop}
\label{p.uniforme-convexite}
The  surfaces  $\Sigma_1^-$  and  $\Sigma_1^+$ are  uniformly  curved.
\end{prop}

\begin{proof}
Fix a  Euclidean metric on $\RR^3$, and  let $\Delta_1^+\subset S_1^+$
be a compact fundamental domain of the action of $\Gamma$ on
$S_1^+$.  
Let $\Delta_0^+$ be  the intersection  of the  past of
$\Delta_1^+$  with  the surface  $S_0^+$.  Note  that $\Delta_0^+$  is
compact  (since $\Delta_1^+$  is compact, and since  $\Delta_1^+$ and
$S_0^+$   are   contained   in   a  globally   hyperbolic   subset   of
$\AdS_3$). Let $\cP(\Delta_0^+)$ be the  set of all the support planes
of $S_0^+$ that meet $S_0^+$ at some point of $\Delta_0^+$. 

\smallskip

\noindent  \emph{Claim  1. There  exists  $R$ such  that,  for  every
  $P\in\cP(\Delta_0^+)$,   the  boundary   of   the  set   $I^-(P)\cup
  I_\epsilon^+(P)$ is more curved than a sphere of radius $R$.}

\smallskip

\noindent On  the one hand,  $\cP(\Delta_0^+)$ is a compact  subset of
the set of all spacelike  totally geodesic subspaces of $\AdS_3$. As a
consequence, there  exists a compact  subset $\cK$ of  $O_0(2,2)$ such
that $\cP(\Delta_0^+)\subset\cK.P_0$. On  the other hand, there exists
$R_0$   such   that   the    boundary   of   the   set   $I^-(P_0)\cup
I_\epsilon^+(P_0)$ is more  curved than a sphere of  radius $R_0$ (see
Remark~\ref{r.epsilon-voisinage}). The claim follows.

\smallskip

\noindent \textit{Claim 2. Every $q\in\Delta_1^+$ is on the boundary of
the set $I^-(P)\cup I_\epsilon^+(P)$ for some $P$ in $\cP(\Delta_0^+)$.}

\smallskip

\noindent Let $q\in\Delta_1^+\subset S_1^+$. By definition of $S_1^+$,
the point $q$ is in the future of the surface $S_0^+$ and the distance
from $q$  to $S_0^+$ is $\epsilon$. Moreover,  since $q$ and
$S_0^+$ are contained in a  globally hyperbolic region of $\AdS_3$, the
distance between  $q$ and $S_0^+$  is realized: there exists  a causal
curve  $\gamma$ of  length $\epsilon$  joining  $q$ to  a point $p\in
S_0^+$. By construction, the point  $p$ is in $\Delta_0^+$. Let $P$ be
any  support  plane   of  $S_0^+$  at  $p$.  Of   course,  $P$  is  in
$\cP(\Delta_0^+)$.  On the one hand, since $\gamma$ is a causal arc of
length $\epsilon$ joining $q$ to a point of $P$, the distance from $p$
to   $P$    is   at   least    $\epsilon$.   On   the    other   hand,
Lemma~\ref{l.convex-comme-intersection} implies that the distance from
$p$ to $P$ must be at most $\epsilon$.  The claim follows. 

Let  $q$  be  a  point  of  $\Delta_1^+$.  By  claim~2,  there  exists
$P\in\cP(\Delta_0^+)$ such that $q$ is on the boundary of the set 
$I^-(P)\cup                    I_\epsilon^+(P)$.                    By
Lemma~\ref{l.convex-comme-intersection},   the   surface  $S_1^+$   is
contained in $I^-(P)\cup  I^+_\epsilon(P)$. Putting these together with
claim~1,  we obtain that  the surface  $S_1^+$ is  more curved  than a
sphere  of radius  $R$  at  $q$. Hence,  the  surface $\Sigma_1^+$  is
uniformly curved.
\end{proof}

This completes the proof of Proposition~\ref{p.uniformisation}. 

\begin{rema}
\label{r.any-surfaces}
All the results of this subsection are still  valid if one replaces
$(\Sigma_0^-, \Sigma_0^+)$ by any other convex trap.
\end{rema}

\begin{rema}
It is well-known that  the boundary of the $\epsilon$-neighbourhood of
any  geodesically convex  subset  of  $\RR^n$ or  $\HH^n$  is a  $C^1$
hypersurface. Unfortunately, this phenomenon  does not  generalize to
Lorentzian setting. In particular, the surfaces $S_1^-$, $S_1^+$,
$\Sigma_1^-$ and $\Sigma_1^+$ are not $C^1$ in general. 
\end{rema}

\subsection{Smoothing the Cauchy surfaces $\Sigma_1^-$ and $\Sigma_1^+$}
\label{s.smoothing}

In order to apply Gerhard's  theorem, we need a \emph{smooth} uniformly 
curved convex trap. The
purpose of this subsection is to prove the following proposition: 

\begin{prop}
\label{p.smoothing}
Arbitrarily close to $\Sigma_1^-$  and $\Sigma_1^+$, there exist some
$C^\infty$  Cauchy  surfaces   $\Sigma^-$  and~$\Sigma^+$,  such  that
$\Sigma^-$ is  convex and uniformly  curved and $\Sigma^+$  is concave
and uniformly curved.
\end{prop}
 
Unfortunately, we    could   not    find   any   simple    proof   of
Proposition~\ref{p.smoothing}  (see  Remark~\ref{r.convolution}).  Our
proof is divided in three steps. In~\ref{ss.step2}, we approximate the
surfaces $\Sigma_1^-$ 
and   $\Sigma_1^+$   by   some   \emph{polyhedral}   Cauchy   surfaces
$\Sigma_2^-$ and $\Sigma_2^+$ (respectively convex and concave). Then,
in~\ref{ss.step3}, we  describe  a  method for  smoothing
convex and  concave polyhedral Cauchy surfaces. Using  this method, we
obtain  two  disjoint  $C^\infty$  Cauchy  surfaces  $\Sigma_3^-$  and
$\Sigma_3^+$,    respectively    convex    and   concave.     Finally,
in~\ref{ss.step4}, using the same trick as in subsection
\ref{s.uniformisation},  we  get  a smooth uniformly convex trap.

\begin{rema}
\label{r.convolution}
The first idea  which comes to in mind for smoothing  a convex surface is
to use some  convolution process. Unfortunately, to make  this kind of
idea  work,  one  needs  a  locally  Euclidean structure\footnote{For
  example,   any   convex   function~$f:\RR^n\rightarrow\RR$  can   be
  approximated by  a smooth convex function $\widehat  f$, obtained as
  the convolution of  $f$ with an approximation of  the unity, but
  the  proof of  the  convexity  of $\widehat  f$  uses the  Euclidean  
  structure of $\RR^{n+1}$.}.  
This is  the reason why this kind  of idea does not  fit our situation
(there is no locally Euclidean structure on the manifold $M$). 
\end{rema}

\subsubsection{Polyhedral convex and concave Cauchy surfaces}
\label{ss.step2}

In  this subsubsection, we  will  define  a  notion of  \emph{polyhedral
surface}  in $M$.   Then,  we will construct two  polyhedral Cauchy  surfaces
$\Sigma_2^-$  and $\Sigma_2^+$  in  $M$, such  that $\Sigma_2^-$  is
convex, $\Sigma_2^+$  is concave, and $\Sigma_2^+$ is  in the future
of $\Sigma_2^-$.

A  subset $\Delta$  of $M$  is a \emph{$2$-simplex}, if there
exists  a projective chart  $\Phi:U\subset M\rightarrow\RR^3$,  such that
$\Delta\subset  U$ and such  that $\Phi(\Delta)$  is a  $2$-simplex in
$\RR^3$. A compact surface $\Sigma$ in $M$ is called \emph{polyhedral}
if it can be decomposed as a finite union of $2$-simplices.

\begin{rema}
\label{r.developpement-polyhedral}
Let   $\Sigma$  be   a   compact  spacelike   surface   in  $M$,   let
$\widetilde\Sigma$ be  a lift of  $\Sigma$ in $\widetilde M$,  and let
$S:=\cD(\widetilde\Sigma)$.          Using        the        embedding
$\Phi_{p_0}:\cA_{p_0}\rightarrow\RR^3$, we can see $S$ as a surface in
$\RR^3$. Then, $\Sigma$ is a polyhedral surface if and only if $S$ can
be  decomposed  as  a  finite   union  of  orbits  (for  $\Gamma$)  of
$2$-simplices of $\RR^3$.
\end{rema}

\begin{rema}
Let $\Sigma$ be a compact  convex spacelike polyhedral surface in $M$.
Then, $\Sigma$ can be decomposed as  a  finite union  of  subsets
$\Sigma:=\Delta_1\cup\dots\cup\Delta_n$, where  each $\Delta_i$ is the
intersection  of $\Sigma$  with one  of its  support planes,  and each
$\Delta_i$  has non-empty  interior (as  a subset  of  $\Sigma$).  The
decomposition is  unique (provided that the  $\Delta_i$'s are pairwise
distinct).  The $\Delta_i$'s are called the \emph{sides} of $\Sigma$. Each
side  of $\Sigma$  is  a finite  union  of $2$-simplices,  but is  not
necessarily a  topological disc  (e.g. in the  case where  $\Sigma$ is
totally geodesic).
\end{rema}

\paragraph{Definition of the set $C_2$, of the surfaces $S_2^-$,
$S_2^+$,    $\Sigma_2^-$    and    $\Sigma_2^+$}   We    consider    a
$\Gamma$-invariant  set  $E$  of  points  of  $\partial  C_1=S_1^-\cup
S_1^+$, such that $\Gamma\backslash E$ is finite (in particular, $E$ is
discrete).   We  denote   by  $C_2$  the  convex  hull   of  $E$.   By
construction, $C_2$ is a $\Gamma$-invariant convex subset of $C_1$. In
particular, $C_2$ is a $\Gamma$-invariant convex subset of $D(S_0)$. So,
by Proposition  \ref{p.boundary-convex-set}, 
the pair of boundary components of $C_2$ in $\AdS_3$, and their projections in
$M$, are convex traps.

\bigskip

Given $\delta>0$, we  say that the set $E$  is \emph{$\delta$-dense in
  the surfaces $S_1^-$ and $S_1^+$}, if every Euclidean   ball of radius
  $\delta$ centered at some  point of $S_1^-$ (resp. $S_1^+$) contains
  some points of  $E$.  The remainder of the  subsection is devoted to
  the proof of the following proposition: 

\begin{prop}
\label{p.polyhedral}
There exists $\delta>0$ such that, if the set $E$ is $\delta$-dense in
the   surfaces   \hbox{$S_1^-$  and   $S_1^+$},   then  the   surfaces
$\Sigma_2^-$ and $\Sigma_2^+$ are polyhedral.
\end{prop}

\begin{rema}
The proof of Proposition~\ref{p.polyhedral} is quite technical.
The reader  who is  not interested in  technical details can  skip the
proof.  Nevertheless, it  should be
noticed that  the boundary  of the  convex hull of  a discrete  set of
points is \emph{not} a  polyhedral surface in general.  In particular,
Proposition  \ref{p.polyhedral}   would  be  false   if  the  surfaces
$\Sigma_1^-$ and $\Sigma_1^+$ were not uniformly curved. 
\end{rema}

Given  a set  $F\subset\RR^3$,  we say  that  an affine  plane $P$  of
$\RR^3$ \emph{splits} the set $F$, if $F$ intersects the two connected
components of $\RR^3\setminus P$.  The  starting point of the proof of
Proposition \ref{p.polyhedral} is the following well-known fact (which
follows from very basic arguments of affine geometry):

\begin{fact}
\label{f.polyhedral}
For every finite set of points $F\subset\RR^3$, the
boundary  of   $\Conv(F)$  is  a  compact   polyhedral  surface;  more
precisely,  the  boundary  of  $\Conv(F)$  is the  union  of  all  the
2-simplices $\Conv(p,q,r)$,  such that the points $p,q,r$  are in $F$,
and such that the plane $(p,q,r)$ does not split~$F$.
\end{fact}

Let $\gamma$ be a continuous curve  in a Euclidean plane, and $p$ be a
point  on $\gamma$.   We say  that  the curve  $\gamma$ is  \emph{more
curved than a circle of radius $R$ at $p$} if there exists a Euclidean
disc  $\Delta$ of  radius $R$,  such that  $p$ is  on the  boundary of
$\Delta$, and  such that $\Delta$  contains a neighbourhood of  $p$ in
$\gamma$.  The  proof of  the  following  lemma  of elementary  planar
geometry is left to the reader:

\begin{lemma}
\label{l.polyhedral}
Given two  positive numbers  $\rho$ and $R$,  there exists  a positive
number $\delta=\delta(\rho,R)$ such that:  for every convex set $D$ in
an Euclidean plane, if there  exists a subarc $\alpha$ of the boundary
of $D$, such that the boundary of  $D$ is more curved than a circle of
radius $R$  at each point of  $\alpha$, and such that  the diameter of
$\alpha$ is bigger than $\rho$,  then $D$ contains a Euclidean ball of
radius $\delta$.
\end{lemma}

\begin{figure}[ht]
\centerline{\input{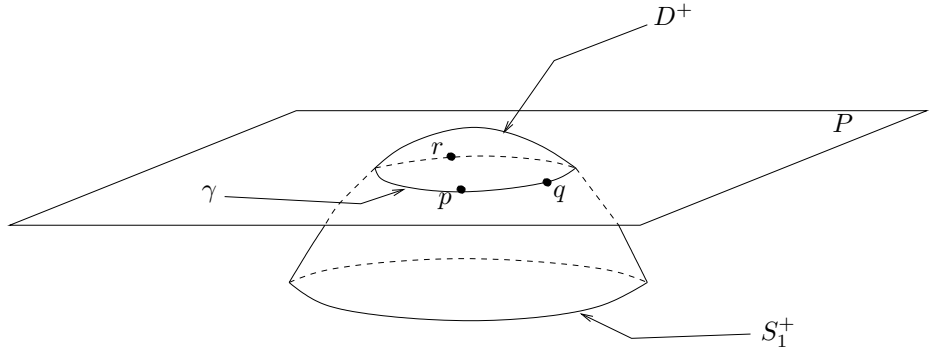}}
\caption{\label{f.preuve-polyhedrale}The situation in the proof of 
Proposition~\ref{p.polyhedral}}
\end{figure}

\begin{proof}[Proof of Proposition \ref{p.polyhedral}]
Consider a compact  fundamental domain $U$ for the  action of $\Gamma$
on $C_1$. Then, consider a  compact neighbourhood $V$ of $U$ in $C_1$,
and  a compact  neighbourhood $W$  of  $V$ in  $C_1$. One  can find  a
positive  number $\rho$  such every  Euclidean   ball  of  radius $\rho$
centered  in  $U$  (resp.   $V$)  is contained  in  $V$  (resp.   $W$).
Moreover, since  $V$ is compact, one  can find a  positive number $R$,
such that  the surface $S_1^-$ (resp.  $S_1^+$) is more  curved than a
sphere of radius $R$ at every point of $S_1^-\cap V$ (resp. $S_1^+\cap
V$).  

From  now on,  we assume  that the  set $E$  is $\delta$-dense  in the
surfaces  $S_1^-$ and  $S_1^+$, where  $\delta=\delta(\rho,R)$  is the
positive  number given  by Lemma  \ref{l.polyhedral}. Up  to replacing
$\delta$  by  $\min(\delta,\rho)$,  we  can assume  that  $\delta$  is
smaller than $\rho$. Under these  assumptions, we shall prove that the
surfaces $S_2^-$ and $S_2^+$ are polyhedral.  

\medskip

\noindent \emph{Claim 1. If $p,q,r$ are three points of $E$, such that
the 2-simplex $\Conv(p,q,r)$ intersects  $U$, and such that the affine
plane $P:=(p,q,r)$ does  not split the set $E$,  then the three points
$p,q,r$ are in $W$.}

\medskip

\noindent To prove this claim, we argue by contradiction: we suppose
that there exists three points $p,q,r$ in $E$, such that the 2-simplex
$\Conv(p,q,r)$ intersects $U$, such that the affine plane $P:=(p,q,r)$
does not  split the  set $E$, and  such that  one of the  three points
$p,q,r$  is  not  in  $W$.   We shall  show  that  these  assumptions
contradict the $\delta$-density of the set $E$.

Since  $P$ does  not  split the  set  $E$, one  of  the two  connected
components of  $\cA_{p_0}\setminus P$ is disjoint from  $E$. We denote
by $H_P$ this connected component. First of all, we observe that $H_P$
does  not intersect  the curve  $\partial S_0$,  since $H_P$  does not
contain any point of $E$,  since $E$ is a non-empty $\Gamma$-invariant
subset of  $D(S_0)$, and since the  curve $\partial S_0$  is the limit
set of the action of $\Gamma$ on $D(S_0)$. Therefore, the 
intersection of $H_P$ with the boundary of $C_1$ 
is contained in one of the two connected components $S_1^-$ and $S_1^+$
of $\partial  C_1\setminus\partial S_0$.  Without  loss of generality,
we assume that $H_P\cap\partial C_1$  is contained in $S_1^+$, and we
consider the set $D^+:=H_P\cap S_1^+$ (see figure~\ref{f.preuve-polyhedrale}). 

We  shall prove  that there  exists an  Euclidean ball  $B$  of radius
$\delta$  centered   at  some  point   of  $D^+$,  such   that  $B\cap
S_1^+\subset  D^+$. Since  $D^+$ must  be disjoint  from  $E$ (because
$D^+\subset  H_P$),  this  will   contradict  the  fact  that  $E$  is
$\delta$-dense in  $S_1^+$.  For that  purpose, we consider  the curve
$\gamma:=P\cap S_1^+$. Observe that this curve $\gamma$ is the boundary of
the topological disc  $D^+$. Moreover, the curve $\gamma$  is also the
boundary of  the convex  subset $D:=P\cap C_1$  of the plane  $P$. The
curve $\gamma$  passes through  the points $p$,  $q$ and $r$,  and the
2-simplex $\Conv(p,q,r)$ is contained in  the convex set $D$.  We shall
distinguish two cases (and get a contradiction in each case):

\emph{First  case:   the  curve   $\gamma$  does  not   intersect  the
neighbourhood $V$.} We consider a point $m$ in $D\cap U$ (such a point
does  exist   since  $\Conv(p,q,r)\cap  U\neq\emptyset$  and
$\Conv(p,q,r)\subset D$),  and we denote  by $m'$ the unique  point of
intersection of $D^+$ with the line passing through $m$ and orthogonal
to the plane $P$. The point $m$ is in $U$, and the curve $\gamma$ does
not intersect $V$; so, by definition of $\rho$, the Euclidean distance
between $m$ and $\gamma$ must  be bigger than $\rho$, and thus, bigger
than $\delta$. Moreover, the Euclidean distance between the point $m'$
and the  curve $\gamma$  is bigger than  the distance between  $m$ and
$\gamma$. So,  we have  proved that the  Euclidean ball $B$  of radius
$\delta$   centered   at   $m'$   does   not   intersect   the   curve
$\gamma$. Hence,  the connected component of  $B\cap S_1^+$ containing
the point $m'$ is contained in $D^+$. Since $D^+$ is disjoint from $E$,
this contradicts the $\delta$-density of $E$ in $S_1^+$.

\emph{Second case: the curve $\gamma$ does intersect the neighbourhood
$V$.} Then, by definition of $\rho$, we can find a subarc $\alpha$ of
the curve $\gamma$, such that  the diameter of $\alpha$ is bigger than
$\rho$, and such  that $\alpha$ is contained in  $W$.  Since $S_1^+$ is
more curved  than a sphere  of radius $R$  at every point of  $V$, the
curve $\gamma$  is more  curved than  a circle of  radius $R$  at each
point of $\alpha$. Thus, by lemma~\ref{l.polyhedral}, there exists a
point $m\in D$ such that  the Euclidean distance between the point $m$
and the curve $\gamma$ is  bigger than $\delta$.  The same argument as
above shows that  this contradicts the $\delta$-density of  $E$ in the
surface $S_1^+$.

In both case, we have obtained a contradiction. This completes 
the proof of claim~1. 

\medskip

\noindent \emph{Claim  2. If $W'$  is a compact subset  of $\cA_{p_0}$
  such that $W\subset  W'$, then the sets $\Conv(E\cap  W')\cap U$ and
  $\Conv(E\cap W)\cap U$ coincide.}  

\medskip

\noindent This claim is a consequence of Claim 1 and fact
\ref{f.polyhedral}. Since $W'$  is a compact subset of  $\AdS_3$, the set
$E\cap W'$ is finite. Hence, the boundary of the set $\Conv(E\cap W')$
is the union of the  2-simplices $[p,q,r]$, such that the three points
$p,q,r$ are  in $E\cap W'$, and  such that the  affine plane $(p,q,r)$
does not split $E\cap W'$. By  claim 1, such a 2-simplex can intersect
$U$ only if the three points $p$,  $q$ and $r$ are in $W$.  Using once
again fact \ref{f.polyhedral}, this implies that the intersection of boundary of
$\Conv(E\cap W')$ with $U$ is contained in the intersection of the boundary of
$\Conv(E\cap  W)$ with $U$.  But,  if the  boundary of  a
convex set  is contained  in the boundary  of another convex  set, then
these two convex sets must coincide. The claim follows.

\medskip

\noindent\emph{End of the proof.} Let us consider an increasing
sequence  $(W_n)_{n\in\NN}$ of  compacts subsets  of $\AdS_3$,  such that
$\bigcup_{n\in\NN}  W_n=\AdS_3$.   On  the  one  hand,  we  clearly  have
$\Conv(E)=\mbox{Closure}(\bigcup_{n\in\NN}\Conv(E\cap   W_n))$.   On  the
other hand, according  to Claim 2, there exists  an integer $n_0$ such
that $\Conv(E\cap  W_n)\cap U=\Conv(E\cap  W)\cap U$ for  every $n\geq
n_0$.  As  a consequence,  we have $\Conv(E)\cap  U=\Conv(E\cap W)\cap
U$. Now, since $E\cap W$ is a finite set, the boundary of $\Conv(E\cap
W)$ is  a compact  polyhedral surface.  Thus,  we have proved  that the
boundary of the set $C_2=\Conv(E)$ coincides, in $U$, with a
polyhedral surface.   Since $U$ contains a fundamental  domain for the
action of $\Gamma$ on $C_2$, this implies each of the surfaces $S_2^-$
and  $S_2^+$  can  be  decomposed  as  a finite  union  of  orbits  of
$2$-simplices. Hence,  the surfaces $\Sigma_2^-$  and $\Sigma_2^+$ are
polyhedral (see Remark~\ref{r.developpement-polyhedral}). 
\end{proof}

\begin{adde}
There exists $\delta>0$ such that, if the set $E$ is $\delta$-dense in
the surfaces $S_1^+,S_1^-$, then  each side of the polyhedral surfaces
$\Sigma_2^-,\Sigma_2^+$ is contained in the domain of an 
projective chart of $M$.
\end{adde}

\begin{proof}
From the proof of  Proposition~\ref{p.polyhedral}, one can extract the
following 
statement: for  every $\rho>0$, there exists $\delta>0$  such that, if
the set $E$ is $\delta$-dense  in the surface $S_1^-$, then, for every
support  plane $P$ of  the surface  $S_2^-$, the  diameter of  the set
$P\cap  S_2^-$ is  less than  $\rho$.  Of  course, there  is  a similar
statement for the surface $S_2^+$. The addendum follows immediately.
\end{proof}

\subsubsection{Smooth convex and concave Cauchy surfaces}
\label{ss.step3}

In this subsubsection, we describe a process for smoothing the polyhedral
Cauchy surfaces $\Sigma_2^-$ and $\Sigma_2^+$.  More precisely, we 
prove the following:

\begin{prop}
\label{p.lissage-surface-polyhedrale}
Let $\Sigma$ be a convex polyhedral Cauchy surface in $M$. Assume that
each side of $\Sigma$ is contained in an affine domain of $M$. Then,
arbitrarily close to $\Sigma$, there exists a $C^\infty$ convex Cauchy
surface.  
\end{prop}

Of  course,  the  analogous  statement  dealing  with  concave  Cauchy
surfaces       is       also       true.      The       proof       of
Proposition~\ref{p.lissage-surface-polyhedrale}    relies    on    the
following technical lemma: 
 
\begin{lemma}
\label{l.lissage-fonction-convexe}
Let  $U$  be some  subset  of  $\RR^2$  and $f:U\rightarrow\RR$  be  a
continuous convex  function. Then, for every $\eta>0$,  there exists a
continuous convex function $\widehat f:U\rightarrow\RR$ satisfying the
following properties:

\smallskip

\noindent $\bullet$ $\widehat f\geq f$, the distance between $f$ and 
$\widehat f$ is less than $2\eta$, and $\widehat f$ coincides with $f$ on
the set $f^{-1}\left([2\eta,+\infty[\right)$;

\smallskip

\noindent   $\bullet$   $\widehat   f$   is  constant   on   the   set
$f^{-1}([0,\eta])$; in  particular, $\widehat f$ is  $C^\infty$ on the
set $f^{-1}\left([0,\eta[\right)$; 

\smallskip

\noindent  $\bullet$  if $f$  is  $C^\infty$  on  some subset  $U$  of
$\Dom(f)$, then $\widehat f$ is also $C^\infty$ on $U$. 
\end{lemma}

\begin{proof}
We   consider  a   $C^\infty$   function  $\phi:[0,+\infty[\rightarrow
[0,+\infty[$   such  that:  $\phi$   is  non-decreasing   and  convex,
$\phi(t)=\frac{3}{2}\eta$ for  every $t\in [0,\eta]$,  and $\phi(t)=t$
for  every  $t\in [2\eta,+\infty[$.  Then,  consider 
$\widehat   f:   U\rightarrow   [0,+\infty[$  defined   by   $\widehat
f:=\phi\circ f$.  This function satisfies all  the required properties.
\end{proof}

We endow $M$ with a Riemannian  metric; this allows us to speak of the
(Riemannian)  $\epsilon$-neighbourhood of  any subset  of $M$  for any
$\epsilon>0$. We say that a surface $\Sigma_1$ is
\emph{$\epsilon$-close} to another surface $\Sigma_2$ if there
exists  a  homeomorphism  $\Psi:\Sigma_1\rightarrow\Sigma_2$ which  is
$\epsilon$-close to the identity. 
The following  Remark will allow us  to see a polyhedral  surface as a
collection of graphs of functions:

\begin{rema}
\label{r.graph}
Let  $\Sigma$ be  a convex  compact  surface in  $M$, let  $\Pi$ be  a
support plane of $\Sigma$  and let $\Delta:=\Sigma\cap\Pi$.  We assume
that $\Delta$  is contained in  an affine domain  of $M$. Then,  we can
find  a  neighbourhood  $V$ of  $F$  in  $M$,  and some  local  affine
coordinates $(x,y,z)$ on $V$, such that:

\smallskip

\noindent -- $\Pi\cap V$ is the plane of equation $(z=0)$, and
$\Sigma\cap V$ is the graph $(z=f(x,y))$ of a non-negative convex function
$f:U\rightarrow [0,+\infty[$ (where $U$ is some convex subset of $\RR^2$).

\smallskip

\noindent -- if $\Sigma'$ is a convex Cauchy surface close
enough to $\Sigma$, then $\Sigma'\cap V$ is the graph $z=f'(x,y)$ of
a convex function $f':U\rightarrow\RR$. 
The function $f'$ depends continuously of the surface
$\Sigma'$. Moreover, if  $\Sigma'$ is in the future  of $\Sigma$, then
$f'\geq f$ (and thus, $f'\geq 0$). 
\end{rema} 

We  denote by  $\Delta_1,\dots,\Delta_n$ the  sides of  the polyhedral
surface                $\Sigma$.                To               prove
Proposition~\ref{p.lissage-surface-polyhedrale},  we will  construct a
sequence  of convex  Cauchy surfaces  $\Sigma_0,\dots,\Sigma_n$, where
$\Sigma_0=\Sigma$, and  where $\Sigma_{k+1}$ is  obtained by smoothing
$\Sigma_k$ in the neighbourhood of $\Delta_{k+1}$.  More precisely, we
will prove the following: 

\begin{prop}
\label{p.lissage-par-recurrence}
For every $k\in\{0,\dots,n\}$ and every $\epsilon>0$ small enough,
there exists a convex Cauchy surface $\Sigma_{k,\epsilon}$ in $M$ such
that:

\smallskip

\noindent -- the surface $\Sigma_{k,\epsilon}$ is in the future of the
surface $\Sigma$.

\smallskip

\noindent -- the surface $\Sigma_{k,\epsilon}$ is $\epsilon$-close to the
  surface $\Sigma$,

\smallskip

\noindent -- the surface $\Sigma_{k,\epsilon}$ is smooth outside
  the $\epsilon$-neighbourhoods of the sides $\Delta_{k+1},\dots,\Delta_n$.
\end{prop}

Notice    that    Proposition~\ref{p.lissage-par-recurrence}   implies
Proposition~\ref{p.lissage-surface-polyhedrale}    (for   $k=n$,   the
surface  $\Sigma_{k,\epsilon}$  is  a  smooth convex  Cauchy  surface,
$\epsilon$-close to the  initial surface $\Sigma$). So, we  are left to
prove Proposition~\ref{p.lissage-par-recurrence}.

\begin{proof}[Proof of Proposition~\ref{p.lissage-par-recurrence}]
First   of  all,  we   set  $\Sigma_{0,\epsilon}:=\Sigma$   for  every
$\epsilon>0$.  Now, let $k\in\{0,\dots,n-1\}$, and let us suppose that
we  have  constructed  the  surface  $\Sigma_{k,\epsilon}$  for  every
$\epsilon>0$   small   enough.    We   will  construct   the   surface
$\Sigma_{k+1,\epsilon}$ for every $\epsilon>0$ small enough. 

Since $\Delta_{k+1}$  is a  side of $\Sigma$,  there exists  a support
plane $\Pi_{k+1}$ of $\Sigma$ such that $\Pi_{k+1}\cap\Sigma=\Delta_{k+1}$.
Using  Remark~\ref{r.graph}, we  find a  compact neighbourhood  $V$ of
$\Delta_{k+1}$ in $M$, and  some local affine coordinates $(x,y,z)$ on
$V$, such that in these coordinates, $\Pi_{k+1}\cap V$ is the plane of
equation  $(z=0)$,  and  the  surface  $\Sigma\cap  V$  is  the  graph
$(z=f(x,y))$ of a non-negative convex function
$f:\Dom(f)\subset\RR^2\rightarrow\RR$.  Moreover, the function  $f$ is
positive in  restriction to $\partial\Dom(f)$, and  thus, the quantity
$\delta:=\inf\{f(x,y)\mid   (x,y)\in\partial\Dom(f)\}$   is   positive
($\partial\Dom(f)$ is compact).

Now,  we  fix   some  $\epsilon>0$  such  that  $\epsilon/3<\delta/2$.
By the second item of Remark~\ref{r.graph}, we can find
$\epsilon'>0$,  such that  $\epsilon'<\epsilon/3$, and  such  that the
surface $\Sigma_{k,\epsilon'}\cap V$ is the graph of a convex function
$g:\Dom(g)=\Dom(f)\rightarrow\RR$.            Moreover,          since
$\Sigma_{k,\epsilon'}$ is in the  future of $\Sigma$, the function $g$
is bigger  than $f$; in particular,  $g$ is non-negative,  and we have
$g(x,y)>\delta$ for every $(x,y)\in\partial\Dom(g)$.

Applying  Lemma~\ref{l.lissage-fonction-convexe} to  the  function $g$
with   $\eta:=\epsilon/3$,   we  obtain   a   convex  function   $\hat
g:\Dom(g)\rightarrow [0,+\infty[$ satisfying the following properties:

\smallskip

\noindent  (a) $\widehat  g\geq g$  and the  distance between  $g$ and
$\widehat g$ is less than $2\epsilon/3$,

\smallskip

\noindent (b) $\widehat g$ is $C^\infty$ on $g^{-1}([0,\epsilon/3])$,

\smallskip

\noindent   (c)  if   $g$   is   smooth  on   some   open  subset   of
$\Dom(g)=\Dom(\widehat g)$, then $\widehat g$ is also smooth on $U$,

\smallskip

\noindent (d) $\widehat g$ coincides with $g$ on
$g^{-1}([2\epsilon/3,+\infty[)$; in  particular, $\widehat g$ coincide
with $g$ on $\partial\Dom(\widehat g)=\partial\Dom(g)$.

\smallskip

We construct the  surface $\Sigma_{k+1,\epsilon}$ as follows: starting
from    the    surface     $\Sigma_{k,\epsilon'}$,    we    cut    off
$\Sigma_{k,\epsilon}\cap V$ (\emph{i.e.} we cut off the graph of $g$),
and we  paste the graph of  $\widehat g$.  This is  possible since the
graphs  of  the functions  $g$  and  $\widehat  g$ coincide  near  the
boundary  of $V$ (property  (d)).  There  is a natural diffeomorphism
$\Psi$    between     the    surfaces    $\Sigma_{k,\epsilon'}$    and
$\Sigma_{k+1,\epsilon}$ defined as  follows: $\Psi$ coincides with the
identity outside $V$, and maps the point of coordinates $(x,y,g(x,y))$
to the point of  coordinates $(x,y,\widehat g(x,y))$. By property (a),
$\Psi$ is  $(2\epsilon/3)$-close to  the identity; hence,  the surface
$\Sigma_{k+1,\epsilon}$   is  $(2\epsilon/3)$-close  to   the  surface
$\Sigma_{k,\epsilon'}$.       Since      $\Sigma_{k,\epsilon'}$     is
$\epsilon'$-close  to $\Sigma$,  and since  $\epsilon'<\epsilon/3$, we
get that $\Sigma_{k+1,\epsilon}$ is $\epsilon$-close to $\Sigma$.

The inequality $\widehat g\geq g$ implies that $\Sigma_{k+1,\epsilon}$
is in  the future of $\Sigma_{k,\epsilon'}$, and  \emph{a fortiori} in
the future  of $\Sigma$.  The  convexity of the function  $\widehat g$
implies that $\Sigma_{k+1,\epsilon}$ admits a support plane at each of
its   points.    By  Proposition~\ref{p.spacelike-support-plane}   and
Remark~\ref{r.spacelike-support-plane}, this implies that 
$\Sigma_{k+1,\epsilon}$    is    a    spacelike    surface.     Hence,
$\Sigma_{k+1,\epsilon}$ is  a Cauchy surface  (every compact spacelike
surface   embedded  in  $M$   is  a   Cauchy  surface).    Now,  since
$\Sigma_{k+1,\epsilon}$  is a  spacelike surface  admitting  a support
plane at  each point,  it is  either convex or  concave; and  since it
coincides  with  $\Sigma_{k,\epsilon}$   outside  $V$,  it  cannot  be
concave. So, $\Sigma_{k+1,\epsilon}$ is a convex Cauchy surface.

It remains to study the smoothness of $\Sigma_{k+1,\epsilon}$. Let $q$
be a point on the surface $\Sigma_{k+1,\epsilon}$, which is not in the
union    of    the     $\epsilon$-neighbourhoods    of    the    sides
$\Delta_{k+2},\dots,\Delta_n$,                 and                 let
$p:=\Psi^{-1}(q)\in\Sigma_{k,\epsilon'}$.  Since the  distance between
the  points $p$  and $q$  is less  than $2\epsilon/3$,  the  point $p$
cannot be in the union of the $\epsilon/3$-neighbourhoods of the sides
$\Delta_{k+2},\dots,\Delta_n$. There are two cases:

\smallskip

\noindent -- \emph{if the point $p$ is in the $\epsilon/3$-neighbourhood
of the  side $\Delta_{k+1}$}, then the  distance between $p$  and the plane
$\Pi_{k+1}$ is  less than  $\epsilon/3$, and thus,  property~(b) implies
that   the   surface   $\Sigma_{k+1,\epsilon}$   is  smooth   in   the
neighbourhood of $\Psi(p)=q$;

\smallskip

\noindent -- \emph{if the point $p$ is not in the
$\epsilon/3$-neighbourhood  of the side  $\Delta_{k+1}$,} then  the surface
$\Sigma_{k,\epsilon'}$ is smooth in the neighbourhood of $p$ (here, we
use   the  inequality  $\epsilon'<\epsilon/3$);   hence,  property~(c)
implies  that the surface  $\widehat\Sigma_{k,\epsilon}$ is  smooth in
the neighbourhood of $\Psi(p)=q$.

\smallskip

\noindent As a consequence, the surface $\Sigma_{k+1,\epsilon}$ is
smooth outside the union of the $\epsilon$-neighbourhoods of
the    sides    $\Delta_{k+2},\dots,\Delta_n$.    Therefore,   the    surface
$\Sigma_{k+1,\epsilon}$ satisfies all the required by properties. 
\end{proof}

Applying    Proposition~\ref{p.lissage-surface-polyhedrale}   to   the
polyhedral Cauchy  surfaces $\Sigma_2^-$ and $\Sigma_2^+$,  we get two
disjoint  $C^\infty$ Cauchy  surfaces  $\Sigma_3^-$ and  $\Sigma_3^+$,
respectively  convex and  concave, such  that $\Sigma_3^+$  is  in the
future of $\Sigma_3^-$.

\subsubsection{Smooth   uniformly  curved   convex  and   concave  Cauchy
  surfaces} 
\label{ss.step4}

The  Cauchy   surfaces  $\Sigma_3^-$  and  $\Sigma_3^+$  are  smooth,
respectively convex  and concave, but not uniformly  curved. Using the
same trick as in subsection~\ref{s.uniformisation}, we will aproximate 
$\Sigma_3^-$ and  $\Sigma_3^+$ by some smooth  uniformly curved Cauchy
surfaces $\Sigma_4^-$ and $\Sigma_4^+$.

\paragraph{Definition   of  the   Cauchy  surfaces   $\Sigma_4^-$  and
  $\Sigma_4^+$.} 
Let $\epsilon$ be a positive  number. Let $\Sigma_4^+$ be the set made
of the points  $p\in M$, such that  $p$ is in the past  of the surface
$\Sigma_3^+$ and  such that the  distance from $p$ to  $\Sigma_3^+$ is
exactly $\epsilon$.  If $\epsilon$ is small  enough, then $\Sigma_4^+$
is a topological  Cauchy surface which is convex  and uniformly curved
(see                   Remark~\ref{r.any-surfaces}                  and
Remark~\ref{r.def-alternative-S_1}).  We construct similarly a
topological  Cauchy surface $\Sigma_4^-$  which is  concave, uniformly
curved, and contained 
in the past  of $\Sigma_3^-$. By construction, $\Sigma_4^+$  is in the
future of $\Sigma_4^-$. 

\begin{prop}
If $\epsilon$ is small enough, the Cauchy surfaces $\Sigma_4^-$ and
$\Sigma_4^+$ are smooth (of class $C^\infty$). 
\end{prop}

\begin{proof}
We denote by $TM$ the tangent bundle of $M$, by $\pi$ the canonical
projection of  $TM$ on $M$,  and by $(\phi^t)_{t\in\RR}$  the geodesic
flow on $TM$.  We consider the subset $T_N\Sigma_3^+$  of $TM$ made of
the  pairs  $(p,\nu)$  such that  $p$  is  a  point of  the  surface
$\Sigma_3^+$ and  $\nu$ is the  future-pointing unit normal  vector of
$\Sigma_3^+$ at $p$.  

Let $p$  be a  point on the  surface $\Sigma_4^+$. By  construction of
$\Sigma_4^+$,  the  distance  from  $p$  to  $\Sigma_3^+$  is  exactly
$\epsilon$.  Since $M$ is  globally hyperbolic, and since $\Sigma_3^+$
is  a smooth  spacelike  surface,  this implies  that  there exists  a
timelike geodesic segment of  length exactly $\epsilon$, orthogonal to
$\Sigma_3^+$,    joining    $\Sigma_3^+$     to    $p$    (see,    for
example,~\cite[page 217]{Haw}).  As a consequence, the surface $\Sigma_4^+$ is
contained in the set $\pi(\phi^\epsilon(T_N\Sigma_3^+))$.  

We are left to  prove that the set $\pi(\phi^\epsilon(T_N\Sigma_3^+))$
is a smooth surface. Since  $\Sigma_3^+$ is a smooth compact spacelike
surface in $M$,  $T_N\Sigma_3^+$ is a smooth compact  surface in $TM$,
nowhere tangent to the fibers of the projection $\pi$, and hence, for
$\epsilon$  small enough,  $\phi^\epsilon(T_N\Sigma_3^+)$ is  a smooth
compact   surface  in  $TM$,   nowhere  tangent   to  the   fibers  of
$\pi$. Therefore, $\pi(\phi^\epsilon(T_N\Sigma_3^+))$  is a smooth surface
in $M$.  
\end{proof}

\subsection{End of the proof of Theorem~\ref{t.CMC-foliation} in the
  case $g\geq 2$}

In the previous paragraph, we have constructed a 
smooth uniformly curved convex trap $(\Sigma_4^{+}, \Sigma_4^{+})$. By
Proposition~\ref{p.link-convexity-curvature-2},       the      surface
$\Sigma_4^-$ have negative curvature and the surface $\Sigma_4^+$ have
positive curvature.  Thus, $(\Sigma_4^-,\Sigma_4^+)$ is a pair of barriers in
$M$. By Theorem~\ref{t.Moncrief} and~\ref{t.Gerhardt}, the existence
of a pair of barriers implies the existence of a CMC time function. 
This completes the proof  of Theorem~\ref{t.CMC-foliation} in the case
where the genus of the Cauchy surfaces is at least $2$.

\section{Proof of Theorem~\ref{t.CMC-foliation} in the case $g=1$}
\label{s.g=1}

The purpose of this  section is to prove Theorem~\ref{t.CMC-foliation}
in the  case where the genus  of the Cauchy surfaces  of the spacetime
under consideration is~$1$.  According to Remark~\ref{revet finis},
after performing some finite covering if necessary, we can reduce this
case to the case where the Cauchy surface is a $2$-torus.

In  subsection~\ref{ss.torus-universe},  we  define a  class  of
spacetimes,  called  Torus  Universes\footnote{These  spacetimes  were
  already      considered       by      several      authors,      see
  Remark~\ref{r.idem-Carlip}}, and  we will prove  that Torus Universes
admit CMC  time functions (actually, we  construct explicitly a
CMC    time   function    on   any    such   spacetime).    Then,   in
subsection~\ref{ss.proof-genus-1},  we  prove that  every maximal
globally  hyperbolic  spacetime, locally  modelled  on $AdS_3$,  whose
Cauchy  surfaces  are $2$-tori,  is  isometric  to  a Torus  Universe.

\subsection{Torus Universes}
\label{ss.torus-universe}

Consider  the  $1$-parameter   subgroup  of  $SL(2,\RR)$  of  diagonal
matrices $(g^{t})_{t \in {\mathbb R}}$ where: 
$$g^{t} = 
\left(\begin{array}{cc}
                 e^{t} & 0 \\
                 0 & e^{-t}\end{array} \right) = 
e^{t\Delta} 
\;\;\;\;  \mbox{where:}\;\;\; \Delta = \left(\begin{array}{cc}
                 1 & 0 \\
                 0     & -1 \end{array}\right)$$
We denote by  $A$ the set of elements  of $SL(2,\RR) \times SL(2,\RR)$
for which both left and  right components belong to the one-parameter
subgroup  $(g^{t})_{t\in\RR}$. Obviously, $A$ is  a free  abelian Lie
subgroup of rank~$2$ of $SL(2,\RR)\times SL(2,\RR)$. This group acts isometrically
on  $AdS_3$  (recall that  the  isometry  group  of $AdS_3$  can  be
identified       with      $SL(2,\RR)\times       SL(2,\RR)$,      see
subsection~\ref{ss.projective-model}).   We  denote  by  $\Omega$  the
union of spacelike $A$-orbits in $AdS_{3}$.  

We will  see below that  $\Omega$ has four connected  components which
are   open   convex   domains    of   $AdS_3$.   For   any   lattice
$\Gamma\subset A$, the action of $\Gamma$ on $\Omega$ is obviously
free  and  properly discontinuous,  and  preserves  each  of the  four
connected components of $\Omega$. 

\begin{defi}
A Torus Universe  is the quotient $\Gamma\backslash U$  of a connected
component $U$ of $\Omega$ by a lattice $\Gamma$ of $A$.
\end{defi}

\begin{theo}
\label{p.torus-universe-GH}
Every Torus Universe is  a globally hyperbolic spacetime, which admits
a CMC time function.  
\end{theo}

To  prove Theorem~\ref{p.torus-universe-GH},   we  will  use  the
$SL(2,{\mathbb R})$-model  of $AdS_{3}$ (see  subsection \ref{ss.SL2}).  We
recall  that  $SL(2,\RR)\times   SL(2,\RR)$  acts  on  $SL(2,\RR)$  by
$(g_L,g_R).g=g_Lgg_R^{-1}$.  

\begin{lemma}
\label{l.theta}
For  every  element  $g\in\Omega$,  the $A$-orbit  contains  a  unique
element of the form 
$$R_\theta=\left(\begin{array}{cc}
  \cos\theta & \sin\theta \\
 -\sin\theta   &   \cos\theta  \end{array}\right)\;\;\;\;\mbox{   with
  }\;\;\;\;\theta\in [0,2\pi[$$ 
When $g$  ranges over $\Omega$, the angle  $\theta$ varies continously
with            $g$,             and            ranges            over
  $]0,\pi/2[\cup]\pi/2,\pi[\cup]\pi,3\pi/2[\cup]3\pi/2,2\pi[$.
\end{lemma}

\begin{proof}
Consider an element $g$ in $AdS_3\simeq SL(2,\RR)$ and write 
$$g= \left(\begin{array}{cc}
                 a & b \\
                 c  & d  \end{array}\right) \;\;\;  \mbox{with} \;\;\;
                 ad-bc = 1$$ 
Then, the elements of the $A$-orbit of $g$ are the matrices 
$$g^{t}gg^{-s} = \left(\begin{array}{cc}
                 ae^{t-s}    & be^{t+s} \\
                 ce^{-(t+s)} & de^{s-t} \end{array}\right)$$
where $s$ and $t$ range over $\RR$.
Thus, the  $A$-orbit of  $g$ is  spacelike if and  only if,  for every
                 $p,q\in\RR$, the determinant of:
$$\left(\begin{array}{cc}
                 (p-q)a    & (p+q)b \\
                 -(p+q)c & (q-p)d \end{array}\right)$$           
is negative, \emph{i.e.} if and only if the quadratic form $(p-q)^{2}ad -
(p+q)^2 bc$ is positive definite.  Since $ad-bc = 1$, it follows that the
$A$-orbit of $g$ is spacelike if and only if:
$$\begin{array}{rcl}
0 < & ad & < 1 \\
-1 < & bc & < 0
\end{array}$$ 
In particular, if the $A$-orbit of $g$ is spacelike, then $abcd\neq 0$.
It follows that, if the $A$-orbit of $g$ is spacelike, then it contains an element of the form 
$$R_{\theta} = \left(\begin{array}{cc}
  \cos\theta & \sin\theta \\
 -\sin\theta & \cos\theta \end{array}\right)$$
(take $s,t$ such that $e^{2(t-s)} = d/a$ and $e^{2(s+t)}=c/b$). 
The angle $\theta$ is obviously unique, it is not a multiple of
$\frac{\pi}{2}$ (since $d\neq 0$ and $c\neq 0$), it varies
continuously with  $g$, and it takes  any value in  $[0,2\pi[$ that is
  not a multiple of $\frac{\pi}{2}$ when $g$ ranges over $\Omega$.
\end{proof}

\begin{rema} 
If $g=\left(\begin{array}{cc}
a & b \\
c & d \end{array}\right)\in\Omega$,  then the unique number $\theta\in
[0,2\pi[$ such that the rotation $R_\theta$ is in the $A$-orbit of $g$
is   characterized    by   the   equalities    $\cos^2\theta=ad$   and
$-\sin^2\theta=bc$ (see the proof of Lemma~\ref{l.theta}).
\end{rema}

Lemma~\ref{l.theta} implies that $\Omega$ has four connected
components    (corresponding    to   $\theta\in    ]0,\frac{\pi}{2}[$,
$\theta\in]\frac{\pi}{2},\pi[$,  $\theta\in ]\pi,\frac{3\pi}{2}[$, and
$\theta\in]\frac{3\pi}{2},2\pi[$). 

\begin{rema}
The four connected components of $\Omega$ are all isometric one to the
other by isometries  centralizing the group $A$.  Hence,  with no loss
of generality, we  may restrict ourselves to Torus Universes that are
obtained as quotients  of  the connected  component corresponding  to
$0<\theta<\pi/2$.  
\end{rema}

\begin{figure}[ht]
\centerline{\input{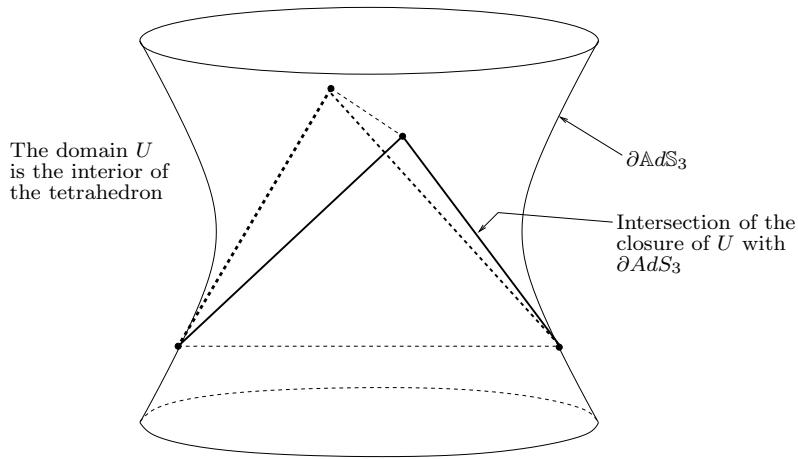}}
\caption{\label{f.torus}The domain $U$ represented in the projective model of $\AdS_3$
(more precisely, here  we use a projective chart  mapping some domain
of $\AdS_3$ in $\RR^3$).} 
\end{figure}

\begin{proof}[Proof of Theorem~\ref{p.torus-universe-GH}]
Denote  by $U$  the connected  component of  $\Omega$  corresponding to
$0<\theta<\frac{\pi}{2}$. Consider a lattice $\Gamma$ in $A$, and
consider    the    associated    Torus   Universe    $M=\Gamma\backslash
U$.  Lemma~\ref{l.theta}  provides   us  with  a  continuous  function
$\theta: U \mapsto ]0, \frac{\pi}{2}[$. By construction, this function
is increasing  with time and  $\Gamma$-invariant: it follows  that the
quotient  manifold  $M=\Gamma\backslash U$  is  equipped  with a  time
function $\bar{\theta}$.  

The    equalities   $\cos^2\theta=ad$   and    $-\sin^2\theta=bc$   (see
Lemma~\ref{l.theta}) imply that the connected component $U$ is exactly
$$\left\{g=\left(\begin{array}{cc}  a  &  b  \\   c  &  d
\end{array}\right)\in SL(2,\RR) \mbox{ such  that } 0<a \;,\;0<b \;,\;
0>c \mbox{ and }
0<d\right\}$$ 
Thus, in the Klein model  of $\AdS_3$, the connected component $U$
is the interior of a simplex which is the convex  hull of four points
in $\partial  {\mathbb A}d{\mathbb  S}_{3}$ (these points  are nothing
but the fixed points of $A$) (see figure~\ref{f.torus}). 
The main information we extract from this observation is that
$U$  is   a  convex  domain  in  ${\mathbb   A}d{\mathbb  S}_{3}$,  in
particular,  its  intersection  with  any geodesic  -  in  particular,
nonspacelike geodesics - is connected. Moreover, geodesics joining two
points of  $\partial U$ satisfying both  $bc = 0$  (respectively $ad =
0$)  are spacelike.  Hence, nonspacelike  segments in  $U$  admits two
extremities in  $\partial U$, one satisfying  $bc = 0$,  and the other
$ad = 0$. The equalities $ad = \cos^{2}\theta$, 
$bc = -\sin^{2}\theta$ imply that $\theta$ restricted to
such a nonspacelike segment takes all values between $0$, and
$\frac{\pi}{2}$. In other words, 
every nonspacelike geodesic in $U$ intersects every
fiber of $\theta$. Hence, every nonspacelike geodesic 
in $M$ intersects every fiber of $\bar{\theta}$: these
fibers are thus Cauchy surfaces, and $M$ is globally hyperbolic.

Since every fiber of $\bar{\theta}$ is a $A$-orbit, it
obviously admits constant mean curvature $\kappa(\bar{\theta})$.
Let us calculate this mean curvature at $R_{\theta}$.
We will need to take covariant derivatives, and here, the
situation is similar to the familiar situation concerning
Riemannian embeddings: if
$X$, $Y$ are vector fields in $M(2, {\bf R})$ both tangent
to $G$, then the covariant derivative $\overline{\nabla}_{X}Y$ in
$G$ is the orthogonal projection on the tangent space to
$G$ of the natural affine covariant derivative $\nabla_{X}Y$
for the affine connection on the ambient linear space.

A straightforward calculation shows that the curve 
$\theta \mapsto R_{\theta}$ is orthogonal to the
$A$-orbits, hence, the unit normal
vector to $AR_{\theta}$ at $R_\theta$ is:
$$n(\theta) = 
\left(\begin{array}{cc}
-\sin\theta & \cos\theta \\
-\cos\theta & -\sin\theta
\end{array}\right)$$
Moreover, this unit normal vector is future oriented
if we consider the orientation of $U$ for which $\theta$
increases with time. Now, for any $p$, $q$, consider the
curve $t \mapsto c(t) = g^{pt}n(\theta)g^{-qt}$. 
Its tangent vector at $t=0$ is:
$$ X_{p,q}=\left(\begin{array}{cc}
(p-q)\cos\theta & (q+p)\sin\theta \\
(q+p)\sin\theta & (q-p)\cos\theta
\end{array}\right)$$
The unit normal vector $n(t)$
to the $A$-orbit at $c(t)=g^{pt}R_{\theta}g^{-qt}$ is 
$$g^{pt}n(\theta)g^{-qt}=  
\left(\begin{array}{cc}
-e^{t(p-q)}\sin\theta & e^{t(q+p)}\cos\theta \\
-e^{-t(q+p)}\cos\theta & -e^{t(q-p)}\sin\theta
\end{array}\right)$$
Hence, the derivative at $t=0$ is:
$$\left(\begin{array}{cc}
(q-p)\sin\theta & (q+p)\cos\theta \\
(q+p)\cos\theta & (p-q)\sin\theta
\end{array}\right)$$
The orthogonal projection of this tangent vector to $AR_{\theta}$
at $R_{\theta}$ is the covariant  derivative of the unit normal vector
along the curve $t \mapsto c(t)$. It follows that
the second fundamental form is:
$$    II(X_{p,q},    X_{p,q})     =    -    \langle    X_{p,q}    \mid
\overline{\nabla}_{X_{p,q}}n(t)\rangle  =  ((p-q)^{2}  -
(p+q)^{2})\sin(2\theta)$$ 
Whereas the first fundamental form, i.e., the metric itself, is:
$$\langle X_{p,q} \mid X_{p,q} \rangle = 
(p-q)^{2}\cos^{2}\theta + (p+q)^{2}\sin^{2}\theta$$
Therefore,  the principal  eigenvalues are  $-2\mbox{cotan}\theta$ and
$2\tan\theta$.   It  follows   that  the   mean  curvature   value  is
$\kappa(\theta)  =  -4\mbox{cotan}(2\theta)$.   The  function  $\kappa
\circ \bar{\theta}$ is then increasing with time: this is the required
CMC time function. 
\end{proof} 

\begin{rema}
\label{r.not-spacelike}
The closure of the domain $U$ meets the conformal boundary at infinity
$\partial  {\mathbb A}d{\mathbb S}_{3}$  on a  topological nontimelike
circle,  but   it  is   not  a  spacelike   curve.  Actually, the intersection 
of the  closure of  $U$ with $\partial\AdS_3$ is the union of four lightlike geodesic
segments (see figure~\ref{f.torus}).    
\end{rema}

\begin{rema}
\label{r.idem-Carlip}
The Torus Universes  as defined above are the  same as those described
in~\cite{Car} in  the case of negative  cosmological constant (this
follows        immediatly       from       the        results       of
subsection~\ref{ss.proof-genus-1} below). 
Observe that the expression of the metric on the $A$-orbit  enables to recover easily the features discussed in \cite{Car}: 
the volume of the slices $\bar{\theta} = Cte$
are proportionnal to $\sin2\theta$, and the conformal classes
of these toroidal metrics describe geodesics in the modular
space $\mbox{Mod}(T)$ of the torus. 
More precisely: on the slice $\bar{\theta} = Cte$, the conformal
class and the second differential form define
naturally a point in the cotangent bundle of
$\mbox{Mod}(T)$, and when the Cte is evolving,
these data describe an orbit of the geodesic flow
on $T^{\ast}\mbox{Mod}(T)$. Conversely, every
orbit of the geodesic flow on $T^{\ast}\mbox{Mod}(T)$
corresponds to a Torus Universe. 
\end{rema}

\subsection{Every  maximal   globally  hyperbolic  spacetime,  locally
  modelled on $\AdS_3$, with closed  Cauchy surfaces of genus~$1$ is a
  Torus Universe} 
\label{ss.proof-genus-1}

In this section, we  consider a maximal globally hyperbolic Lorentzian
manifold $M$, locally modelled  on $AdS_3$, whose Cauchy surfaces are
$2$-tori. We  will prove that such  a spacetime $M$ is  isometric to a
Torus           Universe            (as           defined           in
subsection~\ref{ss.torus-universe}).           Together           with
Theorem~\ref{p.torus-universe-GH}, this  will imply that  $M$ admits a
CMC time  function.

As  in  section~\ref{s.higher-genus},  we  consider a  Cauchy  surface
$\Sigma_0$ in 
$M$, and the lift  $\widetilde\Sigma_0$ of $\Sigma_0$ in the universal
covering  $\widetilde  M$  of   $M$.   We  have  a  locally  isometric
developping  map $\cD:\widetilde  M\rightarrow AdS_3$, and  a holonomy
representation     $\rho$     of     $\pi_1(M)=\pi_1(\Sigma_0)$     in
the isometry group of $AdS_3$. We denote $\Gamma=\rho(\pi_1(M))\subset
SL(2,\RR)\times SL(2,\RR)$ (here, we prefer to see the isometry group
of $\AdS_3$ as $SL(2,\RR)\times  SL(2,\RR)$ rather than $O(2,2)$), and
we     denote      $S_0=\cD(\widetilde\Sigma_0)$.     According to
Proposition~\ref{p.embedded}, $S_0$ is properly embedded in $AdS_3$.  

The surface  $\sigma_0$ is a two-torus~: hence,  the fundamental group
of $\Sigma_0$ is isomorphic to ${\mathbb Z}^{2}$.  Moreover, according
to Proposition~\ref{p.embedded},  $\Gamma=\rho(\pi_1(M)$ is a discrete
subgroup of $SL(2, {\mathbb R}) \times SL(2,{\mathbb R})$.  
Hence, $\Gamma$ is  a lattice in some abelian group  $A = H_{L} \times
  H_{R}$, where $H_{L} = \{ e^{t h_{L}} \}_{t\in\RR}$ 
(resp. $H_{R} = \{ e^{s h_{R}} \}_{s\in\RR}$ is
a one parameter subgroup of $SL(2, {\mathbb R}) \times \{ id \}$
(resp. $\{ id \} \times SL(2, {\mathbb R})$).
Since $A$ is isomorphic to ${\mathbb R}^{2}$, these one-parameter 
groups are either parabolic or hyperbolic. In other words,
up to factor switching and conjugacy, there are only three cases
to consider:\\
- \emph{Hyperbolic - hyperbolic}: 
$$ h_{L} = h_{R} = \left(\begin{array}{cc}
                   1 & 0 \\
                   0 & -1
                   \end{array}\right)$$
- \emph{Parabolic - parabolic}: 

$$h_{L} = h_{R} = \left(\begin{array}{cc}
                   0 & 1 \\
                   0 & 0
                   \end{array}\right)$$
- \emph{Hyperbolic - parabolic}: 
$$h_{L} = \left(\begin{array}{cc}
                   1 & 0 \\
                   0 & -1
                   \end{array}\right) \;\;\mbox{ and }\;\; 
h_{R} = \left(\begin{array}{cc}
                                                      0 & 1 \\
                                                      0 & 0
                                                      \end{array}\right)$$

Let us consider an orbit $O$ of $A$. The restriction to $O$ of the ambient
Lorentzian metric defines a field of quadratic forms on $O$.
Since $A$ is a group of isometries, the quadratic forms appearing in this
field have a well-defined
type:  each of them  is either  positive definite,  negative definite,
Lorentzian, or degenerate. We call such a field of quadratic forms a 
generalized metric.  The  following  lemma describes  all  the
``isometry''  type of  generalized metrics  which can  arise by
this construction: 

\begin{lemma}
\label{isotype}
Every orbit $O$ of $A$ has dimension $1$ or $2$. Moreover:

- If $O$ has dimension $1$, then it is isometric to a line, or to an
isotropic line (i.e. equipped with the trivial null generalized metric). 

- If $O$ has dimension $2$, then it is isometric to the Euclidean  plane, 
the Minkowski plane, or the degenerate plane, i.e. the plane
with coordinates $(x,y)$ equipped with the quadratic form $dx^{2}$.
\end{lemma}

\begin{proof}
If an element $(e^{th_{L}}, e^{-sh_{R}})$ fixes a point $g$ in 
$SL(2, {\mathbb R})$, then $e^{th_{L}} = g e^{sh_{R}} g^{-1}$. Observe
that in the hyperbolic-parabolic case, this implies $s=t=0$: in this case,
every orbit of $A$ is a $2$-dimensional plane.
In  the hyperbolic-hyperbolic  case or  the  parabolic-parabolic case,
this implies
$s=t$ and $g = e^{th_{L}}$: hence, there is no $0$-dimensional orbit,
$1$-dimensional orbits are lines, and $2$-dimensional orbits are planes.

We parametrize the $A$-orbit $O$
of an element $g_{0}$ of $AdS_{3}
\approx SL(2,{\bf R})$ by 
$(s,t) \mapsto e^{t h_{L}}g_{0} e^{-s h_{R}}$.
The differential of this parametrization is: 
$$(h_{L}e^{t h_{L}}g_{0}e^{-s h_{R}})ds - (e^{t h_{L}}g_{0}e^{-s h_{R}}h_{R})dt$$
Since $h_{R}$ and $h_{L}$ commute respectively with their exponential,
and since these exponentials have determinant $1$, the determinant
of this expression reduces to the determinant of:
$$(h_{L}g_{0})ds - (g_{0}h_{R})dt$$
The quadratic form induced on the tangent space of $O$ 
at $(s,t)$ is $-\mbox{det}$ of this expression.

If $O$ has dimension $1$, then $g_{0}h_{R}g^{-1}_{0} = h_{L} = h_{R}$,
thus this determinant is equal to the determinant of $h_{L}ds - h_{L}dt$.
In the parabolic-parabolic case, we obtain identically $0$: $O$ is an
isotropic line. In the hyperbolic-hyperbolic case, we obtain $(d(s-t))^{2}$:
$O$ is a Euclidean  line.

When $O$ has dimension $2$, it is diffeomorphic to the plane.
Observe  that in the  expression above,  $s$ and  $t$ appear  only by
their differentials: 
this means that the generalized metric is actually a parallel field
of quadratic forms. In other words, it is given by the quadratic form
$-\mbox{det}(h_{L}g_{0}ds - g_{0}h_{R}dt)$ on the $2$-plane $O$ with
linear coordinates $(s,t)$.
The lemma follows from the classification of quadratic forms on the plane (the
negative definite case and the case $-(dx)^{2}$ are excluded since the
quadratic form is obtained by the restriction of a Lorentzian quadratic
form). 
\end{proof}

\begin{lemma}
\label{kice}
The surface $S_{0}$ intersects only $2$-dimensional spacelike orbits of $A$.
\end{lemma}

\begin{proof}
Let $O$ be the $A$-orbit of an element $x_{0}$ of $S_{0}$. Assume first
that $O$ has dimension $1$: according to Lemma \ref{isotype}, $O$ is
a line. Observe that $O$ is preserved by the action of $\Gamma$.
Since $\Gamma$ acts freely on $S_{0}$, $x_{0}$ is not fixed by any
element of $\Gamma$. Hence, every $\Gamma$-orbit in $O$ is dense.
It follows that there are $\Gamma$-iterates of $x_0$ arbitrarly
close to $x_0$. This is impossible, since $\Gamma$ acts properly in
a neighbourhood of $S_{0}$. 

Therefore, $O$ has dimension $2$. Assume that $O$ is not spacelike.
According to Lemma \ref{isotype},
it is isometric to the Minkowski plane or the degenerate plane.
Since $S_{0}$ is spacelike, $S_{0}$ and $O$ are transverse.
Their intersection
is a closed $1$-manifold $L$. Moreover,
the ambient Lorentzian metric restricts as a metric
on $L$ which is complete. The argument used in Proposition~\ref{p.embedded}
can then be applied once more: if $O$ is a Minkowski plane, $L$ intersects
every timelike line in $O$ in one and only one point, and if $O$
is degenerate, the same argument proves that $L$ must intersect 
every degenerate line $y=Cte$ in one and only one point (in this
situation, the projection of $L$ on the coordinate $x$ is
an isometry!). 

It follows that in both cases, $L$ is connected.
Therefore, it is isometric to $\mathbb R$. But since $O$ and
$S_{0}$ are both preserved by $\Gamma$, the same is true
for $L$: we obtain that $L \approx \mathbb R$ admits a free
and properly discontinuous isometric action by $\Gamma \approx {\mathbb Z}^{2}$.
Contradiction. 
\end{proof}

According to the lemma, some orbits of $A$ are spacelike, and
this excludes all the cases except the hyperbolic-hyperbolic case. 
Hence, $A$ is precisely the abelian group of isometries studied
in subsection \ref{ss.torus-universe} for  the definition of the Torus
Universes.  
Moreover, Lemma \ref{kice} states precisely 
that $S_{0}$ is contained in a connected component $U$ of the
domain $\Omega$. Since this is true for any Cauchy surface $\Sigma$,
and since $M$ is globally hyperbolic, the image of the developping
map is contained in $U$. Hence, $M$ embedds isometrically in the
Torus Universe $\Gamma\backslash U$. Since $M$ is maximal as a globally
hyperbolic spacetime, $M$ is actually isometric to this quotient.

Thus, we have proved:

\begin{theo}
\label{th.isom-torus}
Every  maximal   globally  hyperbolic  Lorentzian   manifold,  locally
modelled  on  $AdS_3$,  with   closed  oriented  Cauchy  surfaces  of
genus~$1$ is isometric to a Torus Universe.  
\end{theo}

\begin{coro} 
\label{c.torus-maximal}
Torus Universes are maximal as globally hyperbolic spacetimes. 
\end{coro}

\begin{proof}[Proof of Theorem~\ref{t.CMC-foliation} in the case $g=1$]
The      result      follows     from      Theorem~\ref{th.isom-torus}
    and~\ref{p.torus-universe-GH}.  
\end{proof}

\section*{Acknowledgements}
This work has been partially supported by the CNRS and the 
ACI "Structures g\'eom\'etriques et trous noirs".

\end{document}